\documentclass[a4paper,11pt]{amsart}

\usepackage{graphicx,amsmath,amsfonts,amsthm,amssymb,url,dsfont,bbm}

\usepackage[OT2,T1]{fontenc}
\DeclareSymbolFont{cyrletters}{OT2}{wncyr}{m}{n}
\DeclareMathSymbol{\Sha}{\mathalpha}{cyrletters}{"58}

\input xy
\xyoption{all}

\def\]{\textup{\mbox{]\hspace{-.15em}]}}}
\def\[{\textup{\mbox{[\hspace{-.15em}[}}}

\def\got{\mathfrak}

\newenvironment{pf}
{\medskip\noindent {\it Preuve --- \ }}
{\hfill\nobreak $\Box$ \par\bigbreak}

\newcommand{\F}{{\mathbb F}}

\newcommand{\isomo}{\overset{\sim}{\rightarrow}}
\newcommand{\cal}{\mathcal}

\newcommand{\WW}{\mathcal{W}}

\newcommand{\Hom}{\mathrm{Hom}}

\newcommand{\Gal}{\mathrm{Gal}}
\newcommand{\G}{\mathrm{G}}

\newcommand{\pn}{\par \noindent}

\newcommand{\Qpb}{\overline{\mathbb{Q}}_p}

\newcommand{\GL}{\mathrm{GL}}

\newcommand{\ps}{\par \smallskip}

\newcommand{\End}{\mathrm{End}}
\newcommand{\N}{\mathbb{N}}
\newcommand{\Z}{\mathbb{Z}}

\newcommand{\Fpb}{\overline{\mathbb{F}}_p}
\newcommand{\Q}{\mathbb{Q}}

\newcommand{\Qp}{\mathbb{Q}_p}

\newcommand{\rhob}{{\bar\rho}}

\newcommand{\Ker}{\mathrm{Ker}}

\newcommand{\AAA}{\mathbb{A}}
\newcommand{\OO}{\mathcal{O}}

\newcommand{\cA}{{\mathcal A}}

\newcommand{\cB}{{\mathcal B}}

\newcommand{\Ref}{\cal{F}}
\newcommand{\fg}{(\varphi,\Gamma)}
\newcommand{\Ro}{\mathcal{R}}
\newcommand{\uk}{{\bf k}}

\newcommand{\cDc}{\mathcal{D}_{\rm cris}}

\newcommand{\FF}{{\mathbb{F}}}

\newtheorem{thm}[subsection]{Th\'eor\`eme}

\newtheorem{lemma}[subsection]{Lemme}
\newtheorem{remark}[subsection]{Remarque}
\newtheorem{cor}[subsection]{Corollaire}
\newtheorem{prop}[subsection]{Proposition}

\newtheorem{conj}[subsection]{Conjecture}

\newtheorem{definition}[subsection]{D\'efinition}

\newcommand{\Dc}{D_{\rm cris}}
\newcommand{\cE}{{\mathcal{E}}}
\newcommand{\cF}{{\mathcal{F}}}
\newcommand{\cX}{{\mathfrak X}}

\newcommand{\Drig}{D_{\rm rig}}

\newcommand{\Fil}{{\rm Fil}}
\newcommand{\Aff}{{\rm Aff}}

\title{Sur la densit\'e des representations cristallines de ${\rm
Gal}(\Qpb/\Qp)$}

\begin{document}

\author[G. Chenevier]{Ga\"etan Chenevier}
\address{Ga\"etan Chenevier\\ C.N.R.S, Centre de Math\'ematiques Laurent Schwartz, \'Ecole Polytechnique, 91128 Palaiseau Cedex\\ France.  chenevier@math.polytechnique.fr }


\begin{abstract} Let $\cX_d$ be the $p$-adic analytic space classifying the
semisimple continuous representations ${\rm Gal}(\Qpb/\Q_p) \rightarrow
\GL_d(\Qpb)$. We show that the crystalline representations are Zarski-dense in many
irreducible components of $\cX_d$, including the components made of residually
irreducible representations.  This extends to any dimension $d$ previous
results of Colmez and Kisin for $d = 2$.

For this we construct an analogue of the infinite fern of Gouv\^ea-Mazur in
this context, based on a study of analytic families of trianguline
$\fg$-modules over the Robba ring.  We show in particular the
existence of a universal family of (framed, regular) trianguline
$\fg$-modules, as well as the density of the crystalline
$\fg$-modules in this family.  These results may be viewed as a local
analogue of the theory of $p$-adic families of finite slope automorphic
forms, they are new already in dimension $2$.
The technical heart of the paper is a collection of results about the
Fontaine-Herr cohomology of families of trianguline $\fg$-modules.

\end{abstract}

\maketitle 
\section*{Introduction}

Soient\footnote{\date \, L'auteur est financ\'e par le C.N.R.S.} $\FF_q$
un corps fini de caract\'eristique $p>0$, $W$ l'anneau des vecteurs de Witt
de $\FF_q$, et $F=W[1/p]$. Si $K$ est une extension finie de $\Q_p$ on
d\'esigne par $G_K={\rm Gal}(\overline{K}/K)$ son groupe de Galois absolu. Fixons $$\rhob :
G_{\Q_p} \rightarrow
\GL_d(\FF_q)$$ une repr\'esentation continue absolument
irr\'eductible\footnote{Cela signifie que $\rhob$ est irr\'eductible et le
reste apr\`es extension des scalaires \`a $\Fpb$. Rappelons aussi que $\rhob$ est n\'ecessairement d\'efinie sur le sous-corps de $\FF_q$ engendr\'e par les coefficients des $\det(t-\rhob(g))$, $g \in
G_{\Q_p}$ (l'obstruction de Schur est vide pour les corps finis), 
de sorte qu'il n'est pas restrictif de supposer que ce corps est exactement 
$\FF_q$, ce que nous faisons d\'esormais par commodit\'e.} et
d\'esignons par $R$ la $W$-alg\`ebre de d\'eformation universelle de $\rhob$ au sens de
Mazur~\cite{mazurdef}. Supposons enfin que $\rhob \not\simeq \rhob(1)$ (condition
automatiquement satisfaite si $p-1$ ne divise pas
$d$). D'apr\`es Tate~\cite{tate}\footnote{En effet, cela entra\^ine que $H^2(G_{\Q_p},{\rm ad}(\rhob))=0$
et $\dim_{\FF_q} H^1(G_{\Q_p},{\rm ad}(\rhob)) = d^2+1$}, cela
entra\^ine que $R \simeq W[[T_0,\dots,T_{d^2}]]$. 
En particulier, l'espace rigide analytique
$X$ associ\'e par Berthelot \`a $R[1/p]$ est
la boule unit\'e ouverte $\{ (T_0,\dots,T_d), |T_i|<1\}$ de dimension $d^2+1$ sur $F$. \ps

Rappelons que si $L$ est une extension finie de $F$, $X(L)$ param\`etre les classes de $L$-isomorphie de repr\'esentations continues
$G_{\Q_p} \rightarrow
\GL_d(\OO_L)$ dont la r\'eduction modulo $\pi_L$ est isomorphe \`a $\rhob \otimes_{\F_q} k_L$. 
Bien s\^ur, $\OO_L$ d\'esigne ici l'anneau des entiers de $L$, $\pi_L$ en est une uniformisante, et $k_L=\OO_L/\pi_L\OO_L$. Enfin, on dira que $x \in X(L)$ a une certaine propri\'et\'e  si la repr\'esentation associ\'ee $\rho_x :
G_{\Q_p} \rightarrow \GL_d(L)$ a cette propri\'et\'e. Le r\'esultat principal de cet article est le suivant. \ps\ps

{\bf Th\'eor\`eme A :} {\it Il existe une extension finie $L$ de $F$ telle que 
les points cristallins de $X(L)$ sont Zariski-denses et d'accumulation dans
$X$}.\ps\ps

Autrement dit $X(L)$ contient au moins un point cristallin et pour tout $x
\in X(L)$ cristallin et tout ouvert affino\"ide connexe $U \subset X$ contenant $x$ alors
l'ensemble des points cristallins de $U(L)$ est Zariski-dense dans $U$ : 
toute fonction rigide analytique sur $U$ s'annulant sur les points cristallins
de $U(L)$ est identiquement nulle. En fait, on verra que tout $L$ contenant
une extension explicite de $F$ (de degr\'e $\leq d^2$) convient.\ps\ps

Ce r\'esultat est facile quand $n=1$ et d\'emontr\'e ind\'ependamment par Colmez et Kisin quand $n=2$ dans \cite[\S 5]{colmeztri} et \cite[\S 1]{kisinfern}. Leurs preuves sont inspir\'ees de la foug\`ere infinie de Gouv\^ea et Mazur \cite{gm}, qui concerne un analogue global du probl\`eme qui nous int\'eresse, et qui est elle-m\^eme issue des travaux fondateurs de Coleman~\cite{coleman} sur les familles $p$-adiques de formes modulaires. Bien que les techniques employ\'ees par Colmez et Kisin soient diff\'erentes, leurs approches sont tr\`es similaires ; nous nous bornerons ci-dessous \`a d\'ecrire celle de Colmez car c'est son point de vue que nous \'etendrons par la suite. \ps

Quand $n=2$ l'espace $X$ est une boule de dimension $5$. Consid\'erons 
l'ensemble $S$ des paires $(x,t)$ o\`u $x \in X$ est 
{\it trianguline}\footnote{On rappelle, suivant 
Colmez~\cite{colmeztri}, qu'une repr\'esentation 
$\rho : G_{\Q_p} \rightarrow \GL_d(L)$ est dite 
trianguline si son $\fg$-module pris sur l'anneau 
de Robba est extension successive $\fg$-modules de 
rang $1$. Nous appellerons le choix d'une telle extension une triangulation de $\rho$.} 
et o\`u $t$ est la donn\'ee d'une triangulation de $x$. D'apr\`es Colmez, 
$S$ admet une structure "naturelle" d'espace analytique $p$-adique, qui est lisse et de dimension $4$. 
Notons $$\pi : S \rightarrow X$$
l'application ensembliste $(x,t) \mapsto x$. Son image $\pi(S) \subset X$ est appel\'ee la {\it foug\`ere infinie}. Colmez d\'emontre que : \begin{itemize}
\item[(a)] $\pi(S)$ contient tous les points cristallins de $X$. De plus, chaque $x \in X$ cristallin non exceptionnel\footnote{Pr\'ecis\'ement, si le Frobenius cristallin de $D_{\rm cris}(\rho_x)$ est semi-simple.} admet exactement deux ant\'ec\'edents dans $S$. \ps
\item[(b)] Les $(x,t) \in S$ avec $x$ cristallin non exceptionnel sont Zariski-denses et d'accumulation dans $S$.\ps
\item[(c)] $\pi$ est localement analytique : pour (presque\footnote{Le 
presque ici fait r\'ef\'erence au fait que pour Colmez le param\`etre 
de $(s,t)$ doit \^etre $p$-r\'egulier, notion qui sera introduite plus tard
(cette hypoth\`ese n'est d'ailleurs plus n\'ecessaire gr\^ace au th\'eor\`eme 
\ref{thmdegal2}).}) tout $s \in S$ il existe un voisinage affino\"ide de $s$ dans $S$ (une boule ferm\'ee de dimension $4$) surlequel $\pi$ est une immersion analytique. \ps
\item[(d)] si $x \in X$ est cristallin non exceptionnel, alors les deux sous-boules\footnote{Ces deux familles sont les analogues locaux des deux familles de Coleman passant par les deux formes jumelles associ\'ees \`a une forme modulaire propre et $p$-ancienne.} de $X$ passant par $x$ qui sont donn\'ees par (a) et (c) ne sont pas confondues.
\end{itemize}

Comme dans chacune des deux sous-boules mentionn\'ees dans le (d), les points cristallins non exceptionnels sont Zariski-denses par le (b), on obtient une id\'ee de la structure fractale de la foug\`ere infinie dans $X$ (et une justification pour son nom!). Le th\'eor\`eme ci-dessus s'en d\'eduit ais\'ement. Dans 
l'approche de Kisin, disons simplement que l'espace $S$ est remplac\'e dans cet argument par sa version $X_{fs}$ 
construite dans~\cite{kisin}. La diff\'erence essentielle est que Kisin d\'ecoupe $X_{fs}$ dans l'espace $X \times \mathbb{G}_m$ (il doit donc le "majorer", ce qui est assez d\'elicat) 
alors que Colmez le construit explicitement, ce qui est peut-\^etre plus naturel par rapport \`a notre probl\`eme. Quand $n=2$, mentionnons que ce r\'esultat de densit\'e des cristallines joue un r\^ole important dans la preuve par Colmez de la correspondance de Langlands locale $p$-adique \cite{colmezgros}. Ajoutons enfin que toujours quand $n=2$, Nakamura a \'etendu dans \cite{nakamura} l'approche de Kisin et le th\'eor\`eme ci-dessus au cas o\`u $G_{\Q_p}$ est remplac\'e par $G_K$ avec $K/\Q_p$ finie arbitraire. \ps

Notre d\'emonstration reprend pour $d$ quelconque les \'etapes de la preuve
esquiss\'ee ci-dessus.  L'ensemble $S$, ainsi que la foug\`ere infinie, sont
d\'efinis de la m\^eme mani\`ere.  Son \'etude a \'et\'e amorc\'ee par
Bella\"iche et l'auteur dans \cite[\S 2]{bch}.  Une cons\'equence simple des
travaux de Berger est que les
triangulations d'une repr\'esentation cristalline $V$ de dimension $d$
quelconque sont en bijection avec les drapeaux $\varphi$-stables de
$D_{\rm cris}(V)$, de sorte que g\'en\'eriquement (au sens na\"if) il y en a
$d!$.  C'est la g\'en\'eralisation du (a) dont nous aurons besoin ici.  De
plus, mentionnons que bien que l'on n'y munisse pas $S$ d'une structure
analytique, une \'etude du voisinage infinit\'esimal de chaque point est
aussi men\'ee {\it loc.  cit.}, ce qui est un premier pas en direction (c),
tout en \'etant insuffisant\footnote{Mentionnons qu'il semble
difficile d'appliquer le th\'eor\`eme d'approximation d'Artin pour relever
les germes formels de familles de $\fg$-modules triangulins construits
dans~\cite[\S 2]{bch} en des vrais germes de
familles analytiques.  Cela vient entre autres de ce que nous ne savons pas d\'emontrer
que les foncteurs sous-jacents sont de pr\'esentation finie, question qui nous semble assez profonde. } pour le th\'eor\`eme $A$. Les
g\'en\'eralisations de (b) et (c) \`a la dimension $d$ formeront le coeur
technique et nouveau de cet article, et nous y reviendrons plus loin en d\'etail. Un
point assez d\'elicat concerne la g\'en\'eralisation du (d) : en un point
cristallin $x \in X$ assez g\'en\'eral, la foug\`ere infinie admet $d!$
feuilles chacune \'etant de dimension $\frac{d(d+1)}{2}+1$, et le probl\`eme
est de comprendre leurs positions relatives.  Ce probl\`eme est facile \`a
r\'esoudre pour des raisons de dimension quand $d=2$, car une analyse ad-hoc
des param\`etres triangulins sur chacune des deux feuilles de dimension $4$
montre que les deux feuilles ne sont pas confondues. Le cas g\'en\'eral est
nettement plus d\'elicat, mais a d\'ej\`a
\'et\'e r\'esolu par l'auteur dans~\cite{chU3} (en fait dans~\cite{chpeccot}). Le r\'esultat cl\'e est que 
si $x$ cristallin est assez g\'en\'erique en un sens pr\'ecis rappel\'e plus bas, alors la somme des espaces 
tangents en $x$ des $d!$ feuilles de la foug\`ere en $x$ est l'espace tangent de $X$ en $x$ tout entier. 
L'argument d'adh\'erence Zariski employ\'e dans~\cite{chU3} permet alors de
conclure : nous renvoyons \`a 
la section~\ref{reductions} pour l'argument complet.\ps\ps

Revenons sur la structure analytique de l'espace $S$. Nous raisonnerons enti\`erement dans le monde des 
$\fg$-modules sur l'anneau de Robba (non n\'ec\'essairement \'etales), ce qui est essentiellement permis 
par un th\'eor\`eme de Kedlaya et Liu~\cite{kedliu}. Si $A$ est une $\Q_p$-alg\`ebre affino\"ide, on 
note $\Ro_A$ l'anneau de Robba \`a coefficients dans $A$ (voir \S~\ref{defRoA} pour la d\'efinition pr\'ecise).  
Soient $\Aff$ la cat\'egorie des $\Q_p$-alg\`ebres affino\"ides et $$F_d^\square : \Aff \longrightarrow {\rm Ens}$$
le foncteur associant \`a $A$ l'ensemble des classes d'\'equivalence de $\fg$-modules sur $\Ro_A$ qui sont 
triangulins, r\'eguliers et rigidifi\'es : nous renvoyons \`a la section~\ref{San} pour les d\'efinitions de ces
termes.\footnote{Pr\'ecisions tout de m\^eme que dans toute la suite, un $\fg$-module triangulin est la
donn\'ee d'un $\fg$-module {\it muni} d'une filtration (la filtration fait
partie de la donn\'ee).}
\ps\ps

{\bf Th\'eor\`eme B :} {\it Le foncteur $F_d^{\square}$ est repr\'esentable par un espace analytique $p$-adique $S_d^\square$ qui est 
irr\'eductible et lisse sur $\Q_p$, de dimension $\frac{d(d+3)}{2}$. Les $\fg$-modules cristallins sont Zariski-denses et d'accumulation dans $S_d^\square$.}
\ps\ps

Ce th\'eor\`eme peut
\^etre vu comme un analogue local de la th\'eorie des vari\'et\'es de Hecke (ou
"{\it eigenvarieties}").\footnote{Le fait que l'analogue "automorphe" global de ce
r\'esultat \'etait d\'ej\`a connu est aussi la raison pour laquelle nous avions
\'etudi\'e le cas global en premier dans~\cite{chU3}. Ces deux contextes comportent des similarit\'es
\'evidentes mais aussi des diff\'erences importantes, par exemple on ne dispose
malheureusement pas en global d'analogue de la
proposition~\ref{existcrisgen} (ii).}
Notons que sa premi\`ere partie est nouvelle m\^eme pour
$d=2$, o\`u elle compl\`ete des r\'esultats de Colmez et de Kisin. Toujours dans le cas
$d=2$, nous construisons en fait une famille "tautologique" des $\fg$-modules
non n\'ecessairement r\'eguliers (Th\'eor\`eme~\ref{thmdegal2}). 
Mentionnons que la notion de rigidification utilis\'ee fait que l'espace
$S_d^\square$ admet moralement $d-1$ dimensions de plus que sa d\'efinition
na\"ive, en revanche elle permet de traiter sur un pied d'\'egalit\'e tous
les $\fg$-triangulins (y compris ceux qui sont scind\'es par exemple). Le
th\'eor\`eme $B$ a des cons\'equences int\'eressantes concernant la
construction de repr\'esentations cristallines ayant certaines
propri\'et\'es. Par exemple il permet de montrer que {\it si $\rhob : G_{\Q_p} \rightarrow \GL_d(\F_q)$
est une repr\'esentation semi-simple continue quelconque, il existe une extension finie
$L/\Q_p$ et une repr\'esentation $G_{\Q_p} \rightarrow \GL_d(L)$ cristalline
absolument irr\'eductible dont la repr\'esentation r\'esiduelle est
$\rhob$} (Prop.~\ref{liftgeneral}).\ps

Le coeur technique du th\'eor\`eme $B$, et de cet article, est un ensemble de r\'esultats sur les $\fg$-modules sur
$\Ro_A$, notamment sur leur cohomologie \`a la Fontaine-Herr (\cite{herr1},\cite{herr2}). Ceci
fait l'objet de la section~\ref{cohomologie}. Une \'etape de la d\'emonstration est la v\'erification dans ce 
contexte que les complexes
$C_{\varphi,\gamma}$ et $C_{\psi,\gamma}$ sont quasi-isomorphes comme dans
la th\'eorie classique de Herr, ce qui avait notamment \'et\'e conjectur\'e
par Kedlaya dans~\cite[\S 2.6]{kedlayaseul}. Concr\`etement, il s'agit
d'\'etudier la structure de $D^{\psi=0}$ comme $\Gamma$-module quand $D$ est
un $\fg$-module sur $\Ro_A$ ou sur $\Ro_A^+$. Ce point assez technique
repose sur une \'etude pr\'eliminaire des familles de $\Gamma$-modules
effectu\'ee en section~\ref{prelimfamilles}. Notre preuve est directement
inspir\'ee de \cite[\S V]{colmezgros} qui en d\'emontre le cas particulier o\`u
$A$ est un corps. \ps

Nous calculons enfin la cohomologie de $\Ro_A(\delta)$
quand $\delta : \Q_p^\ast \rightarrow A^\ast$ est un caract\`ere
continu, \'etendant des r\'esultats de Colmez (en degr\'es $0$ et $1$,
\cite{colmeztri})) et de Liu (en degr\'e $2$, \cite{liu}) dans le cas
particulier o\`u $A$ est un corps. Notre preuve, bien qu'inspir\'ee de celle de
Colmez, est en fait un peu plus simple que celle
dans~\cite{colmeztri} : l'id\'ee nouvelle essentielle est de remplacer son d\'evissage pour se ramener
au cas o\`u $v(\delta(p)) < 0$ par un argument direct utilisant la
transform\'ee de ... Colmez ! qui est un d\'evissage de l'anneau de Robba.
Nous d\'ecrivons de plus la structure de $\Ro_A(\delta)^{\psi=1}$ 
comme module sur une certaine compl\'etion de $A[\Gamma]$ not\'ee $\Ro_A^+(\Gamma)$.
Nous \'etendons enfin ces r\'esultats \`a tous les $\fg$-modules triangulins
sur $\Ro_A$.
Mentionnons qu'une des sp\'ecificit\'es de cette th\'eorie en famille est
l'absence de dualit\'e. De plus, le coeur (au sens de Fontaine) d'une
famille de $\fg$-modules sur $\Ro_A$ n'est pas n\'ec\'essairement libre sur
$\Ro_A^+(\Gamma)$. Voici un \'echantillon des r\'esultats obtenus.\ps

{\bf Th\'eor\`eme C :} {\it Si $D$ est un $\fg$-module triangulin de rang $d$ sur $\Ro_A$,
alors $H^i(D)$ est de type fini sur $A$ pour tout $i$ et dans le groupe de
Grothendieck des $A$-modules de type fini on a la relation
$[H^0(D)]-[H^1(D)]+[H^2(D)]=-[A^d]$. De plus, la formation des $H^i(D)$ commute \`a
tout changement de base plat affino\"ide. Enfin, $D^{\psi=1}$ contient un
$\Ro_A^+(\Gamma)$-module libre de rang $d$, le quotient \'etant de type fini
sur $A$. \ps
	Si $D$ est r\'egulier, alors $H^0(D)=H^2(D)=0$ et $H^1(D)$ est libre de
rang $d$ sur $A$. La formation des $H^i(D)$ commute alors \`a tout changement
de base affino\"ide. Si $D$ est $p$-r\'egulier, alors $D^{\psi=1}$ est libre de
rang $d$ sur $\Ro_A^+(\Gamma)$.}
\ps\ps

L'auteur remercie chaleureusement Pierre Colmez pour les explications pr\'ecieuses de ses r\'esultats,
dont les \S~\ref{checkgamma} et \S~\ref{cohoRoA} sont tr\`es largement inspir\'es, ainsi
que Laurent Fargues et Olivier Ta\"ibi pour des discussions utiles.

\tableofcontents

\section{$\fg$-modules sur $\Ro_A$}\label{prelimfamilles}
\newcommand{\RR}{\mathcal{R}}

\subsection{Quelques anneaux de fonctions analytiques}\label{defRoA}

Nous renvoyons \`a \cite{bgr} pour les g\'en\'eralit\'es sur les espaces analytiques $p$-adiques au sens de Tate. Si $X$ est un tel espace, nous noterons $\OO(X)$ la $\Q_p$-alg\`ebre de ses fonctions globales. Quand $X$ est
affino\"ide, c'est une alg\`ebre de Banach noeth\'erienne, et on note aussi $X={\rm Sp}(\OO(X))$. 
En g\'en\'eral, on munit $\OO(X)$ de la topologie de la convergence uniforme sur tout ouvert 
affino\"ide, c'est une alg\`ebre de Fr\'echet si $X$ admet un recouvrement d\'enombrable admissible 
par des affino\"ides, ce qui sera le cas pour tous les espaces ci-dessous.
Si $X$ est r\'eunion admissible d'ouverts affino\"ides $X_n$ ($n\geq 0$) 
avec $X_n \subset X_{n+1}$, c'est une alg\`ebre de Fr\'echet-Stein au sens
de~\cite[\S 3]{schneiderteitelbaum}.\footnote{Une cons\'equence de
cette propri\'et\'e est l'existence d'une sous-cat\'egorie pleine,
ab\'elienne, naturelle de celle des $\OO(X)$-modules que
Schneider et Teitelbaum nomment "co-admissibles" : un $\OO(X)$-module est
dit co-admissible si il est la limite projective d'une suite de
$\OO(X_n)$-modules de type fini $M_n$ munis d'isomorphismes $M_{n+1}
\otimes_{\OO(X_{n+1})} \OO(X_n) \rightarrow M_n$ pour tout $n\geq 0$. 
Ces modules ont par d\'efinition une
topologie canonique de $\OO(X)$-module de Fr\'echet. Les $\OO(X)$-modules de pr\'esentation
finie sont admissibles, ainsi que tous leurs sous-modules de type fini et
plus g\'en\'eralement ferm\'es. Nous renvoyons \`a {\it loc. cit} pour plus
de renseignements. Nous n'aurons recours \`a ces r\'esultats que dans la preuve du lemme~\ref{topo} (iv) (lui-m\^eme utilis\'e uniquement pour le th\'eor\`eme~\ref{h1universel}).}

Soit $A$ une $\Q_p$-alg\`ebre affino\"ide. Un {\it mod\`ele} de $A$ est une 
sous-$\Z_p$-alg\`ebre $\mathcal{A} \subset A$ topologiquement de type fini, 
i.e. quotient de $\Z_p\langle t_1,\dots,t_m\rangle$ pour un certain $m\geq  
1$, et telle que $\mathcal{A}[1/p]=A$. Un mod\`ele est ouvert, born\'e, sans
$p$-torsion, et complet s\'epar\'e pour la topologie $p$-adique. Pour toute famille finie 
$x_1$,\dots, $x_r$ d'\'el\'ements de $A$ \`a puissances (positives) born\'ees, il existe toujours un mod\`ele $\mathcal{A}$ de $A$ contenant les
$x_i$. Si $A$ est r\'eduit, l'ensemble de ses \'el\'ements \`a puissances born\'ees est le plus grand mod\`ele
de $A$ (Tate). Si $|.|$ est une norme sur $A$ (sous-entendu, sous-multiplicative et pour laquelle $A$ est complet), alors sa boule unit\'e $A^0$ est un mod\`ele de $A$. \ps

Soit $I \subset [0,1[$ un intervalle d'extr\'emit\'es dans $p^\Q$, on note $B_I$ l'ouvert 
admissible de la droite affine rigide $\mathbb{A}^1$ (de param\`etre $T$) d\'efini par $|T| \in I$. 
C'est un disque si $0 \in I$ et une couronne sinon, il est affino\"ide si $I$ est un segment. On notera 
encore $T \in \OO(B_I)$ le param\`etre tautologique. Soit $A$ une $\Q_p$-alg\`ebre affino\"ide. Si 
$I \subset [0,1[$ est un intervalle on pose $$\cE_A^I:=\OO( {\rm Sp}(A) \times B_I).$$ 
\begin{itemize}\item[-]
Si $0 \in I$, $\cE_A^I \subset A[[T]]$ est aussi la sous-alg\`ebre des s\'eries 
$f=\sum_{n \in \N} a_n T^n$ telles que pour tout $r \in I$ on ait $|a_n|r^n \rightarrow 0$ 
quand $n\rightarrow \infty$. C'est une $A$-alg\`ebre de Fr\'echet pour les normes 
$|f|_{[0,r]}:=\sup_{n \in \N} |a_n|r^n$, si $|.|$ est une norme fix\'ee sur $A$. 
\item[-] Si $0 \notin I$, $\cE_A^I$ est aussi la $A$-alg\`ebre des s\'eries de 
Laurent $f=\sum_{n \in \Z} a_n T^n$ telles que pour tout $[r,s] \subset I$ on 
ait $|a_n|s^n \rightarrow 0$ et $|a_{-n}| r^{-n} \rightarrow 0$ quand 
$n\rightarrow \infty$. C'est une $A$-alg\`ebre de Fr\'echet pour les 
normes $|f|_{[r,s]}:=\sup(\sup_{n \in \N} |a_n|s^n, \sup_{n \in \N} |a_{-n}|r^{-n})$, 
o\`u $[r,s] \subset I$. \end{itemize}
Notons que les descriptions ci-dessus sont classiques quand $I$ est un segment, et 
s'en d\'eduisent en g\'en\'eral en consid\'erant le recouvrement admissible de 
${\rm Sp}(A) \times B_I$ par les ${\rm Sp}(A) \times B_J$ pour $J \subset I$ un 
segment. De plus, $\cE_A^I$ est une alg\`ebre affino\"ide si $I$ est un
segment et de Fr\'echet-Stein en g\'en\'eral (consid\'erer le recouvrement par
les ${\rm Sp}(A) \times B_{I_n}$ o\`u $I_n$ est une suite croissante de
segments recouvrant $I$). \ps

On pose encore $$\Ro_{A,r}=\cE_A^{[r,1[}, \, \, \, \Ro_A=\bigcup_{0<r<1} 
\Ro_{A,r},\, \, \,  {\rm et}\, \, \, \, \Ro^+_A=\Ro_{A,0}.$$ 
Lorsque $A=\Q_p$ on omettra souvent de le mentionner en indice. Par exemple 
$\Ro:=\Ro_{\Q_p}$ et $\cE^I:=\cE_{\Q_p}^I$. De plus, si $X={\rm Sp}(A)$ on remplacera 
parfois $A$ par $X$ dans les notations ci-dessus, de sorte que par exemple $\cE_X^I:=\cE_{\OO(X)}^I$. \ps

Il sera commode par la suite d'introduire certains mod\`eles sur $\Z_p$ des anneaux
ci-dessus.\ps

\begin{lemma}\label{lemmmodel} Soit $I=[p^{-\frac{a}{n}},p^{-\frac{b}{m}}]$ avec
$a,b,m,n$ entiers tels que $0\leq \frac{a}{n} \leq \frac{b}{m}$ (on suppose
les d\'enominateurs non nuls et le fractions r\'eduites). Alors
$\OO^I:=\Z_p\langle \frac{T^m}{p^b},\frac{p^a}{T^n},T \rangle$ est un mod\`ele de $\cE^I$.
\par
De plus, le morphisme surjectif naturel $\Z_p \langle T,U,V \rangle
/(p^bU-T^m,T^nV-p^a) \rightarrow \OO^I$ a pour noyau la $p^\infty$-torsion de $\Z_p
 \langle T,U,V \rangle
/(p^bU-T^m,T^nV-p^a)$ 
\end{lemma}

\begin{pf} En effet, par d\'efinition de $B_I$ on a $\cE^I=\Q_p\langle T, U,
V\rangle/(p^b U-T^m, T^n V-p^a)$, donc les \'el\'ements $\frac{T^m}{p^b}, \frac{p^a}{T^n}
\in \cE^I$ sont \`a puissances born\'ees et $\OO^I$ est un mod\`ele de
$\cE^I$ car image de
$\Z_p\langle T, U,     
V\rangle$. Le derni\`ere assertion suit car le morphisme de l'\'enonc\'e est
un isomorphisme apr\`es avoir invers\'e $p$.
\end{pf}


Si $I$ est un segment on d\'efinit $\OO^I$ comme dans le lemme ci-dessus. Si $\mathcal{A}$ est un mod\`ele de $A$, alors $\OO_\cA^I:=\cA \widehat{\otimes}_{\Z_p} \OO^I$ est un 
mod\`ele de $\cE_A^I$. \ps

Enfin, on pose $\cE_\cA^{\dag, 0}=:\cA[[T]]$, et si $n\geq 1$, on d\'efinit 
$\cE_{\mathcal{A}}^{\dag,n}$ comme \'etant le compl\'et\'e de $\mathcal{A}[[T]][\frac{p}{T^n}]$ pour la topologie $p$-adique. C'est aussi l'anneau des s\'eries de Laurent de la forme $\sum_{k\in \Z} a_k T^k$ telles que $a_k \in \mathcal{A}$ pour tout $k \in \Z$, $a_{-k} \in p^{\alpha_k}\mathcal{A}$ pour $k \geq 0$, o\`u $\alpha_k$ est une suite d'entiers $\geq [\frac{k}{n}]$ et tendant vers l'infini avec $k$. En particulier, $\cE_{\mathcal{A}}^{\dag,n} \subset \Ro_{A,p^{-\frac{1}{n}}}$. On le munit de sa {\it topologie faible} :  
une base de voisinages de $0$ est l'ensemble des $p^\alpha 
\cE_{\mathcal{A}}^{\dag, n} + T^\beta \mathcal{A}[[T]]$, $\alpha,\beta\geq 0$ entiers. Il est aussi complet 
pour cette topologie. Nous renvoyons \`a \cite[\S 17]{schneider} pour les g\'en\'eralit\'es sur les produits 
tensoriels compl\'et\'es. Le (i) ci-dessous fait notamment le lien entre les
d\'efinitions employ\'ees ici et celles de~\cite{kedliu}
et de~\cite{bergercolmez}.\ps

\begin{lemma}\label{topo}
\begin{itemize}
\item[(i)] Pour tout $I$, l'application naturelle $A \widehat{\otimes}_{\Q_p} \cE^I \rightarrow \cE^I_A$ est 
un isomorphisme de Fr\'echets. \ps
\item[(ii)] Pour tout $0<r<1$, la norme $|.|_{[0,r]}$ sur $\Ro_A^+$ induit
sur $\cA[[T]]$ la topologie d\'efinie par l'id\'eal    
$(p,T)\mathcal{A}[[T]]$, ou ce qui revient au m\^eme, la {\it
topologie
faible} dont une base de voisinages de $0$ est donn\'ee par les $p^\alpha
\cA[[T]] + T^\beta \cA[[T]]$ avec $\alpha,\beta \geq 0$. L'injection naturelle $\mathcal{A}[[T]] 
\rightarrow \Ro_A^+$ est d'image ferm\'ee. \ps
\item[(iii)] Pour tout $n\geq 1$ et $p^{-1/n} \leq s <1$, la norme $|.|_{[p^{-1/n},s]}$ 
induit sur $\cE_{\mathcal{A}}^{\dag, n}\subset \Ro_{A,p^{-1/n}}$ la
topologie faible. En particulier, $\cE_{\mathcal{A}}^{\dag, n}$ est ferm\'e
dans $\Ro_{A,p^{-1/n}}$. \ps
\item[(iv)] Si $J$ est un ideal de $A$ et $I$ est un intervalle quelconque
de $[0,1[$, alors l'application
naturelle $\cE_A^I/J\cE_A^I \rightarrow \cE_{A/J}^I$ est un
isomorphisme.\ps
\item[(v)] Pour tout intervalle $I$ de $[0,1[$, le $A$-module 
$\cE_A^I$ est plat. En particulier, $\Ro_A$ est plat sur $A$.
\end{itemize}
\end{lemma}

\begin{pf} Quand $I$ est un segment, $B_I$ est affino\"ide et le (i) est \'evident. En g\'en\'eral, l'injectivit\'e de $A \otimes_{\Q_p} \Q_p^{\Z} \rightarrow A^ \Z$ entra\^ine 
celle de $\iota_I : A \otimes_{\Q_p} \cE^I \rightarrow \cE^I_A$. Si $J \subset I$ est un segment, on a une semi-norme naturelle produit-tensoriel $|.|_J:=|.|\otimes |.|_J$ sur $A \otimes_{\Q_p} \cE^I$ associ\'ee :  
$|x|_J$ est l'infimum sur toutes les \'ecritures $x=\sum_i a_i \otimes f_i$ des ${\rm sup}_i |a_i||f_i|_{J}$. 
En particulier, $|\iota_I(x)|_J\leq |x|_J$ pour tout $x \in A \otimes_{\Q_p} \cE^I$. R\'eciproquement, consid\'erons le diagramme commutatif \'evident 
$$\xymatrix{ A \otimes_{\Q_p} \cE^I \ar@{->}[rr] ^{\iota_I}\ar@{->}[d]_\mu
& & \cE_A^I \ar@{->}[d] \\
 A \otimes_{\Q_p} \cE^J \ar@{->}[rr]^{\iota_J} & & \cE_A^J}.$$
D'apr\`es \cite[Prop. 17.4 (iii)]{schneider}, $|\mu(x)|_J=|x|_J$ pour tout $x \in A \otimes_{\Q_p} \cE^I$. De plus, le cas d'un segment entra\^ine qu'il existe une constante $C_J>0$ telle que $|x|_J \leq C_J
|\iota(x)|_J$ pour tout $x$ dans $A \otimes_{\Q_p} \cE^J$. Ainsi, $|.|_J$ et $|\iota_I(.)|_J$ sont \'equivalentes sur $A \otimes_{\Q_p} \cE^I$. On conclut car $\cE_A^I$ est un Fr\'echet contenant $A \otimes_{\Q_p} \cE^I$ comme sous-espace dense.\ps
Le premier point du (ii) est un exercice classique sans difficult\'e
laiss\'e au lecteur. Comme $\mathcal{A}[[T]]$ est complet 
pour la topologie $(p,T)$-adique, c'est un sous-espace ferm\'e de $\Ro_A^+$ et le (ii) suit. 
Pour le (iii), on constate sur la description donn\'ee plus haut de $\cE_\cA^{\dag, n}$ que
$\cE_\cA^{\dag,n}=T\cA[[T]] \oplus (\oplus_{i=0}^{n-1} T^i \Z_p \langle \frac{p}{T^n}
\rangle)$ est
une somme directe topologique, les termes de droite \'etant respectivement
munis de la topologie $(p,T)$-adique pour le premier et de la topologie $p$-adique
pour le second. Par
d\'efinition des $|.|_J$ sur $\Ro_{A,r}$, cette somme directe est une isom\'etrie pour
chaque $|.|_J$ (les termes de droites \'etant munis de la norme induite). Le premier point du (iii) se 
d\'eduit alors de celui du (ii). Le (iii) suit car $\cE_\cA^{\dag,n}$ est
complet pour la topologie faible. Comme il est complet pour cette topologie, cela
conclut.\ps
 Pour d\'emontrer le (iv) il faut voir si $f \in \cE_A^I$ a tous ses
coefficients dans $J$ (vue comme s\'erie de Laurent), alors $f \in
J\cE_A^I$. Soient $e_1,\dots,e_g$ une famille finie de g\'en\'erateurs de
$J$ comme $A$-module, d'apr\`es Tate la surjection ($A$-lin\'eaire) $A^g
\rightarrow I$ qui s'en d\'eduit est n\'ecessairement ouverte. En
particulier il existe une constante $C>0$ telle que tout \'el\'ement $x \in J$
s'\'ecrive sous la forme $\sum_i x_i
e_i$ avec $|x_i| \leq C |x|$ pour tout $i$. Ainsi, appliquant ceci \`a tous
les coefficients d'un $f \in \cE_A^I$ \`a coefficients dans $J$, il vient
que $f \in J \cE_A^I$.\ps
	Prouvons enfin (v). Si $I$ est un segment, alors $\cE_A^I=A
\widehat{\otimes}_{\Q_p} \cE^I$ est isomorphe \`a $A \langle t \rangle$
comme $A$-module, et il est bien connu que ce dernier est plat sur $A$. En
effet, si $\cA \subset A$ est un mod\`ele de $A$, alors $\cA\langle t
\rangle$ est plat sur $\cA[T]$ comme compl\'et\'e de ce dernier (qui est un
anneau noeth\'erien) pour la topologie $p$-adique, et donc plat sur $\cA$.
On conclut en inversant $p$. En g\'en\'eral, on \'ecrit 
$\cE_A^I$ est comme limite projective de $\cE_A^{I_n}$ pour une suite
arbitraire croissante $I_n$ de segments recouvrant $I$. Soit $0
\longrightarrow P \longrightarrow Q \longrightarrow R \longrightarrow 0$ une
suite exacte de $A$-modules de type fini. Par platitude de $\cE_A^{I_n}$ sur
$A$ pour tout $n\geq 0$ (cas pr\'ec\'edent), on dispose d'un syst\`eme projectif de suites
exactes $0    
\longrightarrow P_n \longrightarrow Q_n \longrightarrow R_n \longrightarrow 0$
o\`u $X_n:=X \otimes_A \cE_A^{I_n}$ (la famille $(X_n)$ d\'efinit donc un
$\cE_A^I$-module coh\'erent au sens de Schneider-Teitelbaum).
D'apr\`es~\cite[\S 3, Thm.]{schneiderteitelbaum}, c'est un fait g\'en\'eral
que cette suite reste exacte apr\`es passage \`a la limite projective sur $n$.
Pour conclure, il suffit de voir que si
$X$ est un $A$-module de type fini, alors l'application naturelle $$X
\otimes_A \cE_A^I \longrightarrow \projlim_n X_n$$
est un isomorphisme. C'est clair si $X$ est libre et cela suit en
g\'en\'eral si l'on choisit une pr\'esentation $L \rightarrow L' \rightarrow X
\rightarrow 0$ avec $L,L'$ libres de type fini sur $A$, et consid\`ere le diagramme
commutatif
$$\xymatrix{ L \otimes_A \cE_A^I \ar@{->}[d] \ar@{->}[rr] & & L' \otimes_A
\cE_A^I
\ar@{->}[d] \ar@{->}[rr] & & X\otimes_A \cE_A^I \ar@{->}[d] \rightarrow 0 \\ \projlim_n L_n \ar@{->}[rr] & & \projlim_n L'_n 
\ar@{->}[rr] & & \projlim_n X_n \rightarrow 0}$$ 
dont les suites horizontales sont exactes (par exactitude \`a droite du
produit tensoriel pourcelle du haut et par~~\cite[\S 3,
Thm.]{schneiderteitelbaum} pour celle du bas) et les deux verticales de
gauche sont des isomorphismes. L'assertion sur $\Ro_A$ suit car c'est une limite inductive filtrante de $\cE_A^I$.\end{pf}

\ps\ps

Le groupe $\Gamma:=\Z_p^\ast$ agit par automorphismes des $B_I$ par la formule \begin{equation}\label{formgamma}
\gamma(T)=(1+T)^\gamma-1=\sum_{n\geq 1} {\gamma \choose {n}} T^n \in \Z_p[[T]],
\end{equation}
de sorte que $|\gamma(t)|=|t|$ si $|t|<1$. Si $A$ est une alg\`ebre affino\"ide sur $\Q_p$, cette action de 
$\Gamma$ s'\'etend donc en une action $A$-lin\'eaire sur $\cE_A^I$, et en fait sur tous les anneaux 
introduits ci-dessus.

\begin{lemma}\label{estimeegamma} Soient $A$ une alg\`ebre affino\"ide sur $\Q_p$ et $\cA \subset A$ un mod\`ele. 
\begin{itemize}
\item[(i)] Si $I$ est un segment, l'application induite $\Gamma \rightarrow {\rm End}_A (\cE_A^I)$ est continue.\ps 
\item[(ii)] Si $1\leq n \leq m$, $I=[p^{-1/n},p^{-1/m}]$ et $\gamma \in 1+2p^M\Z_p$ alors $(\gamma-1)\OO^I_\cA \subset T^M \OO^I_\cA$.\ps
\end{itemize}
\end{lemma}

(Le terme \`a droite dans (i) est l'alg\`ebre des endomorphismes $A$-lin\'eaires
continus du $A$-module de Banach $\cE_A^I$, c'est une $A$-alg\`ebre de Banach pour la norme
d'op\'erateurs).\ps
\begin{pf} L'action de $\Gamma$ \'etant $A$-lin\'eaire, on peut supposer 
$A=\Q_p$ et $\mathcal{A}=\Z_p$ dans (i) et (ii). V\'erifions (i), il suffit 
par multiplicativit\'e de montrer la continuit\'e en $1 \in \Gamma$. \'Ecrivons 
$I=[p^{-a/n},p^{-b/m}]$ (resp. $I=[0,p^{-b/m}]$). Soient $U=\frac{T^m}{p^b}$ et 
$V=\frac{p^a}{T^n}$, ainsi que $\OO^I=\Z_p\langle U, V,T \rangle \subset \cE^I$ (resp. $\Z_p\langle U,T \rangle$)
le mod\`ele de $\cE^I$ d\'efini au lemme~\ref{lemmmodel}. Si $\gamma \in \Gamma$ alors 
$$\gamma(U)=\gamma(T)^m/p^b=\gamma^m U(1+\sum_{k\geq 2}
\frac{{ \gamma \choose k}}{\gamma} T^{k-1})^m \in U\Z_p\langle U \rangle [T]$$
$$\gamma(V)=\frac{p^a}{\gamma(T)^n}=\gamma^{-n} V (1+\sum_{k\geq 2}
{\frac{{\gamma \choose k}}{\gamma}} T^{k-1})^{-n} \in V\Z_p\langle U \rangle[T].$$
(car $\Z_p[[T]] \subset \Z_p\langle U \rangle [T]$). En particulier, $\Gamma$ pr\'eserve
$\OO^I$. Notons que $T^m \in p\OO^I$. Ainsi, si $M\geq 1$ est un entier suffisamment grand pour que
${\gamma \choose i} \in p\Z_p$ si $i=2,\dots,m$ et $\gamma \in 1+p^M\Z_p$, les formules ci-dessus entra\^inent que $1+p^M\Z_p$ agit trivialement sur
$\OO^I/p\OO^I$. L'identit\'e
\begin{equation}\label{identite} W^p-1 \equiv (W-1)^p \bmod p(W-1)\Z[W]\end{equation} 
montre alors que $(\gamma-1)\OO^I \subset p^N\OO^I$ si $\gamma \in
1+2p^{M+N-1}\Z_p$, ce qui termine la preuve du (i).\ps
V\'erifions donc (ii). Les formules pour $\gamma(T)$, $\gamma(U)$ et $\gamma(V)$ 
donn\'ees ci-dessus montrent alors que $\OO^I$ et $T\OO^I$ sont $\Gamma$-stables, puis 
que $\gamma$ agit trivialement sur $\OO^I/T\OO^I$ (car sur les images de $U$ et $V$) 
d\`es que $\gamma \in 1+p\Z_p$. Cela d\'emontre le (ii) pour $M=1$. On en d\'eduit que 
pour tout $i\geq 0$ et tout $\gamma \in 1+p\Z_p$, alors $(\gamma-1)T^i\OO^I \subset T^{i+1}\OO^I$. 
En effet, cela vient par r\'ecurrence sur $i$ de la formule $(\gamma-1)(ab)=(\gamma-1)(a)\gamma(b)+
(\gamma-1)(b)a$ et de ce que $$(\gamma-1)(T)\in (T^2,pT)\Z_p[[T]] \subset T^2 \OO^I$$ 
(car $p=VT^n \in T \OO^I$). En particulier, $(\gamma-1)^jT^i\OO^I \subset
T^{i+j}\OO^I$ pour tout $i,j\geq 0$ et $\gamma \in 1+p\Z_p$. Le (ii) suit alors pour tout $M$, encore par r\'ecurrence sur $M$, en utilisant cette fois-ci l'identit\'e (\ref{identite}) et encore le fait que $p=VT^n \in T \OO^I$. 
\end{pf}

\subsection{Pr\'eliminaire sur les familles de $\Gamma$-modules}\label{checkgamma}
\newcommand{\Mat}{{\rm Mat}}
Soit $A$ une $\Q_p$-alg\`ebre affino\"ide et $N\geq 0$ un entier. On d\'efinit un 
{\it $\Gamma$-module sur $\cE_{\cA}^{\dag,N}$} comme \'etant un $\cE_{\cA}^{\dag,N}$-module $D$ 
libre de rang fini muni d'une action semi-lin\'eaire de
$\Gamma$ qui soit continue, {\it i.e.} telle qu'il existe une base $e_1,\dots,e_d$ de $D$ sur 
$\cE_{\cA}^{\dag,N}$ pour laquelle l'application 
$\gamma \in \Gamma \mapsto \Mat(\gamma) \in M_d(\cE_{\cA}^{\dag,N})$ associant \`a $\gamma$ sa 
matrice dans la base $(e_i)$ soit continue coefficient par coefficient. On v\'erifirait ais\'ement 
\`a l'aide des lemmes~\ref{estimeegamma} (i) et~\ref{topo} (ii) et (iii) que si cela vaut pour une base alors cela vaut pour toutes. \ps

Pour tout segment $I \subset [p^{-\frac{1}{N}},1[$, on a un morphisme naturel $\cE_\cA^{\dag,N} \rightarrow 
\OO_\cA^I$, de sorte qu'il y a un sens \`a consid\'erer $D^I=D \otimes_{\cE_\cA^{\dag,N}} \OO_\cA^I$. 

\begin{lemma}\label{estimeeD} Soit $s\geq 0$ un entier. Il existe un entier $M\geq 1$ tel que $\forall \gamma \in 1+p^M\Z_p$ : \ps \begin{itemize}
\item[(a)] $(\gamma-1)D^I \subset T^s D^I$ pour tout $I=[p^{-1/m},p^{-1/m'}]$ avec $1\leq m \leq m'$ et $m\geq
N$,\ps \ps
\item[(b)] et de plus $(\gamma-1)D^{[0,p^{-1/2}]} \subset p^s D^{[0,p^{-1/2}]}$ si $N=0$.\end{itemize}\end{lemma}

\begin{pf} Supposons d'abord $N\geq 1$. Par continuit\'e de $\Gamma$, on peut choisir $M\geq s$ tel que pour tout 
$\gamma \in 1+p^{M}\Z_p$ on ait $\Mat(\gamma)-{\rm id} \in p^sM_d(\cE_{\cA}^{\dag,n})+T^sM_d(\cA[[T]])$. Mais 
$p \in T^N\cE_{\cA}^{\dag,N}$, donc $\Mat(\gamma)-{\rm id} \in T^sM_d(\cE_{\cA}^{\dag,N})$. En particulier, 
si $1\leq m \leq m'$, $m\geq N$ et $I=[p^{-1/m},p^{-1/m'}]$, alors $$(\gamma-1)(e_i)
\in T^s D^I,\,\, \,\,\,\forall \gamma \in 1+p^M\Z_p, \,\,\forall i=1,\dots,d.$$
On conclut alors le (a) par l'identit\'e $$(\gamma-1)(\sum_i x_i e_i)=\sum_i
(\gamma-1)(x_i)\gamma(e_i)+x_i(\gamma-1)(e_i),$$
et le lemme~\ref{estimeegamma} (ii). Si $N=0$, le fait que l'on puisse choisir $M$ de sorte que (a) soit satisfait d\'ecoule du cas pr\'ec\'edent en consid\'erant le $\Gamma$-module sur $\cE_\cA^{\dag,1}$ obtenu par extension des scalaires. Pour v\'erifier (b) on proc\`ede de m\^eme que ci-dessus en remarquant par exemple que $T^2 \in p\OO_I$ 
si $I=[0,p^{-1/2}]$ et en utilisant le lemme~\ref{estimeegamma} (i).
\end{pf}

Terminons par une proposition qui jouera un r\^ole important dans la suite.
Fixons $D$ un $\Gamma$-module sur $\cE_{\cA}^{\dag,N}$. Si $\gamma \in \Gamma$, consid\'erons l'application 
$$G_\gamma : D \rightarrow D, \, \, x \mapsto Tx+(1+T) \cdot (\gamma-1)(x),$$
c'est un endomorphisme $A$-lin\'eaire de $D$. On consid\`ere le morphisme de
$A$-alg\`ebres
$A[G_\gamma] \rightarrow \Ro_{A,p^{-1/n}}$ envoyant $G_\gamma$ sur
$T$. \ps 

\begin{prop}\label{propcle} Il existe un entier $M\geq 1$ tel que pour tout $\gamma \in 1+p^M\Z_p$, et pour tout $n\geq N$, le $A[G_\gamma]$-module $D^{(n)} = D \otimes_{\cE_{\cA}^{\dag,N}} \Ro_{A,p^{-1/n}}$ s'\'etend de mani\`ere unique en un $\Ro_{A,p^{-1/n}}$-module not\'e\footnote{Comme $A$-module, on a donc
$D^{'(n)}=D^{(n)}$...} $D^{'(n)}$ tel que l'application structurelle 
$\Ro_{A,p^{-1/n}} \times D^{'(n)} \rightarrow D^{'(n)}$ soit continue. De plus, 
si $(e_i)$ est une base de $D$ sur $\cE_{\cA}^{\dag,N}$ alors $(e_i \otimes 1)$ est une base de $D^{'(n)}$ sur $\Ro_{A,p^{-1/n}}$. 
\end{prop}

La topologie sous-entendue sur $D^{'(n)}$ dans cet \'enonc\'e est celle de
$\Ro_{A,p^{-1/n}}$-module de Fr\'echet $D^{(n)}$ sous-jacent (qui est libre de rang
fini par d\'efinition).

\begin{pf} Appliquons le lemme pr\'ec\'edent pour $s=2$. Fixons une fois pour toutes $\gamma \in 1+p^M \Z_p$, o\`u $M$ est donn\'e par ce lemme. Soit $I=[p^{-1/m},p^{-1/m'}]$ avec $m\geq N$ et $m\geq 1$, et consid\'erons $D^I$. C'est un $\OO_\cA^I$-module libre de rang $d$, il est en particulier complet pour la topologie $p$-adique, et aussi pour la topologie $T$-adique (car $\OO_\cA^I$ l'est, puisque $T^{m'}\OO_\cA^I \subset p\OO_\cA^I$). L'anneau $\cal{B}^I=\End_\cA(\OO_\cA^I)$ est donc lui aussi complet pour les topologies $p$-adique et $T$-adique. En particulier, si on pose $$\psi(x):=\frac{1+T}{T}(\gamma-1)(x)$$
alors $\psi(D^I) \subset T D^I$ par choix de $M$, en particulier $\psi \in \cal{B}^I$, mais aussi $\sum_{k\geq 0} (-1)^k\psi^k$ converge dans $\cal{B}^I$, vers un inverse de $1+\psi$. De plus, on a \'evidemment $\psi(TD^I) \subset TD^I$, donc $\psi$ induit l'endomorphisme nul sur $D^I/TD^I$. Si on note $m_u$ la multiplication par $u$, nous avons donc montr\'e que : \begin{itemize}
\item[(i)] $G_\gamma=m_T \cdot (1+\psi)$ dans $\cal{B}^I$ et $G_\gamma \equiv m_T$ dans $\End_{\OO^I/(T)}(D^I/TD^I)$.
\item[(ii)] $\frac{G_\gamma^{m'}}{p}=m_{\frac{T^{m'}}{p}}\cdot (1+\psi)^{m'} \in \cal{B}^I$ et $\frac{G_\gamma^{m'}}{p} \equiv m_{\frac{T^{m'}}{p}}$ dans $\End_{\OO^I/(T)}(D^I/TD^I)$.\ps
\item[(iii)] $\frac{p}{G_\gamma^m}=m_{\frac{p}{T^m}} \cdot (\sum_{k\geq 0} (-1)^k \psi^k)^m \in \cal{B}^I$ et $\frac{p}{G_\gamma^{m}} \equiv m_{\frac{p}{T^m}}$ dans $\End_{\OO^I/(T)}(D^I/TD^I)$.
\end{itemize}
Comme $\cal{B}^I$ est complet pour la topologie $p$-adique, 
et comme $\OO^I$ est par d\'efinition de quotient de $\Z_p\langle
T,U,V\rangle/(pU-T^m,T^{m'}V-p)$ par sa $p$-torsion (lemme~\ref{lemmmodel}), 
il d\'ecoule de (i), (ii) et (iii) que le morphisme naturel $\cA[\G_\gamma] \rightarrow \cal{B}^I$ 
s'\'etend en un morphisme continu $\cE_\cA^I \rightarrow \cal{B}^I$. On note $D^{'I}$ le groupe 
ab\'elien $D^I$ muni de cette structure de $\cE_\cA^I$-module. Comme $1+\psi$ est inversible dans 
$\cal{B}^I$, notons que $T^kD^{'I}=G_\gamma^k D^I=T^k D^I$ pour tout $k\geq 0$, puis que $D^{'I}$ est 
complet pour la topologie $T$-adique et sans $T$-torsion. Consid\'erons l'application identit\'e 
$$D^I/TD^I \rightarrow D^{'I}/TD^{'I}.$$
Le (i), (ii) et (iii) ci-dessus assurent que cette application est un morphisme de
$\OO_\cA^I$-modules. En particulier, $D^{'I}/(T)$ est libre comme $\cE_\cA^I/(T)$-module. 
Comme $D^{'I}$ et $\cE_\cA^I$ sont complets pour la topologie $T$-adique et sans $T$-torsion, 
un raisonnement standard montre que $D^{'I}$ est libre sur $\cE_\cA^I$, et qu'une famille $
(e_i)$ est une base de $D^I$ si et seulement si c'est une base de $D^{'I}$. Quand $n=0$, 
on traite le cas de $I=[0,p^{-1/2}]$ par un raisonnement enti\`erement analogue en ne 
consid\'erant que la topologie $p$-adique (et non pas $T$-adique, c'est en fait seulement 
plus simple il n'y a pas de condition de type (iii) \`a v\'erifier).
\ps

Pour terminer, il ne reste qu'\`a "recoller" les $D^I[1/p]$. Pour cela, fixons $n\geq N$ et 
consid\'erons l'ensemble $\mathcal{I}$ des intervalles de $[p^{-1/n},1[$ de la forme $[p^{-1/m},p^{-1/m'}]$ avec $m\leq m'$, ou de la forme $[0,p^{-1/2}]$ (ce qui ne se produit que si $n=0$). Si $J \subset I$ sont dans $\mathcal{I}$ on a bien s\^ur un morphisme de restriction $$r_{I,J}: D^I[1/p] \rightarrow D^J[1/p].$$ 
Si $I,J \in \mathcal{I}$, alors $I\cap J=\emptyset$ ou $I\cap J \in \mathcal{I}$. Comme 
les $B_I$ avec $I\in \mathcal{I}$ recouvrent admissiblement $B_{[p^{-1/n},1[}$, il vient 
que $D^{(n)}$ (resp. $\Ro_{A,p^{-1/n}}$) s'identifie \`a la limite projective sur $\mathcal{I}$ 
des $\cE_A^I$-modules $D^I[1/p]$ (resp. des $\cE_A^I$). Il est de plus imm\'ediat sur la construction 
ci-dessus que si $f \in \cE_A^I$ et $v \in D^{I}[1/p]$, alors $r_{I,J}(f \ast v)=r_{I,J}(f) \ast r_{I,J}(v)$ 
o\`u $\ast$ d\'esigne ici la structure de module de $D^{'I}[1/p]$ et $D^{'J}[1/p]$ sur $\cE_A^I$ et $\cE_A^J$. 
Cela nous permet d'une part de munir $D^{(n)}$ d'une structure de $\Ro_{A,p^{-1/n}}$-module, disons $D^{'(n)}$, 
en posant $(f_I)\ast (v_I):=(f_I\ast v_I)$. D'autre part, si $(e_i)$ est une base de $D$
sur $\cE_{\cA}^{\dag,N}$ nous avons vu que 
$e_i \otimes 1$ est une base de $D^{'I}[1/p]$ sur $\cE_A^I$, on en d\'eduit que $D^{'(n)}$ 
est libre sur $\Ro_{A,p^{-1/n}}$ de base $e_i \otimes 1$.\ps
	Le choix d'une base de $D$ munit $D^{(n)} = \oplus_i \Ro_{A,p^{-1/n}} e_i$ d'une structure d'espace 
de Fr\'echet qui ne d\'epend pas du choix de la base $(e_i)$. La continuit\'e de $\Ro_{A,p^{-1/n}} \times 
D^{'(n)} \rightarrow D^{'(n)}$ se d\'eduit alors de celle des $\cE_A^I \times D^{'I}[1/p] \rightarrow 
D^{'I}[1/p]$. L'assertion d'unicit\'e vient plus pr\'ecis\'ement de ce qu'il existe au plus une structure 
de $\Ro_{A,p^{-1/n}}$-module telle que pour tout $v \in D^{(n)}$, l'application $f \mapsto f\ast v$ soit 
continue, car $A[T,T^{-1}]$ (resp. $A[T]$) est dense dans $\Ro_{A,p^{-1/n}}$ si $n\geq 1$ (resp. si $n=N=0$).
\end{pf}

\subsection{Familles de $\fg$-modules sur l'anneau de Robba} 

Soit $A$ une $\Q_p$-alg\`ebre affino\"ide. Les anneaux $\Ro^+_A$ et $\Ro_A$ sont munis d'un endomorphisme 
d'anneaux $\varphi$ d\'efini par $$\varphi(f)(T)=f((1+T)^p-1)$$ et qui commute \`a l'action de $\Gamma$.
Plus pr\'ecis\'ement, supposons $r=0$ ou $r>p^{-\frac{1}{(p-1)}}$, on dispose d'un morphisme analytique 
$\varphi_\ast : B_{[r,1[} \rightarrow B_{[r^p,1[}$ d\'efini par $t \mapsto (1+t)^p-1$ et $\varphi$
est par d\'efinition le morphisme $\Ro_{A,r^p} \rightarrow \Ro_{A,r}$ qui s'en
d\'eduit. Il est \'evident sur la formule pour $\varphi_\ast$ que $\varphi$
commute \`a l'action de $\Gamma$. \ps

Un $\fg$-module sur $\Ro_A$ est un $\Ro_A$-module $D$ libre et de rang fini muni d'actions semi-lin\'eaires 
de $\varphi$ et $\Gamma$ qui commutent et satisfaisant les axiomes suivants. D'une part on demande que 
$\varphi$ envoie une $\Ro_A$-base de
$D$ sur une $\Ro_A$-base de $D$, i.e. $\Ro \varphi(D)= D$. D'autre part, on demande que l'action de 
$\Gamma$ soit continue au sens suivant : il existe 
une $\Ro_A$-base $e_1,\dots,e_d$ de $D$ et $r \in [0,1[$ tels que si $\gamma \mapsto \Mat(\gamma)$ 
d\'esigne la matrice de $\gamma \in \Gamma$
dans cette base, alors $M(\gamma) \in M_d(\Ro_{A,r})$ pour tout $\gamma \in \Gamma$ et l'application 
$\Gamma \rightarrow
M_d(\Ro_{A,r})$, $\gamma \mapsto M(\gamma)$, est continue (coefficient par 
coefficient). On dira que $D$ est {\it $\Gamma$-born\'e} si on peut trouver un mod\`ele $\mathcal{A} 
\subset A$, un entier $n\geq 0$, et une $\Ro_A$-base $e_i$ de $D$ dans laquelle 
$\Mat(\gamma) \in M_d(\cE_{\cA}^n)$ pour tout $\gamma \in \Gamma$. Il
r\'esulte du lemme~\ref{topo} (iii) que le $\Gamma$-module $\oplus_i
\cE_{\cA}^n e_i$ qui s'en d\'eduit est bien un $\Gamma$-module sur $\cE_{\cA}^n$
au sens du \S~\ref{checkgamma}. \ps
\newcommand{\FG}{{\rm FG}}

Les $\fg$-modules sur $\Ro_A$ forment un cat\'egorie $A$-lin\'eaire $\FG_A$ s'il on
consid\`ere pour $\Hom_{\FG_A}(D_1,D_2)$ les applications $\Ro_A$-lin\'eaires qui commutent \`a $\Gamma$ et 
\`a $\varphi$. \ps 

L'int\'er\^et des $\fg$-modules sur $\Ro_A$ r\'eside dans leurs liens avec les repr\'esentations 
continues $\Gal(\Qpb/\Qp) \rightarrow \GL_d(A)$ (Fontaine,
Cherbonnier-Colmez, Kedlaya, Berger-Colmez, Kedlaya-Liu) pour lequel nous renvoyons \`a
Berger-Colmez~\cite{bergercolmez} et Kedlaya-Liu~\cite{kedliu}. La condition $\Gamma$-born\'ee introduite 
ci-dessus n'est pas standard, et appara\^it ici pour des raisons techniques. Bien que nous n'utiliserons 
pas ce r\'esultat, mentionnons que la m\'ethode de
Berger-Colmez~\cite{bergercolmez}, g\'en\'eralisant un r\'esultat de 
Colmez-Cherbonnier, assure que si $M$ est un $\cA$-module libre muni d'une action $\cA$-lin\'eaire continue 
de\footnote{Rappelons que d'apr\`es~\cite[lemme 3.18]{chhecke}, une repr\'esentation continue
$\rho : G \rightarrow \GL_d(A)$
d'un groupe profini $G$ \'etant donn\'ee, on peut toujours trouver un
recouvrement fini de $X$ par des affino\"ides $U_i$, ainsi que des mod\`eles
$\cal{A}_i \subset \OO(U_i)$, tels que pour tout $i$ la repr\'esentation $\rho \otimes_A \OO(U_i)$
de $G$ sur $\OO(U_i)^d$ stabilise un sous-$\cA_i$-module libre $L_i$ de rang
$d$ tel que $L_i[1/p]=\OO(U_i)^d$.} 
 $\Gal(\Qpb/\Q_p)$, alors on peut lui associer un $\fg$-module sur $\Ro_A$ qui est $\Gamma$-born\'e.
\ps

Enfin, d\'efinissons un $\fg$-module sur $\Ro_A^+$ comme \'etant un $\Ro^+_A$-module $D$ libre de rang fini 
muni d'actions semi-lin\'eaires de $\varphi$ et $\Gamma$ qui commutent, telles que $\varphi(D)$ contienne 
une $\Ro^+_A$-base de $D$, et telles qu'il existe une $\Ro^+_A$-base $e_1,\dots,e_d$ de $D$ dans laquelle $\gamma \mapsto \Mat(\gamma)$ (la matrice de $\gamma \in \Gamma$
dans cette base) soit continue. On dira encore que $D$ est {\it $\Gamma$-born\'e} si on peut 
trouver un mod\`ele $\mathcal{A} \subset A$ et une $\Ro^+_A$-base $e_i$ de $D$ dans laquelle 
$\Mat(\gamma) \in M_d(\cA[[T]])$ pour tout $\gamma \in \Gamma$. 
 Il
r\'esulte du lemme~\ref{topo} (ii) que le $\Gamma$-module $\oplus_i
\cA[[T]] e_i$ qui s'en d\'eduit est bien un $\Gamma$-module sur $\cA[[T]]$
au sens du \S~\ref{checkgamma}. \ps
\ps\ps

\newcommand{\cT}{{\mathcal T}}
Les $\fg$-modules sur $\Ro_A$ qui nous int\'eressent principalement dans cet article sont les $\fg$-modules triangulins, variante en famille d'une notion introduite par Colmez dans~\cite{colmeztri}. Soit $\cT$ l'espace analytique $p$-adique param\'etrant les caract\`eres continus de $\Q_p^\ast$ : pour toute alg\`ebre affino\"ide $A$, $\cT(A)$ est l'ensemble des morphismes continus $\Qp^\ast \rightarrow A^\ast$. Il est bien connu que cet espace est
isomorphe au produit de $\mathbb{G}_m$ par l'espace $\WW$ des 
caract\`eres continus $p$-adiques de $\Z_p^\ast$, lui-m\^eme
\'etant une r\'eunion disjointe finie de boules unit\'es ouvertes. Si
$\delta \in \cT(A)$, i.e. si $\delta : \Q_p^\ast \rightarrow A^\ast$ est un
morphisme continu de groupes, on d\'efinit un $\fg$-module
$\Ro_A(\delta)$ de rang $1$ sur $\Ro_A$, disons $\Ro_A(\delta)=\Ro_A e$, 
par la formule $\gamma(e)=\delta(\gamma)e$
pour tout $\gamma \in \Gamma$ et $\varphi(e)=\delta(p)e$. On d\'efinit de m\^eme $\Ro_A^+(\delta)$ de 
mani\`ere \'evidente.

\begin{definition} Un $\fg$-module triangulin sur $\Ro_A$ est la donn\'ee d'un $\fg$-module $D$ sur $\Ro_A$ et d'une suite croissante $(\Fil_i(D))_{i=0,\dots,d}$, $d={\rm rg}_{\Ro_A}(D)$, de sous $\Ro_A$-modules de $D$ stables par $\varphi$ et $\Gamma$ telle que
$\Fil_0(D)=0$, $\Fil_d(D)=D$, et telle que pour chaque $i=1,\dots,d-1$, 
$\Fil_i(D)/\Fil_{i-1}(D) \simeq \Ro_A(\delta_i)$ pour un certain
$\delta_i \in \cT(A)$. 
\end{definition}

Nous verrons plus bas (Lemme~\ref{classrg1}) que la suite des $\delta_i$ est 
uniquement d\'etermin\'ee par $(\Fil_i(D))$, nous l'appelons le {\it param\`etre} de $D$.  Notons 
aussi que $\Fil_i(D)$ est libre de rang $i$ sur
$\Ro_A$, et facteur direct dans $D$ comme $\Ro_A$-module. \ps

Terminons ce paragraphe par un sorite sur l'extension des scalaires. Soit $B$ une $A$-alg\`ebre affino\"ide. 
On dispose pour chaque $0<r<1$ d'un morphisme continu d'anneaux $\Ro_{A,r} \rightarrow \Ro_{B,r}$ induisant 
\`a la limite un morphisme $\Ro_A \rightarrow \Ro_B$. Si $D$ est un $\fg$-module sur $\Ro_A$ on note 
$$D \widehat{\otimes}_A B$$ le $\Ro_B$-module $D \otimes_{\Ro_A} \Ro_B$ : c'est un $\fg$-module sur 
$\Ro_B$ de mani\`ere naturelle. On pourrait justifier cette notation en le voyant comme un produit 
tensoriel compl\'et\'e mais cela ne sera pas n\'ecessaire. On a d\'efini ainsi un foncteur 
$ -\,\widehat{\otimes}_A B : \fg/A \rightarrow \fg/B$. Si $I$ est un id\'eal de
$A$, le lemme~\ref{topo} (iv) entra\^ine que $\Ro_A/I\Ro_A=\Ro_{A/I}$ et donc que $D \widehat{\otimes}_A A/I = D/ID$ 
pour tout $\fg$-module $D$ sur $\Ro_A$. 
Si $I=m_x$ est l'id\'eal maximal correspondant \`a $x \in {\rm Sp}(A)$, on posera aussi $$D_x:=D/m_xD,$$ c'est 
un $\fg$-module sur $k(x):=A/m_x$ (une extension finie de $\Q_p$). \ps\ps

Par exemple, si $\delta \in \cT(A)$ alors $\Ro_A(\delta) \widehat{\otimes}_A B = \Ro_B(\delta')$ o\`u $\delta' 
\in \cT(B)$ est le caract\`ere $\delta$ compos\'e par $A^\ast \rightarrow B^\ast$. Notons que si $0 \rightarrow D 
\rightarrow D' \rightarrow D'' \rightarrow 0$ est une suite exacte de $\fg$-modules sur $\Ro_A$, elle est 
scind\'ee comme suite de $\Ro_A$-modules, et induit donc une suite exacte de $\fg$-modules sur $\Ro_B$ apr\`es 
extension des scalaires. En particulier, si $(D,\Fil_\bullet(D))$ est $\fg$-module triangulin sur $\Ro_A$, alors 
$D\widehat{\otimes}_A B$ est triangulin sur $\Ro_B$ pour la filtration $\Fil_i(D) \widehat{\otimes} B$.\ps\ps

\begin{lemma}\label{classrg1} Si $\delta,\delta' \in \cT(A)$, alors $\Ro_A(\delta) \simeq \Ro_A(\delta')$ si et 
seulement si $\delta=\delta'$. 
\end{lemma}

\begin{pf} En effet,  si $A$ est artinien c'est~\cite[Prop. 2.3.1]{bch}. En g\'en\'eral, on remarque que si $I$ est un 
id\'eal de $A$ alors $\Ro_{A}(\delta) \otimes_A A/I = \Ro_{A/I}(\delta \bmod I)$. On conclut car si $A$ est une 
$\Q_p$-alg\`ebre affino\"ide alors l'intersection de ses id\'eaux de codimension finie est
nulle par le th\'eor\`eme d'intersection de Krull, et donc $A$ se plonge dans le produit des $A/I$ avec $I$ de
codimension finie.
\end{pf}

Terminons par une question, dont une r\'eponse (affirmative) ne
semble connue que lorsque $A$ est artinien :\ps
\pn
{\bf Question:} Est-ce que tout $\fg$-module de rang $1$ sur $\Ro_A$
est isomorphe \`a un $\Ro_A(\delta)$ ?
\ps
\ps

\section{Cohomologie des $\fg$-modules triangulins sur $\Ro_A$}\label{cohomologie}

L'objectif de cette partie est de calculer la cohomologie des $\fg$-modules triangulins sur $\Ro_A$. \ps\ps

\subsection{G\'en\'eralit\'es sur la cohomologie des $\fg$-modules}\label{gencoho} Si $D$ est un 
$\fg$-module sur $\Ro_A$, et si $\gamma \in \Gamma$ est un g\'en\'erateur 
topologique\footnote{Un tel g\'en\'erateur n'existe bien s\^ur que pour $p>2$. Quand
$p=2$, on choisit pour $\gamma \in \Gamma$ un \'el\'ement engendrant
topologiquement $\Gamma/\Gamma_{\rm tors}$ o\`u $\Gamma_{\rm tors}=\{\pm
1\}$. On d\'efinit ensuite $C(D)^\bullet$ de la 
m\^eme mani\`ere \`a ceci pr\`es que de $D$ y est partout
remplac\'e par ses invariants  $D^{\Gamma_{\rm tors}}$ sous le groupe
fini $\Gamma_{\rm tors}$, ce qui n'alt\`ere aucun des arguments qui
suivent.} on rappelle que suivant Fontaine et
Herr~\cite{herr1} on dispose du complexe  $C_{\varphi,\gamma}(D)^\bullet$ : $$0 \rightarrow D \overset{x\mapsto (\varphi-1)x+(\gamma-1)x}{\longrightarrow} 
D \times D \overset{(x,y)\mapsto (\gamma-1)x-(\varphi-1)y}{\longrightarrow} D \rightarrow 0,$$
le premier $D$ \'etant plac\'e en degr\'e $0$. On d\'esigne par $H^i(D)$ la cohomologie de ce complexe, ce sont donc des $A$-modules nuls en degr\'e $i \notin \{0,1,2\}$. 
Par d\'efinition, $H^0(D)=\Hom_{\FG_A}(\Ro_A,D)$. De plus, \ps\ps

\begin{lemma} $H^1(D)$ est canoniquement isomorphe \`a ${\rm Ext}_{\FG_A}(\Ro_A,D)$.
\end{lemma}

\begin{pf} Donner une action de $\varphi$ et $\gamma$ sur $D\oplus \Ro_A$ \'etendant 
la structure de $\fg$-module de $D$ est \'equivalent \`a donner $x:=(\varphi-1)(e) \in D$ 
et $y:=(\gamma-1)(e) \in D$, $\varphi$ et $\gamma$ commutant si et seulement si $(x,y)$ est 
dans $Z^1(C_{\varphi,\gamma}(D))$. L'action de $\gamma$ s'\'etend automatiquement en une 
action continue de $\Gamma$. En effet, il d\'ecoule du Lemme~\ref{estimeegamma}
(i) que si $f \in \Ro_{A,r}$ alors la suite $(1+\gamma+\cdots+\gamma^{p^n-1})f$ tend vers $0$. Cela 
vaut donc aussi si $f \in M_d(\Ro_{A,r})$, et donc tout $1$-cocycle $\gamma^\Z \rightarrow 
M_d(\Ro_{A,r})$ s'\'etend en un cocycle continu $\Gamma \rightarrow
M_d(\Ro_{A,r})$. On v\'erifie imm\'ediatement que deux $1$-cocycles donnent des extensions isomorphes si et 
seulement si ils diff\`erent d'un cobord.
\end{pf}

Il se trouve que toujours suivant Fontaine et Herr, un autre complexe est
reli\'e \`a la 
cohomologie de $D$. Supposons $r=0$ ou $r>p^{-\frac{1}{p(p-1)}}$. L'application
$(\Ro_{A,r^p})^p \rightarrow \Ro_{A,r}$ d\'efinie par
$$(f_0,\dots,f_{p-1}) \mapsto \sum_{i=0}^{p-1} (1+T)^i\varphi(f_i)$$
est alors un isomorphisme topologique. En effet, c'est un r\'esultat
standard quand $A=\Q_p$ et le lemme~\ref{topo} (i) nous y ram\`ene en
g\'en\'eral. En particulier, $\varphi$ est fini et plat de degr\'e $p$, continu et
injectif. On d\'efinit alors $\psi : \Ro_{A,r} \rightarrow \Ro_{A,r}$ par la
formule $\varphi \psi = \frac{1}{p} {\rm
trace}_{\Ro_{A,r}/\varphi(\Ro_{A,r^p})}(\varphi)$. Sur $A\otimes_{\Q_p}\Q_p(\mu_p)$ on a 
donc la formule $\varphi\psi(f)= \frac{1}{p}\sum_{\zeta^p=1}f(\zeta (1+T)-1)$.
Le calcul de la trace des $(1+T)^i$ assure que si 
$f=\sum_{i=0}^{p-1} (1+T)^i\varphi(f_i)$ alors $\psi(f)=f_0$. En particulier, 
$\psi : \Ro_{A,r} \rightarrow \Ro_{A,r^p}$ est continu.  \ps

Soit $D$ un $\fg$-module sur $\Ro_A$ ou $\Ro_A^+$. Comme $D$ a par d\'efinition une
$\Ro_A$ ou $\Ro_A^+$-base dans $\varphi(D)$, on a encore une d\'ecomposition
$$D=\oplus_{i=0}^{p-1} (1+T)^i \varphi(D).$$
On peut donc d\'efinir un op\'erateur $\psi: D \rightarrow D$ par la formule $$\psi(\sum_{i=0}^{p-1} (1+T)^i \varphi(x_i))=x_0.$$ Il est $A$-lin\'eaire surjectif, commute \`a l'action de $\Gamma$, et satisfait $\psi \varphi = {\rm id}$. Enfin, si $u \in \Hom_{\FG_A}(D_1,D_2)$, alors $u \cdot \psi = \psi \cdot u$. Le complexe $C_{\psi,\gamma}(D)^\bullet$ est alors d\'efinit de la m\^eme mani\`ere que $C_{\varphi,\gamma}(D)^\bullet$ \`a ceci pr\`es que $\varphi$ est remplac\'e par $\psi$. Un calcul imm\'ediat montre que l'on dispose d'un morphisme $$\eta : C_{\varphi,\gamma}(D)^\bullet \rightarrow C_{\psi,\gamma}(D)^\bullet$$
qui vaut l'identit\'e en degr\'e $0$, $(x,y)\mapsto (-\psi(x),y)$ en degr\'e $1$, et $-\psi$ en degr\'e $2$. Le morphisme $\eta$ est surjectif car $\psi$ l'est, son noyau \'etant le complexe 
$$0 \longrightarrow  0 \longrightarrow D^{\psi=0} \overset{\gamma-1}{\longrightarrow} D^{\psi=0} \longrightarrow 0.$$ 
En particulier, $C_{\psi,\gamma}(D)^\bullet$ et $C_{\varphi,\gamma}(D)^\bullet$ sont quasi-isomorphes si $\gamma-1$ est bijectif sur $D^{\psi=0}$. La proposition suivante est imm\'ediate.

\begin{prop}\label{devissagecoho} Soit $D$ un $\fg$-module sur $\Ro_A$. Si $\gamma-1$ est bijectif sur $D^{\psi=0}$ alors on a des
identifications naturelles $H^0(D)=D^{\psi=1,\gamma=1}$,
$H^2(D)=D/(\psi-1,\gamma-1)$, ainsi qu'une suite exacte naturelle $$
0 \rightarrow D^{\psi=1}/(\gamma-1) \overset{y \mapsto (0,y)}{\longrightarrow} H^1(D) \overset{(x,y)\mapsto \overline{x}}{\longrightarrow} (D/(\psi-1))^{\gamma=1} \rightarrow
0.$$ 
Enfin, si on pose $C(D)=(\varphi-1)D^{\psi=1} \subset D^{\psi=0}$, on a une
suite exacte naturelle $$0 \rightarrow D^{\varphi=1}/(\gamma-1)\rightarrow
D^{\psi=1}/(\gamma-1) \rightarrow C(D)/(\gamma-1) \rightarrow 0.$$
\end{prop}

En th\'eorie des $\fg$-modules de Fontaine classique, un r\'esultat de
Herr~\cite[Thm. 3.8]{herr1} assure que $\gamma-1$ est toujours bijectif sur $D^{\psi=0}$. Dans le cadre ci-dessus, un r\'esultat de Colmez assure aussi que c'est toujours le cas si $A$ est un corps \cite[Prop. 5.1.19]{colmezgros}.  Nous allons d\'emontrer que c'est aussi le cas en g\'en\'eral sous une hypoth\`ese assez faible sur $D$. \ps

\begin{thm}\label{invgamma} Soit $D$ un $\fg$-module sur $\Ro_A$ qui est $\Gamma$-born\'e. Alors $\gamma-1$ est bijectif sur $D^{\psi=0}$.\par
Plus g\'en\'eralement, supposons que $D$ est un $\fg$-module sur $\Ro_A$ poss\'edant une suite croissante $D_1 \subset D_2 \subset \cdots \subset D_s$ de sous-$\Ro_A$-modules qui sont facteurs directs comme $\Ro_A$-modules, et de plus stables par $\varphi$ et $\Gamma$. Supposons enfin que les $D_{i+1}/D_i$ sont $\Gamma$-born\'es. Alors $\gamma-1$ est bijectif sur $D^{\psi=0}$. 
\end{thm}

\begin{pf} En effet, si $0 \rightarrow D_1 \rightarrow D_2 \rightarrow D_3 \rightarrow 0$ alors la surjectivit\'e de $\psi$ (sur $D_1$) entra\^ine que la suite associ\'ee $$ 0 \rightarrow D_1^{\psi=0} \rightarrow D_2^{\psi=0} \rightarrow D_3^{\psi=0} \rightarrow 0$$
est exacte. Comme elle est $\gamma$-\'equivariante, il vient que si $\gamma-1$ est bijectif sur $D_i^{\psi=0}$ pour $i=1,3$ alors il l'est aussi pour $i=2$, de sorte que le second cas suit du premier, que nous consid\'erons maintenant.\par
Soient $e_1,\dots,e_d$ une $\Ro_A$-base de $D$, $\cA \subset A$ un mod\`ele, et $N\geq 0$, tels que $\cal{D}:=\oplus_i \cE_\cA^{\dag,N} e_i$ soit stable par $\Gamma$. Choisissons un entier $M$ comme dans la proposition~\ref{propcle}. Comme $\gamma-1$ divise $\gamma^{(p-1)p^{M-1}}-1$ dans $\Z[\gamma]$, il suffit de montrer que $\gamma-1$ est bijectif sur $D^{\psi=0}$ si $\gamma \in 1+p^M\Z_p^\ast$. Posons $\gamma_0=1+p^M \in \Gamma$. La relation $\gamma_0(1+T)=(1+T)\varphi^M(1+T)$ entra\^ine pour tout $x$ dans $D$ $$(\gamma_0-1)((1+T)\varphi^M(x))=(1+T)\varphi^M((1+T)\gamma_0(x)-x)=(1+T)\varphi^M(G_{\gamma_0}(x)).$$ 
Mais $D$ est la r\'eunion des $\cal{D} \otimes_{\cE_\cA^{\dag,n}} \Ro_{A,p^{-1/n}}$ pour $n\geq N$. La proposition~\ref{propcle} assure que le $A[G_{\gamma_0}]$-module $\cal{D} \otimes_{\cE_\cA^{\dag,n}} \Ro_{A,p^{-1/n}}$ s'\'etend en un $\Ro_{A,p^{-1/n}}$-module via $G_{\gamma_0} \mapsto T$. Comme $T$ est inversible dans $\Ro_{A,p^{-1/n}}$ si $n>0$ il vient que $G_{\gamma_0}$ est inversible sur $\cal{D} \otimes_{\cE_\cA^{\dag,n}} \Ro_{A,p^{-1/n}}$. Pour conclure, il reste \`a remarquer deux choses. Premi\`erement, si $u \in \Z_p^\ast$ alors $\gamma_0^u-1$ agit sur $\cal{D} \otimes_{\cE_\cA^{\dag,n}} \Ro_{A,p^{-1/n}}$ via l'\'el\'ement  $u(G_{\gamma_0}) \in \Ro_{A,p^{-1/n}}^\ast$, qui est aussi inversible. Cela montre que si $\gamma' \in 1+p^M\Z_p^\ast$ alors $\gamma'-1$ est bijectif sur $(1+T)\varphi^M(D)$. Deuxi\`emement, cela vaut encore si on remplace $(1+T)$ par $(1+T)^a$ pour $a \in \Z_p^\ast$. En effet, pour tout $a \in \Z_p^\ast$ et $\gamma' \in 1+p^M\Z_p^\ast$ l'action de $a$ sur $D$ induit un isomorphisme $\gamma'$-\'equivariant $$(1+T)\varphi^M(D) \isomo (1+T)^a\varphi^M(D).$$
On conclut car $D^{\psi=0}=\bigoplus_{1\leq i \leq (p-1)p^{M-1}, (p,i)=1} (1+T)^i\varphi^M(D)$.
\end{pf}

\begin{cor}\label{cortri1} Si $D$ est triangulin sur $\Ro_A$ alors $\gamma-1$ est bijectif sur $D^{\psi=0}$.
\end{cor}

\begin{pf}  Il suffit de v\'erifier que $\Ro_A(\delta)$ est $\Gamma$-born\'e si $\delta \in \cT(A)$. Soit $e$ une base de $\Ro_A(\delta)$ telle que $\gamma(e)=\delta(\gamma)e$ pour tout $\gamma \in \Gamma$. Comme $\Gamma$ est compact et $\delta$ est continu, les \'el\'ements $\delta(\gamma)$ et $\delta(\gamma)^{-1}$ de $A$ sont \`a puissances positives born\'ees pour tout $\gamma \in \Gamma$. Comme $\Gamma$ est topologiquement de type fini, on peut donc trouver un mod\`ele $\cA \subset A$ tel que $\delta(\Gamma) \subset \cA^\ast$. \`A fortiori, $\Gamma.e \subset \cA[[T]]e$, ce qui conclut.
\end{pf}

\newcommand{\Res}{{\mathrm{Res}}}
\newcommand{\LA}{{\mathrm{LA}}}

\subsection{Cohomologie de $\Ro_A(\delta)$ partie I : calcul de $H^0$ et $H^2$} \label{cohoRoA} Nous allons maintenant calculer les $H^i(\Ro_A(\delta))$. Rappelons que lorsque $A$ est un corps, ce calcul est d\^u \`a Colmez pour $i=0,1$ \cite{colmeztri}, et Liu pour $i=2$ \cite{liu}. Dans le cas g\'en\'eral, nous allons proc\'eder de mani\`ere l\'eg\`erement diff\'erente \`a celle de Colmez en utilisant un d\'evissage de l'anneau de Robba (aussi d\^u \`a Colmez) que nous rappelons maintenant. Si $f=\sum_{n\in \Z} a_n T^n \in
\Ro_A$, on note $\Res(f) \in A$ le r\'esidu en $0$ de la forme
diff\'erentielle $ f(T) \frac{dT}{1+T}$, c'est \`a dire l'\'el\'ement $a_{-1}$ dans l'\'ecriture $\frac{f(T)}{1+T}=\sum_{n \in \Z}a_n T^n$. 

\begin{definition}\label{transfocolmez} Si $f \in \Ro_{A}$, la transform\'ee de Colmez de $f$ est la fonction $C(f) : \Z_p \rightarrow A$ d\'efinie par la formule
$$C(f)(x)={\rm Res}(f(1+T)^x)\,\,\,\,
\forall x \in \Z_p.$$ 
\end{definition}
Si $h\geq 0$ est un entier, notons $\LA_h(\Z_p,A)$ le $A$-module des fonctions $A$-valu\'ees et 
$h$-analytiques sur $\Z_p$, i.e. telles que pour tout $x \in \Z_p$, la
fonction $f_{x,h}(t):=f(x+p^ht)$ est dans $A\langle t \rangle $. 
C'est un espace de Banach pour la norme $|f|_h={\rm sup}_{x \in \Z_p}|f_{x,h}|$
o\`u $A\langle t\rangle$ est muni de la norme du sup. des coefficients. On a
de plus $|af|_h\leq |a| |f|_h$ si $a \in A$ et $f \in \LA_h(\Z_p,A)$. On munit $\LA_h(\Z_p)$ d'actions de $\psi$ et $\Gamma$ par les formules ($f \in \LA_h(\Z_p,A)$) $$\forall \gamma \in \Gamma \, \, \, \, \,  \gamma(f)(x)=\gamma f(\gamma^{-1}x), \, \, \psi(f)(x)=f(px).$$
On prendra garde que l'action de $\Gamma$ d\'efinie ci-dessus n'est pas
l'action na\"ive. De plus, $\LA(\Z_p,A)=\cup_{h\geq 0} \LA_h(\Z_p,A)$ est muni d'une action de $\varphi$ si l'on pose $\varphi(f)(x)=0$ si $x \in \Z_p^\ast$, $\varphi(f)(x)=f(x/p)$ si $x \in p\Z_p$.
\ps

\begin{prop}\label{transfcolmez} $C$ induit une suite exacte commutant aux
actions de $\varphi,\psi$ et $\Gamma$ : 
$$ 0 \longrightarrow \Ro_A^+ \longrightarrow \Ro_A \overset{C}{\longrightarrow} \LA(\Z_p,A) \longrightarrow
0.$$
\end{prop}

\begin{pf}  Quand $A$ est un corps, c'est le th\'eor\`eme I.1.3 de
\cite{colmezgros}, l'argument est similaire en g\'en\'eral. En effet, si $x
\in \Z_p$ on a par d\'efinition $(1+T)^x=\sum_{n\geq 0} {{x}\choose{n}} T^n
\in \Z_p[[T]]$. L'application $f \mapsto {\rm Res}(f)$ \'etant clairement
continue sur chaque $\cE^I_A$, on a la formule $$C(f)(x+1)=\sum_{n\geq 0}
a_{-1-n} {{x}\choose{n}} \in A$$ o\`u $f=\sum_{n \in \Z} a_n T^n$.
Appliquant ceci \`a $x=0,1,2,\dots$, il vient que ${\rm Ker}\,\,C(f)=\Ro_A^+$.
De plus, on obtient que ${\rm Im}\,\,C(f)$ est exactement l'ensemble des
fonctions continues $\Z_p
\rightarrow A$ de la forme $\sum_{n\geq 0} c_n {{x}\choose{n}}$ telles qu'il
existe $r > 1$ tel que $|c_n| r^n \rightarrow 0$ quand $n\rightarrow
\infty$. D'autre part, un th\'eor\`eme d'Amice assure que $\LA_h(\Z_p,\Q_p)$ pour
$h\geq 0$ est un
espace de Banach sur $\Q_p$ ayant pour base orthonorm\'ee
$[\frac{n}{p^h}]!{{x}\choose{n}}$. Via l'isom\'etrie naturelle $A\langle t \rangle = \Q_p \langle t
\rangle \widehat{\otimes}_{\Qp} A$, on dispose d'un isomorphisme naturel
$\LA_h(\Z_p,\Q_p)\widehat{\otimes}_{\Q_p} A \isomo \LA_h(\Z_p,A)$, de sorte
que $[\frac{n}{p^h}]!{{x}\choose{n}}$ est aussi une base orthonorm\'ee du
$A$-module de Banach $\LA_h(\Z_p,A)$. Pour conclure la surjectivit\'e de $C$, il suffit donc de voir
que pour une suite $(c_n) \in A^\N$,
il y a \'equivalence entre satisfaire $|c_n|r^n \rightarrow 0$ pour un
certain $r>1$
et satisfaire $\frac{|c_n|}{|[\frac{n}{p^h}]!|} \rightarrow 0$ pour un
entier $h$ assez grand. Mais $v_p([\frac{n}{p^h}]!)=\frac{n}{p^h(p-1)}+O({\rm
log}(n))$ et $v_p([\frac{n}{p^h}]!)\geq \frac{n}{p^h(p-1)}$, donc si on pose $r_h=p^{-\frac{1}{p^h(p-1)}}$, alors
pour $h\geq 0$ fix\'e et tout $n$ assez grand on a $r_h^n \leq |[\frac{n}{p^h}]!| \leq
r_{h+1}^n$. \ps
Il ne reste qu'\`a voir que $C$ est \'equivariant pour $\varphi,\psi$ et
$\Gamma$. Remarquons pour cela que les estim\'ees ci-dessus montrent que
pour tout entier $h\geq 0$ on a $C(\Ro_{A,r_h}) \subset \LA_h(\Z_p,A)$ et
$$C_{|\Ro_{A,r_h}} : \Ro_{A,r_h} \rightarrow \LA_h(\Z_p,A)$$ est continue.
La commutation \`a $\varphi,\psi$ et $\Gamma$ se v\'erifie alors sur la
partie dense $A.\Ro_{\Q_p,r_h}$, soit encore dans le
cas $A=\Q_p$, o\`u elle est d\'emontr\'ee dans~\cite[Prop. I.2.2]{colmezgros}.

\end{pf}

Notons $x \in \LA_0(\Z_p,\Z_p)$ la fonction identit\'e. 
Si $N\geq 0$, notons ${\rm Pol}_{\leq N}(\Z_p,A) \subset \LA_0(\Z_p,A)$ le sous-$A$-module (libre de rang $N+1$) des fonctions polynomiales de degr\'e $\leq N$, et ${\rm
Pol}(\Z_p,A)=A[x]=\cup_N {\rm Pol}_{\leq N}(\Z_p,A)$. On a pour tout entier $h\geq 0$ (et pour $h=\emptyset$) $${\rm Pol}_{\leq N}(\Z_p,A)\oplus
x^{N+1}\LA_h(\Z_p,A) = \LA_h(\Z_p,A),$$ les deux sous-espaces \'etant stables par $\psi$ et $\Gamma$. 
La $A$-base des mon\^omes $x^i$ de ${\rm Pol}_{\leq N}(\Z_p,A)$ est propre pour $\psi$ et $\Gamma$ : 
pour $i \leq N$ on a $\psi(x^i)=p^ix^i$ et $\gamma(x^i)=\gamma^{1-i} x^i$ pour tout $\Gamma$.\ps

On pose enfin $t=\log(1+T)=\sum_{k\geq 1} (-1)^{k+1} \frac{T^k}{k} \in \Ro^+$. On a $\varphi(t)=pt$, donc $\psi(t)=p^{-1} t$, et $\gamma(t)=\gamma t$ pour tout $\gamma \in \Gamma$.

\begin{lemma}\label{mainlemme} Soient $\lambda \in A^\ast$ et $N\geq 0$ un entier.\ps \begin{itemize}

\item[(i)] 0n a une d\'ecomposition $\varphi$-stable $\Ro_A^+=( \bigoplus_{0\leq i < N} A t^i ) \oplus T^N\Ro_A^+$. \ps
\item[(ii)] $1-\lambda\varphi$ est injectif sur $\LA(\Z_p,A)$, 
$\Ro_A^{\lambda\varphi=1}=(\Ro_A^{+})^{\lambda\varphi=1}$, et si 
$|\lambda p^N|<1$ alors $1-\lambda\varphi$ est bijectif sur $T^N\Ro_A^+$, d'inverse continu.\ps
\item[(iii)] Si $|p^{N+1}\lambda|<1$, alors $1-\lambda \psi$ est bijectif sur
$x^{N+1}\LA_h(\Z_p,A)$ pour tout $h\geq 0$.  
\ps
\item[(iv)] $\oplus_{i=0}^N A.T^{-(i+1)} \subset \Ro_{A}$ 
est un sous-$A$-module $\psi$-stable sur lequel $C$ induit un isomorphisme avec ${\rm Pol}_{\leq  N}(A,\Z_p)$. \ps
\item[(v)] $(1-\lambda \psi)\Ro_A^+=\Ro_A^+$ et la transform\'ee de Colmez induit une suite exacte 
$$ 0 \longrightarrow (\Ro_A^+)^{\lambda \psi=1} \longrightarrow \Ro_A^{\lambda\psi=1} \overset{C}{\longrightarrow} 
\LA(\Z_p,A)^{\lambda\psi=1} \longrightarrow 0$$ ainsi qu'un isomorphisme $$C: \Ro_A/(\lambda\psi-1) \isomo 
\LA(\Z_p,A)/(\lambda\psi-1).$$ 
\item[(vi)] Si $|\lambda p^N|<1$ alors $1-\lambda\varphi$ induit une bijection
$\Gamma$-\'equivariante $$(T^N\Ro_A^+)^{\lambda^{-1}\psi=1} \isomo (\Ro_A^+)^{\psi=0}\cap
T^N\Ro_A^+.$$ \ps
\end{itemize}
\end{lemma}

\begin{pf} Le (i) d\'ecoule de ce que $t \in T+T^2\Ro_A^+$ et $\varphi(T) \in
T\Ro_A^+$. V\'erifions le premier point du (ii). Tout d'abord, l'injectivit\'e de
$1-\lambda\varphi$ sur $\LA(\Z_p,A)$ entra\^ine que
$\Ro_A^{\lambda\varphi=1}=(\Ro_A^{+})^{\lambda\varphi=1}$ via la suite
exacte de la proposition~\ref{transfcolmez}. Soit donc $f \in \LA(\Z_p,A)$ telle que 
$\lambda \varphi(f)=f$.  Il vient que $\lambda^n \varphi^n(f)= f$ pour tout $n\geq 1$. 
En particulier, $f$ est nulle sur $p^{n-1}\Z_p^\ast$ pour tout $n\geq 1$ :
$f=0$, ce que l'on voulait d\'emontrer. Le second point du (ii) est l'argument de Colmez~\cite[Lemme A.1]{colmeztri}, que l'on rappelle par commodit\'e pour le lecteur. 
D'une part, pour tout $0<r<1$ et tout $f \in \Ro_A^+$, $|\varphi(f)|_{[0,r]}\leq |f|_{[0,r]}$
(se ramener \`a $A=\Q_p$ auquel cas cela d\'ecoule de l'interpretation de
$|.|_r$ comme norme sup. sur $B_{[0,r]}$). D'autre part, pour tout $0<r<1$ 
il existe $C_r>0$ tel que pour tout $i\geq 0$, $|\varphi^i(T^N)|_{[0,r]}\leq \frac{C_r}{p^{Ni}}$
(idem). Ainsi, si $|\lambda p^N|<1$, alors pour tout $f \in T^N\Ro_A^+$,  on a 
$|\lambda^k \varphi^k(f)|_{[0,r]} \leq C_r |f|_{[0,r]}|\lambda p^N|^k$ et
$$|\sum_{k\geq 0}\lambda^k \varphi^k(f)|_{[0,r]}\leq \frac{C_r |f|_{[0,r]}}{1-|\lambda p^N|}.$$ 
On a donc construit un inverse continu de $1-\lambda \varphi$ sur $T^N\Ro_A^+$. \ps

 Montrons le (iii). Remarquons que si $f \in \LA_h(\Z_p,A)$ et $h\geq 1$ alors 
$\psi(f) \in \LA_{h-1}(\Z_p,A)$. Ainsi, $\psi(x^i\LA_{h}(\Z_p,A)) \subset
x^i\LA_{h-1}(\Z_p,A)$ pour tout $i\geq 0$. Enfin, si $f \in
x^{N+1}\LA_{0}(\Z_p,A)$ alors $|\psi(f)|_0\leq \frac{|f|_0}{p^{N+1}}$. Ainsi, pour tout $m\geq h$ on a 
$$|\psi^m(f)|_0\leq \frac{1}{p^{(N+1)(m-h)}} |\psi^h(f)|_0.$$ Sous l'hypoth\`ese sur $N$, la
s\'erie $\sum_{m\geq 0} \lambda^m\psi^m$ converge donc dans les
endomorphismes de $x^{N+1}\LA_{h}(\Z_p,A)$, vers un inverse continu de 
${\rm id}-\lambda\psi$.\ps
Pour le (iv), un calcul  sans difficult\'e montre que $\psi(T^{-i-1})=Q(T)T^{-i-1}$ o\`u $Q(T) \in \Q[T]$ 
est un polyn\^ome de degr\'e $<i+1$ (et tel que $Q(0)=p^{-i}$). On conclut
car $C(T^{-i-1})(x+1)={{x}\choose{i}}$. \ps
Montrons maintenant le (v).  Si $N$ est suffisament grand de sorte que $|\lambda^{-1} p^N|<1$, on a vu au (ii) que  $\sum_{k\geq 0} \lambda^{-k} \varphi^k$ converge normalement sur $T^N\Ro_A^+$ vers un inverse de $1-\lambda^{-1}\varphi$. L'op\'erateur $\psi$ \'etant continu
sur $\Ro_A^+$,  la relation formelle $(1-\lambda\psi)(\sum_{k\geq
0}\lambda^{-k}\varphi^k)=-\lambda  \psi$ a donc un sens sur $T^N\Ro_A^+$, puis $$(1-\lambda\psi) T^N \Ro_A^+=\psi(T^N\Ro_A^+) \supset \psi(\varphi(T^N)\Ro_A^+)=T^N\psi(\Ro_A^+)=T^N\Ro_A^+.$$ 
Enfin on a $\psi(T)=-1$, donc si $i \geq 0$ on a l'identit\'e
$$(\lambda\psi-1)(T\varphi(T^i))=-\lambda T^i-T\varphi(T)^i.$$ Comme $T\varphi(T)^i \in p^iT^{i+1}+T^{i+2}\Ro_A^+$ une 
r\'ecurrence descendante montre que $T^i \in (\lambda\psi-1)\Ro_A^+$ pour tout $i\leq N$. Le (v) suit en appliquant 
$\lambda \psi = 1$ \`a la suite exacte de la transform\'ee de Colmez, et d'apr\`es le (i).\ps

La derni\`ere assertion d\'ecoule de (i) et (ii), et de ce que si $(1-\lambda \varphi)x=y$ alors $\lambda(\lambda^{-1}\psi-1)x=\psi(y)$, donc $\psi(y)=0$ si et seulement si $(1-\lambda^{-1}\psi)x=0$.
\end{pf}


Une cons\'equence de (i), (ii), (iii) et (v) du lemme ci-dessus est la proposition suivante, qui donne une description compl\`ete de $H^i(\Ro_A(\delta))$ pour $i=0$ et $i=2$. \ps

\begin{prop}\label{corH0H2} Soit $\delta \in \cT(A)$. On a 
$$\Ro_A(\delta)^{\varphi=1}=\Ro_A^{\delta(p)\varphi=1}=A[t]^{\delta(p)\varphi=1}$$
En particulier, $H^0(\Ro_A(\delta))=A[t]^{\delta(p)\varphi=1, \delta(\gamma)\gamma=1}$. De plus, 
la transform\'ee de Colmez induit un isomorphisme $\Gamma$-\'equivariant 
$$\Ro_A(\delta)/(\psi-1)=\Ro_A/(\delta(p)^{-1}\psi-1) \isomo {\rm Pol}(\Z_p,A) /(\delta(p)^{-1}\psi-1).$$
En particulier, $H^2(\Ro_A(\delta))={\rm Pol}(\Z_p,A) /(\delta(p)^{-1}\psi-1,\delta(\gamma)\gamma-1)$.
\end{prop}

Suivant Colmez, d\'esignons par $x : \Qp^\ast \rightarrow \Qp^\ast$ le caract\`ere identit\'e et par $\chi :\Qp^\ast \rightarrow \Z_p^\ast$ le caract\`ere cyclotomique, c'est \`a dire $\chi=x|x|$. Un corollaire imm\'ediat de la proposition ci-dessus est le suivant.\ps

\begin{cor}\label{corsuiteH0H2} Soit $\delta \in \cT(A)$. \begin{itemize}
\item[(i)] $H^0(\Ro_A(\delta)) \neq 0$ si et seulement si il existe $i \geq 0$ et $f \neq 0 \in A$ tels que $f\cdot(\delta-x^{-i})$ est identiquement nul sur $\Q_p^\ast$.\ps
\item[(ii)] $H^2(\Ro_A(\delta)) \neq 0$ si et seulement si il existe $i\geq 0$ tel que le ferm\'e de ${\rm Sp}(A)$ sur lequel $\delta = \chi(x)x^i$ soit non vide.\ps
\end{itemize}
\end{cor}

Il reste \`a d\'ecrire $H^1(\Ro_A(\delta))$. D'apr\`es le d\'evissage de la cohomologie d\'emontr\'e plus haut il nous faut nous int\'eresser au $A[\Gamma]$-module $\Ro_A(\delta)^{\psi=1}$. Une \'etape clef sera alors la structure du $A[\Gamma]$-module $\Ro_A^+(\delta)^{\psi=0}$.  \ps

\subsection{Structures sur $\Ro_A^+(\Gamma)$} Soit $C$ un groupe profini isomorphe \`a $\Z_p$ et $c$ 
un g\'en\'erateur topologique de $C$. Si $\mathcal{A}$ est complet s\'epar\'e pour la topologie $p$-adique, 
et en particulier si c'est un mod\`ele d'une alg\`ebre affino\"ide $A$, il est connu depuis Iwasawa que 
l'application $$\cA[[T]] \rightarrow \cA[[C]]:=\projlim_n \cA[C/p^nC]$$ envoyant $T$ sur $[c]-1$ est un 
isomorphisme. Cela permet de d\'efinir un anneau $\Ro_A^+(C)$ en rempla\c{c}ant simplement la variable $T$ 
dans la d\'efinition de $\Ro_A^+$ par $[c]-1$. Cette d\'efinition ne d\'epend pas du choix de $c$ car $\Gamma$ agit par automorphismes sur $\Ro_A^+$. On remarque de plus que $\Ro_A^+(pC)=\varphi(\Ro_A(C))$, et donc que $\Ro_A^+(pC) \otimes_{A[pC]}A[C]=\Ro_A^+(C)$, car les $(1+T)^i$ pour $0\leq i \leq p-1$ forment une base $\Ro_A^+$ sur $\varphi(\Ro_A)^+$. On pose enfin $$\Ro_A^+(\Gamma):=\Ro_A^+(1+p^M\Z_p) \otimes_{A[1+p^M\Z_p]}A[\Z_p^*]$$
qui ne d\'epend par du choix de $M\geq 1$ par ce que l'on vient de dire. On dispose d'inclusions naturelles denses $$A[\Gamma] \subset (\cA[[\Gamma]])[1/p] \rightarrow \Ro_A^+(\Gamma).$$
On posera aussi $A[[\Gamma]]_b=(\cA[[T]])[1/p] \subset \Ro_A^+(\Gamma)$. Il ne d\'epend pas du choix de $\cA$. 
Quand $A$ est un corps le r\'esultat suivant est essentiellement d\^u \`a Berger (\cite[\S
5]{berger1}).

\begin{prop} Soit $D$ un $\fg$-module sur $\Ro_A$ ou $\Ro_A^+$. L'action $A$-lin\'eaire de $\Gamma$ sur $D$ s'\'etend 
de mani\`ere unique en une structure de $\Ro_A^+(\Gamma)$-module sur $D$ qui est continue au sens suivant : si $(e_i)$ 
est une $\Ro_A$-base de $D$ telle que $\Gamma(D_r) \subset D_r$ o\`u $D_r=\oplus_i \Ro_{A,r}e_i$, alors pour tout $f \in D_r$  l'application orbite  $\Ro_A^+(\Gamma) \rightarrow D_r=\Ro_{A,r}^d$, $u \mapsto u(f)$, est continue. Cette action de $\Ro_A^+(\Gamma)$ commute \`a $\varphi$ et $\psi$. \ps
\end{prop}

\begin{pf} En effet, il suffit de voir que si $I$ est un intervalle ferm\'e de $[0,1[$ et si $D=(\cE^I_A)^n$ est muni d'une action semi-lin\'eaire de $\Gamma$ qui soit continue 
dans le sens que la matrice $M_\gamma \in M_n(\cE^I_A)$ de $\gamma \in \Gamma$ dans la base canonique
$(e_i)$ d\'epende contin\^ument de $\gamma$, alors l'application naturelle $$A[\Gamma] \rightarrow \End_A(D)$$
se prolonge en une application continue $\Ro_A^+(\Gamma) \rightarrow \End_A(D)$ (un tel prolongement \'etant n\'ecessairement unique s'il existe). 
Fixons un tel $I$ ainsi qu'un mod\`ele $\cA \subset A$. Le sous-espace $\cal{D}=(\OO^I_\cA)^n
\subset D$ est un ouvert (pour la topologie de module de Banach sur
$\cE_A^I$ de ce dernier), ainsi que $p\cal{D}$, de sorte qu'il existe un entier $M\geq 1$ tel
que pour tout $\gamma \in 1+p^M\Z_p$ et tout entier $i$ on ait $(M_\gamma-1)(e_i) \in
pM_n(\OO^I_\cA)$. Le lemme~\ref{estimeegamma} (i) permet de supposer de plus que
$(\gamma-1)\cE_\cA^I \subset p\cE_\cA^I$ pour tout $\gamma \in 1+p^M\Z_p$.
La relation $$(\gamma-1)(\sum_i a_i e_i)=\sum_i
(\gamma-1)(a_i)\gamma(e_i)+a_i(\gamma-1)(e_i)$$
assurent alors que $(\gamma-1)\cal{D} \subset p\cal{D}$ pour tout $\gamma
\in 1+p^M\Z_p$. On en d\'eduit que le morphisme $A[1+p^M\Z_p] \rightarrow \End_A(D)$ s'\'etend contin\^ument 
en un morphisme $$\Ro_A^+(1+p^M\Z_p) \rightarrow \End_A(D)$$
ainsi donc qu'\`a $\Ro_A^+(\Gamma)=\Ro_A^+(1+p^M\Z_p) \otimes_{A[1+p^M\Z_p]} A[\Z_p^\times]$. Remarquons que 
jusqu'ici nous n'avons pas utilise la structure de $\varphi$-module sur $D$. La commutation de $\Ro_A^+(\Gamma)$ 
\`a $\varphi$ et $\psi$ vient de leur commutation \`a $\Gamma$ et de ce que pour tout $r$ assez grand, $\varphi : D_{r^p}
\rightarrow D_r$ et $\psi : D_r \rightarrow D_r$ sont continues (ceci
d\'ecoulant du cas $D=\Ro_A$ consid\'er\'e au~\S~\ref{gencoho}). \end{pf}

Un point crucial est que la structure de $D^{\psi=0}$ comme $\Ro_A^+(\Gamma)$-module est particuli\`erement simple.

\begin{prop}\label{propRplus} Si $D$ est un $\fg$-module sur $\Ro^+_A$ de rang $d$ qui est $\Gamma$-born\'e alors $D^{\psi=0}$ est libre de rang $d$ sur $\Ro_A^+(\Gamma)$. \ps

Plus pr\'ecis\'ement, si $(e_i)$ est une $\Ro_A^+$-base de $D$ telle que $\oplus \cA[[T]] e_i$ est stable par $\Gamma$ pour un certain mod\`ele $\cA \subset A$, alors  $(1+T)\varphi(e_i)$ est une base de $D^{\psi=0}$ sur $\Ro_A^+(\Gamma)$. 
\end{prop} 

Avant de d\'emontrer ce r\'esultat nous devons dire un mot sur la topologie consid\'er\'ee sur $D^{\psi=0}$. Tout d'abord, le choix d'une base $(e_i)$ fournit une \'ecriture $D=\oplus_i \Ro_A^+ e_i$ et donc une structure d'espace de Fr\'echet sur $D$. Cette structure ne d\'epend pas du choix de la base et l'application structurale $\Ro^+_A \times D \rightarrow D$ est continue. Dans la base $\varphi(e_i)$, on a $\psi(\sum_i x_i \varphi(e_i))=\sum_i \psi(x_i)e_i$ donc la continuit\'e de $\psi$ sur $\Ro_A^+$ entra\^ine sa continuit\'e sur $D$; celle de $\varphi$ est imm\'ediate. En particulier, $D^{\psi=0}$ est ferm\'e, ainsi que $(1+T)\varphi^M(D)$ pour tout $M\geq 1$, et l'application $x \mapsto (1+T)\varphi^M(x)$ est un hom\'eomorphisme de $D$ sur $(1+T)\varphi^M(D)$ d'inverse $y \mapsto \psi^M(\frac{y}{1+T})$.

\begin{pf} Soit $(e_i)$ comme dans l'\'enonc\'e, $\cal{D}=\oplus_i \cA[[T]]e_i$ et $M$ associ\'e \`a $\cal{D}$ comme dans la proposition~\ref{propcle}. Soit $\gamma_0=1+p^M \in \Gamma$ et $x \in D$. Comme on l'a d\'ej\`a vu, on a la relation
$$(\gamma_0-1)((1+T)\varphi^M(x))=(1+T)\varphi^M(G_\gamma(x)).$$ 
La proposition~\ref{propcle} assure donc que l'action de $1+p^M\Z_p$ sur $(1+T)\varphi^M(D)$ s'\'etend en une structure de $\Ro_A^+(1+p^M\Z_p)$-module, qui est de plus libre de base $(1+T)\varphi^M(e_i)$ d'apr\`es la proposition~\ref{propcle}, et telle que $\Ro_A^+(1+p^M\Z_p) \times (1+T)\varphi^M(D) \rightarrow (1+T)\varphi^M(D)$ soit continue. Cette structure est n\'ecessairement celle de la proposition pr\'ec\'edente par unicit\'e de cette derni\`ere. Les identit\'es
$$D^{\psi=0}=\bigoplus_{1\leq i \leq p^{M-1}(p-1), (p,i)=1}(1+T)^i\varphi^M(D), $$
$\Ro_A^+(\Gamma)=\Ro_A^+(1+p^M\Z_p)\otimes_{A[1+p^M\Z_p]}A[\Gamma]$, et pour $a \in \Z_p^\ast$, $a((1+T)\varphi^M(D))=(1+T)^a\varphi^M(D)$, concluent la d\'emonstration. Le dernier point vient de ce que $$\Ro_A^+(\Gamma)/(\gamma-1)=\Ro_A^+(1+p\Z_p)/(\gamma^{p-1}-1)=\Ro_A^+/(T)=A.$$
\end{pf}

\begin{remark}{\rm Dans le m\^eme genre, on d\'eduirait ais\'ement de la proposition~\ref{propcle} que si $D$ est un $\fg$-module sur $\Ro_A$ qui est $\Gamma$-born\'e alors l'action de $\Gamma$ sur $D^{\psi=0}$ s'\'etend en une structure de $\Ro_A(\Gamma)$-module (voir \cite[V \S 3]{colmezgros} pour la d\'efinition), qui est libre 
sur $\Ro_A(\Gamma)$ de rang le rang de $D$ sur $\Ro_A$.}
\end{remark}

\subsection{Cohomologie de $\Ro_A(\delta)$, partie II : 
structure de $\Ro_A(\delta)^{\psi=1}$.} Retournons au 
calcul de $H^1(\Ro_A(\delta))$. \'Etant donn\'e le d\'evissage 
donn\'e par la proposition~\ref{devissagecoho}, il convient 
d'\'etudier tout d'abord $\Ro_A(\delta)^{\psi=1}$. Si $X$ est 
un $\fg$-module sur $\Ro_A$ ou $\Ro_A^+$, on pose suivant Fontaine $C(X)=X^{\psi=0}\cap
(1-\varphi)X=(1-\varphi)X^{\psi=1}$, de sorte que l'on dispose d'une suite
$\Ro_A^+(\Gamma)$-\'equivariante tautologique $$ 0 \longrightarrow
X^{\varphi=1} \longrightarrow X^{\psi=1} \overset{1-\varphi}{\longrightarrow} C(X)
\longrightarrow 0.$$

\begin{lemma}\label{devissedplus} Soient $\delta \in \cT(A)$, $D=\Ro_A(\delta)$ et $D^+=\Ro_A^+(\delta)$. On a des suites exactes naturelles $\Ro_A^+(\Gamma)$-\'equivariantes :\begin{itemize}
\item[(i)] $0 \longrightarrow (D^+)^{\psi=1} \longrightarrow D^{\psi=1} \overset{C}{\longrightarrow} {\rm
Pol}(\Z_p,A)^{\delta(p)^{-1}\psi=1} \rightarrow 0$,\ps
\item[(ii)] $0 \longrightarrow C(D^+) \longrightarrow C(D)
\overset{(1-\delta(p)\varphi)^{-1} C}{\longrightarrow} 
{\rm Pol}(\Z_p,A)^{\delta(p)^{-1}\psi=1} \longrightarrow 0$.
\end{itemize}
De plus, $(A[t])^{\delta(p)\varphi=1}=(D^+)^{\varphi=1}=D^{\varphi=1}$ et
$D/(\psi-1)={\rm Pol}(\Z_p,A)/(\delta(p)^{-1}\psi-1)$.
\end{lemma}

\begin{pf} Pour le (i) on applique $\delta^{-1}(p)\psi=1$ \`a la suite exacte d\'efinie par la 
transform\'ee de Colmez et on note que $(D^+)/(\psi-1)=0$ et 
$\LA(\Z_p,A)^{\delta(p)^{-1}\psi=1}={\rm Pol}(\Z_p,A)^{\delta(p)^{-1}\psi=1}$ d'apr\`es le
Lemme~\ref{mainlemme} (iv) et (iii). Le (ii) d\'ecoule du (i) et de l'injectivit\'e
$1-\delta(p)\varphi$ sur $\LA(\Z_p,A)$
(lemme~\ref{mainlemme} (ii)). En effet, cette injectivit\'e entra\^ine d'une
part que $C(D^+)=C(D)\cap D^+$, puis que la derni\`ere fl\`eche de l'\'enonc\'e est bien
d\'efinie ; elle est surjective par le (i). La derni\`ere assertion a d\'ej\`a \'et\'e d\'emontr\'ee
(corollaire \ref{corH0H2}).
\end{pf}

La structure de $C(D)^+$ s'av\`ere int\'eressante, avant de la d\'ecrire nous avons besoin d'un lemme sur $\Ro_A^+(\Gamma)$. 

\newcommand{\deltat}{\widetilde{\delta}}
\begin{lemma} \label{lemmeker} \begin{itemize}
\item[(i)] Si $\delta : \Z_p^\ast \rightarrow A^\ast$ est un caract\`ere continu, 
il s'\'etend de mani\`ere unique en un morphisme de $A$-alg\`ebres continu 
$\deltat: \Ro_A^+(\Gamma) \rightarrow A$. Tout morphisme de $A$-alg\`ebres continu $\Ro_A^+(\Gamma) \rightarrow A$ est 
de cette forme.\ps 
\item[(ii)] Soit $\gamma_0 \in \Gamma$ un g\'en\'erateur topologique de 
$\Gamma$ (resp. de $1+4\Z_2$ si $p=2$), soit $T_\delta:=[\gamma_0]-\delta(\gamma_0) \in \Ro_A^+(\Gamma)$. 
On a ${\rm Ker}(\deltat)=(T_\delta)$ o\`u $(T_\delta,[-1]-\delta(-1))$ selon 
que $p>2$ ou non. \ps
\item[(iii)] La multiplication par $T_\delta$ est injective sur $\Ro_A^+(\Gamma)$ et $\Ro_A^+(\Gamma)=T_\delta\Ro_A^+(\Gamma) \oplus A$
(resp. $\Ro_A^+(\Gamma)=T_\delta\Ro_A^+(\Gamma) \oplus A[\{\pm 1\}]$
si $p=2$). 
\end{itemize}
\end{lemma}

\begin{pf} Soit $\cA$ un mod\`ele de $A$ contenant $\delta(\Gamma)$. Pour $M\geq 1$ 
assez grand, on a $\delta(1+p^M\Z_p) \subset 1+p\cA$. Rappelons que $\Ro_A^+(1+p^M\Z_p)$ s'identifie \`a $\Ro_A^+$ si l'on envoit $[c]-1$ vers $T$, $c$ \'etant un g\'en\'erateur quelconque de $1+p^M\Z_p$. Il est alors imm\'ediat que 
$\delta: A[1+p^M\Z_p] \rightarrow A$ s'\'etend \`a $\Ro_A^+(1+p^M\Z_p)$, ainsi donc 
qu'\`a $\Ro_A^+(\Gamma)$ par extension des scalaires. La r\'eciproque d\'ecoule ais\'ement de ce que l'application 
naturelle $\Z_p[[\Gamma]] \rightarrow \Ro_A^+(\Gamma)$ est continue
(lemme~\ref{topo} (ii)) : cela d\'emontre le (i). \par 
Pour le (ii), il est clair que $T_\delta$ et $[-1]-\delta(-1)$ sont 
dans ${\rm Ker}(\deltat)$, et que $[\gamma_0^n]-\delta(\gamma_0)^n \in 
(T_\delta)$ pour tout $n \in \Z$. Soit $n \in \Z$ tel que $c=\gamma_0^n$
engendre topologiquement $1+2p\Z_p$. $\Ro_A^+(1+2p\Z_p)$ s'identifie 
\`a l'anneau de Robba $A$-valu\'e positif sur la variable $U=[c]-1$, et il vient que 
$(T_\delta) \supset (U-(\delta(c)-1))\Ro_A^+(1+2p\Z_p)$. Comme
$(\delta(c)-1)^m$ tend vers $0$ quand $m$ tend vers l'infini, on a
$\Ro_A^+(1+2p\Z_p)=A \oplus ([c]-1)-(\delta(c)-1))\Ro_A^+(1+2p\Z_p)$. Le
(ii) suit car l'application canonique $A[\Gamma_{\rm tors}]
\rightarrow \Ro_A^+(\Gamma)/\Ro_A^+(1+2p\Z_p)=A[\Gamma/(1+2p\Z_p)]$ est un
isomorphisme. 
\par 
Le premier point du (iii) suit de l'analyse ci-dessus. En effet, 
$\Ro_A^+(\Gamma)$ est libre de rang fini sur $\Ro_A^+(1+2p\Z_p)$ et 
la multiplication par $U-(\delta(c)-1) \in (T_\delta)$ est injective 
sur l'anneau de Robba positif en $U$. Le second point d\'ecoule du (ii).  \end{pf}

Consid\'erons pour tout $k\geq 0$ l'application $$J_k : \Ro_A^+ \rightarrow \Ro_A^+/T^k\Ro_A^+=A[T]/T^k.$$
Notons que $T^k\Ro_A^+$ est ferm\'e dans $\Ro_A^+$ et qu'il est stable par $\Gamma$ (et donc 
$\Ro_A^+(\Gamma)$) et $\varphi$. Ainsi, $J_k$ est \'equivariante sous $\Ro_A^+(\Gamma)$ et 
$\varphi$. De plus, les images de $1,t,\dots,t^{k-1}$ forment une $A$-base de 
$\Ro_A^+/T^k\Ro_A^+$ propre pour $\varphi$ et $\Gamma$ de valeurs propres \'evidentes.

\begin{lemma}\label{surjjk} $J_k$ induit une surjection $(\Ro_A^+)^{\psi=0} \rightarrow
\Ro_A^+/T^k\Ro_A^+$.
\end{lemma}

\begin{pf} En effet, l'\'ecriture formelle "$1+T={\rm exp}(t)$" assure que $J_k(1+T)=\sum_{i=0}^{k-1} \frac{t^i}{i!}$, l'important pour ce qui suit \'etant que le coefficient de chaque $t^i$ est non nul. Comme pour $i=0,\dots,k-1$ les
caract\`eres $\gamma \mapsto \gamma^i$, $\Gamma \rightarrow \Q_p^\ast$, sont
lin\'eairement ind\'ependants sur $\Q_p$ (car distincts), il vient que
$J_k(A[\Gamma](1+T))=\Ro_A^+/T^k\Ro_A^+$. Cela conclut car $1+T \in
(\Ro_A^+)^{\psi=0}$. 
\end{pf}

\begin{prop}\label{struccplus} Le $\Ro_A^+(\Gamma)$-module $C(D^+)$ est isomorphe \`a
l'id\'eal $$\cap_{i\geq 0}(1-\delta(p)p^i,\Ker(\widetilde{x^i\delta}))$$
de $\Ro_A^+(\Gamma)$. Plus pr\'ecis\'ement, pour tout $k$ assez grand $J_k$ induit
une suite exacte $\Ro_A^+(\Gamma)$-\'equivariante $$0 \rightarrow C(D^+)
\longrightarrow \Ro_A^+(\Gamma)\cdot (1+T) \overset{J_k}{\longrightarrow} \oplus_{i=0}^{k-1} 
A/(1-\delta(p)p^i)
\,\,\cdot \widetilde{\delta x^i} \rightarrow 0.$$ 
\end{prop}

(Comme $1-\delta(p)p^i$ est inversible dans $A$ pour $i$ assez grand,
l'intersection ci-dessus est finie. De plus, le terme centrale est libre de
rang $1$ sur $\Ro_A^+(\Gamma)$.)

\begin{pf} En effet, soit $k$ suffisament grand de sorte que
$|\delta(p)p^k|<1$. Le lemme~\ref{mainlemme} (ii) assure que
$$(1-\delta(p)\varphi)\Ro_A^+=J_k^{-1}((1-\delta(p)\varphi)\Ro_A^+/T^k\Ro_A^+).$$
Comme $\Ro_A^+(\Gamma)$-module on a
$$(1-\varphi)\Ro_A^+(\delta)/T^k\Ro_A^+(\delta)=\oplus_{i=0}^{k-1} (1-\delta(p)p^i)A(
\widetilde{\delta x^i}).$$
Mais $\Ro_A^+(\delta)^{\psi=0}$ est libre de rang $1$ sur $\Ro_A^+(\Gamma)$
engendr\'e par $1+T$ d'apr\`es la proposition~\ref{propRplus}. La
proposition suit alors du lemme~\ref{surjjk}.
\end{pf}

Fixons un g\'en\'erateur topologique $\gamma_0$ de $\Gamma$ (resp. de
$1+4\Z_2$ si $p=2$). On rappelle l'\'el\'ement $T_\delta \in 
A[\Gamma]$ d\'efini dans le lemme~\ref{lemmeker} (ii). Pour $i\in \Z$ on
note $T_i \in A[\Gamma]$ l'\'el\'ement $T_\delta$ o\`u
$\delta(\gamma)=\gamma^i$ pour tout $\gamma \in \Gamma$. 

\begin{definition}\label{defpresque} On dira qu'un
$\Ro_A^+(\Gamma)$-module $D$ est presque libre de rang $d$ si il existe
une suite exacte de $\Ro_A^+(\Gamma)$-modules de la
forme $$0 \rightarrow \Ro_A^+(\Gamma)^d \rightarrow D \rightarrow Q
\rightarrow 0$$
telle que $Q$ est de type fini sur $A$ et annul\'e par un mon\^ome en des
$T_i$ pour $i \in \Z$. On dira que $D$ est presque nul si il est presque
libre de rang $0$, c'est \`a dire de type fini sur $A$ et annul\'e par un
m\^onome en des $T_i$ pour $i \in \Z$.
\end{definition}

\begin{lemma}\label{lemmepl} 
\begin{itemize} \item[(i)] 
Les $\Ro_A^+(\Gamma)$-modules presque libres sont de type fini. \ps
\item[(ii)] Soit $0 \rightarrow D' \rightarrow D \rightarrow D'' \rightarrow
0$ une suite exacte de $\Ro_A^+(\Gamma)$-modules. Si $D''$ (resp. $D'$) est
presque nul et $D$ est presque libre de rang $d$, alors $D'$ (resp. $D''$) 
est presque libre de rang $d$. \ps
\item[(iii)] Enfin, si on a une suite exacte longue de
$\Ro_A^+(\Gamma)$-modules $$ D_1 \longrightarrow D_2 \longrightarrow D_3
\longrightarrow D_4 \longrightarrow D_5$$ avec $D_1$ et $D_5$ presque nuls,
$D_2$ et $D_4$ presque libres de rang respecifs $d_2$ et $d_4$, alors $D_3$
est presque libre de rang $d_2+d_4$. 
\end{itemize}
\end{lemma}

\begin{pf} Le premier point est \'evident. Pour le second, soient $L \subset
D$ libre de rang $d$ et $M$ un mon\^ome en les $T_i$ tel que $MD \subset
L$. Si $M'$ est un autre tel mon\^ome tel que $M'D \subset D'$ alors le
lemme~\ref{lemmeker} (ii) assure que $M'L$ est libre de rang $d$ et que
$L/M'L$ est de type fini sur $A$. Il vient que $D/M'L$ est de type
fini sur $A$ annul\'e par $MM'$, ainsi donc que son sous-module $D'/M'L$ par 
noeth\'erianit\'e de $A$, et donc $D'$
est presque libre de rang $d$ ce qui
prouve le (ii) dans le premier cas. Dans le second cas $D'$ est presque nul,
donc $D' \cap
L =0$ par le lemme~\ref{lemmeker} (ii). Ainsi, la projection $\pi : D \rightarrow
D''$ est injective sur $L$ et $D''/\pi(L)$ est un quotient du $A$-module
presque nul $D/L$, ce qui conclut le (ii).\par
Pour le (iii), on peut supposer $D_1=D_5=0$ par le (ii). Si $D_3$ est
libre sur $\Ro_A^+(\Gamma)$ alors la suite est scind\'ee et le r\'esultat
suit. En g\'en\'eral, quitte \`a
remplacer $D_2$ par l'image inverse dans $D_2$ d'un sous-
$\Ro_A^+(\Gamma)$-module libre de $D_3$ on peut donc supposer que $D_3$ est
presque nul, auquel cas l'affirmation est \'evidente. 
\end{pf}

\begin{thm}\label{mainiwa} Soit $D$ un $\fg$-module triangulin de rang $d$ sur $\Ro_A$. 
Alors les $\Ro_A^+(\Gamma)$-modules $D^{\varphi=1}$ et $D/(\psi-1)$ sont presque
nuls,
et $D^{\psi=1}$ et $C(D)$ sont presque libres de rang $d$. \ps
	Soit $(\delta_i) \in \cT^d(A)$ le param\`etre de $D$ et supposons de
plus que
pour tout $i=1,\dots,d$ et pour tout $j\in \N$ alors $1-\delta_i(p)p^j$ est non
diviseur de z\'ero dans $A$ et $1-\delta_i(p) p^{-j} \in A^\times$. Alors
$D^{\varphi=1}=D/(\psi-1)=0$ et $D^{\psi=1} \isomo C(D)$ est libre de rang $d$ sur
$\Ro_A^+(\Gamma)$.
\end{thm}

\begin{pf} On proc\`ede par r\'ecurrence sur $d\geq 1$. On a une suite
exacte dans $\fg/A$ $$ 0 \longrightarrow D' \rightarrow D \rightarrow \Ro_A(\delta_d)
\rightarrow 0$$
avec $D'$ triangulin de rang $d-1$. Par la suite exacte longue de
cohomologie associ\'ee et par le lemme~\ref{lemmepl} (iii), on peut supposer
$d=1$, i.e. $D=\Ro_A(\delta)$. Dans ce cas, il d\'ecoule de la
proposition~\ref{corH0H2} que $D^{\varphi=1}$ et $D/(\psi-1)$ sont presque
nuls car $1-\delta(p)^{\pm 1}p^i$ est inversible dans
$A$ pour tout entier $i$ assez grand, c'est aussi \'evidemment le cas de ${\rm
Pol}(A,\Z_p)^{\delta(p)^{-1}\psi=1}$. De plus, le premier de ces trois
modules est nul si et seulement si $1-\delta(p)p^i$ est non diviseur de $0$ dans $A$ pour
tout $i\geq 0$, et les deux autres le sont si et seulement si
$1-\delta(p)p^{-i} \in A^\times$ pour tout $i\geq 0$. Sous ces hypoth\`eses, le
lemme~\ref{devissedplus} assure aussi que $D^{\psi=1} \isomo C(D)=C(D^+)$ et la
proposition~\ref{struccplus} montre que
$C(D^+)=(D^+)^{\psi=0}=\Ro_A^+(\Gamma)\cdot (1+T)$ est libre de rang $1$.\end{pf}

Mentionnons qu'en g\'en\'eral, $C(D)$ n'est pas libre sur $\Ro_A^+(\Gamma)$,
et ce m\^eme si $D=\Ro_A(\delta)$. En effet, consid\'erons le cas 
particulier o\`u $\delta \in \cT(A)$ est tel que $1-\delta(p)p^{-i}$ est non diviseur de $0$ dans $A$ pour tout $i\geq
0$, auquel cas $C(D)=C(D^+)$. On conclut par le r\'esultat g\'en\'eral
suivant. 

\begin{prop}\label{critlibre} Soit $\delta \in \cT(A)$ tel que $1-\delta(p)p^i$ est non
diviseur de $0$ dans $A$ pour tout $i\geq 0$. Alors $C(\Ro_A^+(\Gamma))$ est
projectif comme $\Ro_A^+(\Gamma)$-module si et seulement si $1-\delta(p)p^i \in
A^\times$ pour tout $i\geq 0$, auquel cas il est en fait libre de rang $1$. \end{prop}

\begin{pf} Remarquons que si deux id\'eaux de type fini $I$ et $J$ d'un anneau
commutatif $B$ sont tels
que l'id\'eal $IJ$ est projectif comme $B$-module, et si de plus $I+J=B$,
alors $I$ et $J$ sont des $B$-modules projectifs. En effet, \^etre projectif
de type fini est une
propri\'et\'e locale sur ${\rm Spec}(B)$, mais si $x
\in {\rm Spec}(B)\backslash V(I)$ alors $I_x=B_x$, et si $x \in V(I)$ alors $x \notin
V(J)$ et donc $J_x=B_x$ puis $(IJ)_x=I_xJ_x=I_x$. D'apr\`es le
lemme~\ref{comax} ci-dessous et la proposition~\ref{struccplus}, il s'agit
de voir que si $I_\delta$ est projectif de rang fini, avec $1-\delta(p)$ non
diviseur de $0$ dans $A$, alors $1-\delta(p) \in
A^\times$ (et donc $I_\delta=\Ro_A^+(\Gamma)$). Fixons donc $\delta \in
\cT(A)$ avec $1-\delta(p)$ non diviseur de $0$ et supposons $p>2$ pour
simplifier, de sorte que ${\rm Ker}(\deltat)=(T_\delta)$. La
suite de $\Ro_A^+(\Gamma)$-modules
$$0 \longrightarrow \Ro_A^+(\Gamma) \overset{u \mapsto (T_\delta
u,(1-\delta(p))u)}{\longrightarrow}
\Ro_A^+(\Gamma)^2 \overset{(x,y) \mapsto (1-\delta(p))x-T_\delta
y}{\longrightarrow}   
I_{\delta} \longrightarrow 0$$
est exacte. En effet, si $x,y \in \Ro_A^+(\Gamma)$ satisfont $x T_\delta = y
(1-\delta(p))$ alors en appliquant $\deltat$ il vient que
$0=\deltat(y)(1-\delta(p)) \in A$ donc $y=T_\delta u$ avec $u \in
\Ro_A^+(\Gamma)$ unique par
le lemme~\ref{lemmeker}, puis $x=u (1-\delta(p))$. Comme la multiplication
par chaque $T_{\delta'}$ est injective
sur $I_{\delta} \subset \Ro_A^+(\Gamma)$ on en d\'eduit que la suite
ci-dessus est reste exacte modulo $T_{\delta'}$ pour tout $\delta'$.
L'isomorphisme naturel $\deltat' : \Ro_A^+(\Gamma)/(T_{\delta'}) \isomo
A$ envoie tout $\gamma \in \Gamma$ sur $\delta'(\gamma)$. Si $\gamma$ engendre
topologiquement $\Gamma$ on obtient donc une suite exacte 
$$0 \longrightarrow A \overset{u \mapsto ((\delta'(\gamma)-\delta(\gamma))u,
(1-\delta(p))u)}{\longrightarrow}
A \overset{(x,y) \mapsto (1-\delta(p))x-(\delta'(\gamma)-\delta(\gamma))y}{\longrightarrow}   
I_\delta/T_{\delta'}I_\delta \longrightarrow 0,$$ 
et donc $$I_\delta/T_{\delta'}I_\delta \simeq
A^2/A(1-\delta(p),\delta'(\gamma)-\delta(\gamma)).$$
On conclut en appliquant ceci \`a $\delta'=\delta$ : si le
$\Ro_A^+(\Gamma)$-module $I_\delta$ est
projectif il en va de m\^eme du $A$-module $A/(1-\delta(p)) \times A$. Comme $1-\delta(p)$ est non diviseur de $0$
dans $A$, cela
entra\^ine qu'il est inversible, ce qui conclut. L'argument est similaire pour $p=2$. \end{pf}

\begin{lemma}\label{comax} Pour $\delta \in \cT(A)$, notons $I_\delta$ l'id\'eal $(1-\delta(p),{\rm
Ker}(\deltat))$ de $\Ro_A^+(\Gamma)$. Alors pour tout $i \in \Z \backslash
\{0\}$ on a $I_\delta+I_{\delta x^i}=\Ro_A^+(\Gamma)$. En particulier, $C(\Ro_A^+(\delta))=\prod_{i\geq 0}
I_{\delta x^i}$.
\end{lemma}

\begin{pf} En effet, $(1-\delta(p)p^i)-(1-\delta(p))=\delta(p)(1-p^i) \in
A^\times$ si $i\neq 0$. \end{pf}

De la discussion pr\'ec\'edent la proposition~\ref{critlibre} on d\'eduit le :

\begin{cor} Soit $\delta \in \cT(A)$ tel que $1-\delta(p)p^i$ est non diviseur
de $0$ dans $A$ pour tout $i \in \Z$. Alors le $\Ro_A^+(\Gamma)$-module
$C(\Ro_A(\delta))$
est projectif si et seulement si $1-\delta(p)p^i \in A^\times$ pour tout
$i\geq 0$, auquel cas il est libre de rang $1$.
\end{cor}

\subsection{Cohomologie de $\Ro_A(\delta)$ partie III : calcul du $H^1$}
Le calcul de $H^1(\Ro_A(\delta))$ est maintenant une formalit\'e. On
rappelle que $x \in \cT(\Q_p)$ d\'esigne le caract\`ere tautologique
identit\'e et que $\chi \in \cT(\Q_p)$ est le caract\`ere tel que
$\chi(p)=1$ et $\chi_{|\Z_p^\times}=x_{|\Z_p^\times}$ ("caract\`ere
cyclotomique").

\begin{definition} On d\'esigne par $\cT^{\rm reg}
\subset \cT$ l'ouvert compl\'ementaire de l'ensemble discret\footnote{Ils sont discrets au sens que tout ouvert affino\"ide de $\cT$ ne rencontre qu'un nombre fini de tels points.} des
points $\Q_p$-rationnels de la forme $x^{-i}$ ou $\chi x^i$ pour $i\geq
0$.\par
Un caract\`ere $\delta \in \cT(A)$ est dit r\'egulier si il est dans
$\cT^{\rm reg}(A)$, ce qui revient \`a dire que pour tout $z \in {\rm Sp}(A)$ le caract\`ere $\delta_z : \Qp^\times
\rightarrow k(z)^\times$ obtenu par \'evaluation en $z$ n'est pas de la forme $x^i$
ou $\chi x^{-i}$ pour $i\geq 0$ entier. 
\end{definition}

On note $K(A)$ le groupe de Grothendieck des $A$-modules de type fini, 
et si $M$ est un tel $A$-module on note $[M]$ sa classe dans $K(A)$. On
rappelle qu'une alg\`ebre affino\"ide est noeth\'erienne.

\begin{thm}\label{thmH1} Soit $D$ un $\fg$-module triangulin de rang $d$ sur $\Ro_A$.
Alors $H^i(D)$ est de type fini sur $A$ pour tout entier $i$ et on a la
relation dans $K(A)$ $$[H^0(D)]-[H^1(D)]+[H^2(D)]=-[A^d].$$\ps 
Si de plus le param\`etre $(\delta_i)$ de $D$ est dans $\cT^{\rm
reg}(A)^d$, alors $H^0(D)=H^2(D)=0$ et $H^1(D)$ est libre de rang $d$ sur $A$. Enfin, si $D=\Ro_A(\delta)$ avec $\delta \in
\cT^{\rm reg}(A)$, les quatre
morphismes naturels $$ (D^+)^{\psi=0}/(\gamma-1)
\leftarrow C(D^+)/(\gamma-1) \rightarrow C(D)/(\gamma-1) \leftarrow D^{\psi=1}/(\gamma-1)
\rightarrow H^1(D)$$
sont des isomorphismes, un g\'en\'erateur de $(D^+)^{\psi=0}/(\gamma-1)$
\'etant donn\'e par la classe de $1+T$.
\end{thm}

\begin{pf} Par r\'ecurrence sur $d$ et en utilisant la suite longue de cohomologie,
on peut supposer que $d=1$, i.e. $D=\Ro_A(\delta)$. On a d\'ej\`a vu que
$D^{\varphi=1}$ et $D/(\psi-1)$ sont de type fini sur $A$, et donc \`a plus
forte raison que $H^0(D)$ et $H^2(D)$ le sont aussi. Pour v\'erifier le
premier point du th\'eor\`eme, et compte tenu du
d\'evissage donn\'e par le th\'eor\`eme~\ref{devissagecoho}, il suffit donc de d\'emontrer que
$C(D)/(\gamma-1)$ est de type fini sur $A$ et de classe $[A^d]$ dans $K(A)$.
(On rappelle que si $X$ est un $A$-module de type fini et $u \in \End_A(X)$
alors la multiplication par $u$ sur $X$ implique l'identit\'e 
$[X^{u=0}]=[X/u(X)]$ dans $K(A)$.) Les deux suites exactes donn\'ees par le lemme~\ref{devissedplus} (ii) et la
proposition~\ref{struccplus} induisent des suites exactes 
\begin{equation}\label{sexacte1} \begin{split} 0 \longrightarrow {\rm
Pol}(\Z_p,A)^{\delta(p)^{-1}\psi=1,\gamma=1} \longrightarrow
C(D^+)/(\gamma-1) \longrightarrow \\
C(D)/(\gamma-1) \overset{(1-\varphi)^{-1} C}{\longrightarrow} 
{\rm Pol}(\Z_p,A)^{\delta(p)^{-1}\psi=1}/(\gamma-1) \longrightarrow 0.
\end{split}
\end{equation}
\begin{equation}\label{sexacte2}\begin{split} 0 \longrightarrow (\oplus_{i=0}^{k-1}
A/(1-\delta(p)p^i)
\,\,\cdot \widetilde{\delta x^i})^{\gamma=1} \longrightarrow
C(D)^+/(\gamma-1)
\longrightarrow \\ \Ro_A^+(\Gamma)\cdot (1+T)/(\gamma-1) \overset{J_k}{\longrightarrow}
(\oplus_{i=0}^{k-1}
A/(1-\delta(p)p^i)
\,\,\cdot \widetilde{\delta x^i})/(\gamma-1) \longrightarrow 0
\end{split}
\end{equation}
car $\gamma-1$ est injectif sur $C(D)$ et $(\Ro_A^+)^{\psi=0}$
(th\'eor\`eme~\ref{invgamma}). On conclut car $\Ro_A^+(\Gamma)/([\gamma]-1)=A$
(lemme~\ref{lemmeker} pour $\delta=1$). \par
Supposons maintenant $\delta$ r\'egulier. Dans ce cas, les noyaux et
conoyaux de $\gamma-1$ sur $D^{\varphi=1}$, $D/(\psi-1)$, ${\rm
Pol}(\Z_p,A)^{\delta(p)^{-1}\psi=1}$ et $\oplus_{i\geq 0} A/(1-\delta(p)p^i)
\widetilde{\delta x^i}$ sont tous nuls par hypoth\`ese. En effet, si $a, b
\in A$ sont tels que $(a,b)=A$, alors la multiplication par $a$ est
bijective sur $A/bA$ et $A[b]$ (la $b$-torsion dans $A$): si $au+bv=1$ la multiplication par $u$ en est un
inverse. Cela montre que $H^0(D)=H^2(D)=0$ puis que les quatre fl\`eches de
l'\'enonc\'e sont des isomorphismes. En particulier, $H^1(D)$ est isomorphe
\`a $\Ro_A^+(\Gamma)/(\gamma-1)=A$. 
\end{pf}

\begin{thm}\label{changementdebase} Soient $A \rightarrow B$ un morphisme de $\Q_p$-alg\`ebres
affino\"ides et $D$ un $\fg$-module triangulin sur $\Ro_A$. Supposons que 
$A \rightarrow B$ est plat ou que le param\`etre de $D$ est r\'egulier.
Alors pour tout entier $i$ l'application naturelle $$H^i(D) \otimes_A
B \longrightarrow H^i(D \widehat{\otimes}_A B)$$ est un isomorphisme.
\end{thm}

\begin{pf} L'application de l'\'enonc\'e est celle d\'eduite du $A$-morphisme
de complexes $C_{\varphi,\gamma}(D)^\bullet \rightarrow C_{\varphi,\gamma}(D
\widehat{\otimes}_A B)^\bullet$. En particulier, si on a une suite exacte de
$\fg$-modules sur $\Ro_A$ disons $0
\rightarrow D_1 \rightarrow D_2 \rightarrow D_3 \rightarrow 0$, et donc une suite exacte
$0
\rightarrow D_1 \widehat{\otimes}_A B \rightarrow D_2 \widehat{\otimes}_A B \rightarrow D_3
\widehat{\otimes}_A B \rightarrow 0$, on dispose d'une morphisme
$A$-lin\'eaire naturel entre les suites exactes longues de cohomologie.
Il vient que par d\'evissage et par le lemme des $5$, on peut supposer
$D=\Ro_A(\delta)$ est de rang $1$. 

Soit $\delta \in \cT(A)$, on note $\delta_B$ l'image de $\delta$ dans
$\cT(B)$. Remarquons que si $F(\delta)$ d\'esigne
le $A[\varphi,\Gamma]$-module
$A[t](\delta)$ et
${\rm Pol}(\Z_p,A)(\delta)$, alors l'application naturelle $F(\delta) \otimes_A B
\rightarrow F(\delta_B)$ est un isomorphisme. En particulier, pour
$\ast \in \{\varphi,\psi\}$ et si
$C_{\ast,\gamma}(F(\delta))$ d\'esigne le complexe \`a trois termes \'evident,
dont on d\'esignera par $H^i_\ast(F(\delta))$ le $A$-module de cohomologie, 
et si $A \rightarrow B$ est plat, alors 
$$H^i_\ast(F(\delta)) \otimes_A B \isomo H^i_\ast(F(\delta)\otimes_A B) \isomo
H^i_\ast(F(\delta_B)).$$
Mais on a d\'efini (prop.~\ref{corH0H2}) des isomorphismes naturels
$H^0_\varphi(F(\delta)) \isomo
H^0(\Ro_A(\delta))$ et $H^2(\Ro_A(\delta)) \isomo H^2_\psi(F(\delta))$ ($F$ valant
respectivement $A[t]$ dans le premier cas et ${\rm Pol}$ dans le second), le
th\'eor\`eme en d\'ecoule pour $i=0,2$. 

Supposons toujours $D=\Ro_A(\delta)$. Par un argument similaire \`a celui ci-dessus utilisant les suites
exactes~\eqref{sexacte1} et~\eqref{sexacte2}, et le fait que la formation des modules
$(D^+)^{\psi=0}$, $D^{\psi=1}$, $C(D)$ et $C(D^+)$ est fonctorielle en $A$,
ainsi donc que leurs quotients par $\gamma-1$, le th\'eor\`eme suit dans le
cas $i=1$ si l'on montre que pour tout morphisme $A \rightarrow B$ l'application naturelle
$$f : (\Ro_A^+(\delta))^{\psi=0}/(\gamma-1) \otimes_A B \rightarrow
(\Ro_B^+(\delta_B))^{\psi=0}/(\gamma-1)$$ 
est un isomorphisme. Mais on a d\'ej\`a vu que $\Ro_A^+(\delta)^{\psi=0}$ est libre sur
$\Ro_A^+(\Gamma)$ engendr\'e par $1+T$, de sorte que son quotient par
$\gamma-1$ est libre sur $A$ engendr\'e par la classe $\overline{1+T}$ de
$1+T$. Mais par construction $f$ provient par quotient du morphisme \'evident
$\Ro_A^+(\delta)^{\psi=0} \rightarrow \Ro_B^+(\delta_B)^{\psi=0}$ qui envoie $1+T$
sur $1+T$.\end{pf}

Remarquons que jusqu'ici nous n'avons pas donn\'e d'\'enonc\'e exact sur la
structure de $H^1(\Ro_A(\delta))$ quand $\delta$ n'est pas r\'egulier. Si
l'on concat\`ene les suites exactes donn\'ees par la
proposition~\ref{devissagecoho} et les formules~\eqref{sexacte1}
et~\eqref{sexacte2} nous en obtenons un d\'evissage explicite, bien 
que peu rago\^utant en
g\'en\'eral. Ce d\'evissage se simplifie dans le cas utile suivant.

\begin{definition} On dit que $\delta \in \cT(A)$ est bien pla\c{c}\'e si :
\begin{itemize}\item[(i)] pour tout
$i \in \Z$, alors $1-\delta(p)p^i$ est non diviseur de $0$ dans $A$,\ps
\item[(ii)] pour tout $i \geq 0$, l'image de $1-\delta(\gamma)\gamma^{1-i}$
dans\footnote{Si $p=2$, il faut remplacer $A/(1-\delta(p)p^{-i})$ par
$A/(\delta(-1)+1,1-\delta(p)p^{-i})$.} $A/(1-\delta(p)p^{-i})$ est non
diviseur de $0$. 
\end{itemize}
\end{definition}

\begin{prop}\label{h1bienplace} Si $\delta \in \cT(A)$ est bien pla\c{c}\'e et
$D=\Ro_A(\delta)$, alors $H^0(D)=0$ et les
morphismes naturels $$ C(D^+)/(\gamma-1) \rightarrow C(D)/(\gamma-1) \leftarrow
D^{\psi=1}/(\gamma-1)
\rightarrow H^1(D)$$
sont des isomorphismes. On a de plus une suite exacte naturelle de $A$-modules $$ 0
\longrightarrow \prod_{i\geq 0}
(A/(1-\delta(p)p^i))[1-\delta(\gamma)\gamma^i] \longrightarrow H^1(D)
\longrightarrow \bigcap_{i\geq 0} (1-\delta(p)p^i,1-\delta(\gamma)\gamma^i)
\longrightarrow 0.$$

\end{prop}
L'intersection dans le terme de droite de cette derni\`ere suite est sous-entendue \`a l'int\'erieur
de $A$ (c'est donc un id\'eal de $A$). 

\begin{pf} En effet, l'annulation de $D^{\varphi=1}$,
$D/(\psi-1))^{\gamma=1}$ et ${\rm Pol}(\Z_p,A)^{\delta(p)^{-1}\psi=1}$
\'equivaut au caract\`ere bien plac\'e de $\delta$, et on conclut la
premi\`ere assertion par la
proposition~\ref{devissagecoho} et la suite exacte~\eqref{sexacte1}.
D'apr\`es la proposition~\ref{struccplus}, on a une suite exacte de
$\Ro_A^+(\Gamma)$-modules $$0 \longrightarrow C(D^+) \longrightarrow
\Ro_A^+(\Gamma) \longrightarrow \prod_{i\geq 0} A/(1-\delta(p)p^i)(\delta)
\longrightarrow 0$$
dont la derni\`ere assertion se d\'eduit en appliquant $\gamma=1$, en
utilisant que $(\gamma-1)$ est injectif sur $\Ro_A^+(\Gamma)$
(lemme~\ref{lemmeker} pour $\delta=1$). (Quand $p=2$ l'argument ci-dessus
et l'\'enonc\'e ne sont \'evidemment pas tout \`a fait corrects)
\end{pf}

Terminons par un cas particulier important concernant le caract\`ere universel. Si $U \subset
\cT$ est un ouvert affino\"ide, on d\'esigne par $\delta_U \in \cT(U)$
le caract\`ere tautologique.

\begin{thm}\label{h1universel} Pour tout ouvert affino\"ide $U \subset \cT$, 
$H^0(\Ro_U(\delta_U))=0$ et le $\OO(U)$-module $H^1(\Ro_U(\delta_U))$
s'identifie naturellement \`a l'id\'eal de $\OO(U)$ constitu\'e des fonctions qui
s'annulent en tous les points de $U$ param\'etrant les caract\`eres de la forme
$x^{-i}$
pour $i\geq 0$. 

De plus, pour tout $z \in U$, de caract\`ere associ\'e
$\delta_z$, l'application naturelle
	$$H^1(\Ro_U(\delta_U))\otimes_{\OO(U)}k(z) \longrightarrow H^1(\Ro_{k(z)}(\delta_z))$$
est un isomorphisme, \`a moins que $\delta_z$ ne soit de la forme $\chi
x^i$ avec $i\geq 0$, auquel cas cette application est nulle.

Plus pr\'ecis\'ement, supposons que le seul point non r\'egulier de $U$,
disons $u \in U$, param\`etre un caract\`ere de la forme $\chi x^i$ avec $i\geq 0$ et soit 
$m=m_u \subset \OO(U)$ l'id\'eal maximal des fonctions s'annulant en ce
point. On a $k(u)=\Q_p$ et on consid\`ere l'espace tangent
$T_u=\Hom_{\Qp}(m/m^2,\Q_p)$. Alors : 
\begin{itemize}
\item[(i)] Le morphisme naturel $H^1(m \Ro_U(\delta)) \rightarrow
H^1(\Ro_U(\delta_U))$ est un isomorphisme entre $\OO(U)$-module libres de
rang $1$, \par
\item[(ii)] Le $\OO(U)$-morphisme canonique $m \rightarrow m/m^2$ induit une injection
$$H^1(m \Ro_U(\delta_U))
\otimes_{\OO(U)} k(u) \longrightarrow \Hom_{\Q_p}(T_u,H^1(\Ro(\chi x^i)))$$
dont l'image est une droite constitu\'ee d'isomorphismes. En particulier,
cette droite induit une isomorphisme canonique $ \mathbb{P}(T_u) \isomo
\mathbb{P}(H^1(\Ro(\chi x^i))$ entre espaces projectifs sur $\Q_p$ de dimension $1$. 
\end{itemize}
\end{thm}

On notera le r\^ole non sym\'etrique jou\'e ici par les points singuliers
de la forme $x^{-i}$ et ceux de la forme $\chi x^i$. Nous pr\'eciserons un peu plus loin ce th\'eor\`eme en introduisant
l'\'eclat\'e de $U$ aux points non r\'eguliers, ce qui nous permettra
notamment de comprendre compl\`etement la structure analytique de l'espace
des $\fg$-modules triangulins de rang $2$, y compris au voisinage des points
de param\`etre singulier  : c'est exactement la structure sugg\'er\'ee par
Colmez dans sa d\'efinition de l'espace des triangulines en rang $2$.

\begin{pf} Le caract\`ere $\delta_U$ est \'evidemment bien plac\'e. De
m\^eme, pour tout $i \geq 0$ l'\'el\'ement $1-\delta(\gamma)\gamma^i$ n'est
pas diviseur de z\'ero sur l'anneau (localement int\`egre)
$\OO(U)/(1-\delta(p)p^i)$. La premi\`ere assertion d\'ecoule donc de la
proposition~\ref{h1bienplace}. Pour v\'erifier la seconde assertion, il
r\'esulte de la commutativit\'e de $H^1$ au changement de base plat (ici une
immersion ouverte) que l'on peut supposer que $U$ ne contient qu'un seul 
point singulier, disons $z$. 
Si $m \subset \OO(U)$
d\'esigne l'id\'eal maximal des fonctions qui s'annulent en $z$, alors il
est classique que si $a=\delta_U(p)-\delta_z(p)$ et
$b=\delta_U(\gamma)-\delta_z(\gamma)$ alors $m=(a,b)$ et on a 
des suites exactes $$0 \longrightarrow \OO(U) \overset{f\mapsto
(af,bf)}{\longrightarrow}
\OO(U)^2 \overset{(f,g) \mapsto af-bg}{\longrightarrow} 
m \longrightarrow 0,$$
et $$ 0 \longrightarrow m \longrightarrow \OO(U) \longrightarrow k(z)
\longrightarrow 0.$$
Comme $\Ro_U$ est plat sur $\OO(U)$ (lemme~\ref{topo} (vi)), ces
suites restent exactes apr\`es $- \otimes_{\OO(U)} D$, $D=\Ro_U(\delta_U)$, de sorte
qu'en prenant la suite longue de cohomologie on obtienne des suites exactes de
$\OO(U)$-modules 
$$ 0 \longrightarrow H^0(mD) \longrightarrow H^1(D) \overset{(a,b)}{\longrightarrow} H^1(D)^2
\longrightarrow H^1(mD) \longrightarrow H^2(D)
\overset{(a,b)}{\longrightarrow} H^2(D)^2 \longrightarrow H^2(mD) \longrightarrow 0,$$
$$ 0 \longrightarrow H^0(D_z) \longrightarrow H^1(mD) \longrightarrow H^1(D)
\longrightarrow H^1(D_z) \longrightarrow H^2(mD) \longrightarrow H^2(D)
\longrightarrow H^2(D_z) \longrightarrow 0.$$
On a utilis\'e que $D \otimes_{\OO(U)} m = mD$ (platitude de $\Ro_U$ sur
$\OO(U)$).

Si $\delta_z=x^{-i}$, alors la
proposition~\ref{corH0H2} assure que $H^2(D)=H^2(D_z)=0$. 
On d\'eduit des suites ci-dessus qu'alors $H^2(mD)=0$ puis que $H^1(D)
\rightarrow H^1(D_z)$ est surjectif. Comme $H^1(D) \simeq m$ et que $m/m^2$
et $H^1(D_z)$ sont de dimension $2$ sur $k(z)$, l'application $H^1(D)
\otimes_{\OO(U)}k(z) \rightarrow H^1(D_z)$ est un isomorphisme.

Si $\delta_z=\chi x^i$, la proposition~\ref{corH0H2} assure que $H^2(D)
\simeq k(z)$ est tu\'e par $m$, de sorte que la fl\`eche $H^2(D) \rightarrow
H^2(D)^2$ ci-dessus est nulle.  Ainsi, $H^2(D)^2 \isomo H^2(mD)$ et ce
dernier est isomorphe \`a $k(z)^2$. Comme d'autre part $H^2(D) \rightarrow
H^2(D_z)$ est un isomorphisme (car surjectif), on en d\'eduit que $H^1(D_z) \longrightarrow
H^2(mD)$ est surjective : c'est donc un isomorphisme pour des raisons de
dimension. Ainsi, $H^1(D) \rightarrow H^1(D_z)$ est nul. Comme $H^0(D_z)=0$
on en d\'eduit enfin $H^1(mD)=H^1(D)$ : il ne reste qu'\`a prouver le (ii)
du th\'eor\`eme.

Comme $m/m^2$ est annul\'e par $m$, il vient que 
$D \otimes_{\OO(U)}(m/m^2)= D_z \otimes_{\Q_p} (m/m^2)$ (avec action de
$\varphi$ et $\Gamma$). On en d\'eduit pour tout $i$ un isomorphisme
canonique de $\OO(U)$-modules $$H^i(D
\otimes_{\OO(U)}(m/m^2)) \isomo H^i(D_z) \otimes_{\Q_p} (m/m^2) \isomo \Hom_{\Q_p}(T_z,H^i(D_z)).$$
Remarquons que $T_z$ et $H^1(D_z)$ sont des $\Q_p$-espace vectoriel de dimension $2$.
Consid\'erons le morphisme $$\mu : H^1(mD) \rightarrow
H^1(D\otimes_{\OO(U)}m/m^2) = \Hom_{\Q_p}(T_z,H^1(D_z))$$
d\'eduit du morphisme naturel $mD \rightarrow
mD/m^2D=(m/m^2)\otimes_{\OO(U)}D$ (on rappelle que $\Ro_U$ est plat sur
$\OO(U)$). Comme $H^1(mD)=H^1(D)$ est libre de rang $1$, son image dans
$\Hom_{\Q_p}(T_z,H^1(D_z))$ est soit nulle soit une $\Q_p$-droite. 
Si $L \neq 0 \subset T_z$, $\mu(L)$ est par d\'efinition
l'\'el\'ement de $H^1(D_z)$ obtenu comme compos\'e de $\mu$ et de
l'applicatio naturelle $(L \otimes 1) : (m/m^2) \otimes_{\Q_p} D_z
\rightarrow D_z$, ou ce qui revient au m\^eme comme image de l'application
$\OO(U)$-lin\'eaire naturelle $mD \rightarrow D_z$ apr\`es passage au $H^1$.
Soit $m^2 \subset J \subset m$ le noyau de $L$. La suite exacte de
$\OO(U)$-modules $0 \rightarrow J \rightarrow m \rightarrow m/J=k(z) \rightarrow
0$ reste exacte apr\`es extension des scalaires \`a $\Ro_U$, de sorte que
l'on dispose d'une suite longue de $\OO(U)$-modules 
$$ H^1(mD) \rightarrow H^1(D_z) \rightarrow H^2(JD) \rightarrow H^2(mD)
\rightarrow H^2(D_z) \longrightarrow 0.$$
L'image de la premi\`ere fl\`eche est le $\Q_p$-module engendr\'e par $\mu(L)$. Pour
conclure il suffit donc de voir que le noyau de la fl\`eche $H^1(D_z) \rightarrow H^2(JD)$ est
une doite. On a d\'ej\`a vu que $H^1(D_z) \simeq H^2(mD)$ et $H^2(D_z)$ sont de
$\Q_p$-dimensions
respectives $2$ et $1$, il suffit donc de voir que $H^2(JD)$ est de
dimension $2$ sur $\Q_p$. Mais tout comme $m$, l'id\'eal $J$ a deux
g\'en\'erateurs et admet une
pr\'esentation de la forme $$ 0 \longrightarrow \OO(U)
\overset{g}{\longrightarrow} \OO(U)^2 \longrightarrow J \longrightarrow 0$$
avec $g \otimes_{\OO(U)} k(u) = 0$, de sorte qu'un argument d\'ej\`a donn\'e
plus haut montre que $H^2(D)^2 \simeq H^2(JD)$, ce qui conclut la preuve du
th\'eor\`eme.\end{pf}

\section{L'espace des $\fg$-modules triangulins sur $\Ro_A$}\label{San}

%
%
%

\subsection{$\fg$-modules triangulins r\'eguliers rigidifi\'es.}

Fixons $d\geq 1$ un entier. L'espace  $\cT^d$ est muni d'une famille universelle $(\widetilde{\delta}_i)$ de 
caract\`eres $\deltat_i : \Q_p^\ast \rightarrow \OO(\cT^d)^\ast$. On notera $$\cT^{\rm reg}_d \subset \cT^d$$
l'ouvert admissible (en fait, de Zariski) d\'efini par les relations
$\widetilde{\delta}_i/\widetilde{\delta}_j \in \cT^{\rm reg}$ pour tout
$1 \leq i < j \leq d$. Par d\'efinition, si $A$ est une alg\`ebre affino\"ide 
alors $\cT_d^{\rm reg}(A)$ est donc le sous-ensemble des $(\delta_i) \in \cT(A)^d$ tels que 
$\delta_i/\delta_j \in \cT^{\rm reg}(A)$ pour tout $1\leq i < j\leq d$. \ps

\begin{definition} Un $\fg$-module triangulin {\it r\'egulier rigidifi\'e} est un triplet $(D,\Fil_\bullet(D),\nu)$ sur $\Ro_A$ 
o\`u : \begin{itemize}\ps
\item[-] $(D,\Fil_\bullet(D))$ est un $\fg$-module triangulin sur $\Ro_A$ dont le param\`etre $(\delta_i)$ est dans $\cT_d^{\rm reg}(A)$ (condition de r\'egularit\'e), \ps
\item[-] $\nu=(\nu_i)$ est une famille d'isomorphismes $\nu_i :
\Fil_{i+1}(D)/\Fil_i (D) \isomo \Ro_A(\delta_i)$ dans $\fg/A$ pour $i=0,\dots,{\rm rang}_{\Ro_A}(D)-1$ (rigidification). \end{itemize}
Deux tels
triplets $(D,\Fil_\bullet(D),\nu)$ et $(D',\Fil_\bullet(D'),\nu')$
seront dit {\it \'equivalents} si il existe un isomorphisme $f: D \rightarrow D'$ dans $\fg/A$ envoyant 
$\Fil_i(D)$ sur $\Fil_i(D')$ pour tout $i$ (auquel cas $D$ et $D'$ ont m\`eme param\`etre) et tel que $\nu'_i \,\,{\rm o}
\,\, f = \nu_i$ pour tout $i=0,\dots,{\rm rang}_{\Ro_A}(D)-1$. 
\end{definition}

Consid\'erons le foncteur $$F_d^\square : \Aff \longrightarrow {\rm Ens}$$
de la cat\'egorie $\Aff$ des $\Q_p$-alg\`ebres affino\"ides vers celle ${\rm Ens}$ des ensembles associant \`a chaque objet $A$ l'ensemble\footnote{Pour les raisons usuelles il ne s'agit pas vraiment d'un ensemble. Pour contourner ce probl\`eme il suffit de rajouter une fois pour toutes dans la d\'efinition d'un $\fg$-module triangulin r\'egulier rigidifi\'e de rang $d$ sur $\Ro_A$ que le $\Ro_A$-module sous-jacent est $\Ro_A^d$ (plut\^ot que simplement, isomorphe \`a $\Ro_A^d$). }  des classes
d'\'equivalence de $\fg$-modules triangulins $p$-r\'eguliers rigidifi\'es sur $\Ro_A$.
Si $u=(D,\Fil_\bullet(D),\nu) \in F_d^\square(A)$ et si $f: A \rightarrow B$ est un morphisme dans $\Aff$, on pose bien 
entendu $$F_d^\square(f)(u)=(D\widehat{\otimes}_A B,(\Fil_i(D) \widehat{\otimes}_A B), (\nu_i \otimes_{\Ro_A} \Ro_B))  
\in F_d^\square(B),$$ ce qui fait bien de $F_d^\square$ un foncteur covariant. Le param\`etre fournit un morphisme de 
foncteurs $\delta : F_d^\square \rightarrow \cT_d^{\rm reg}$. 
\newcommand{\cS}{\mathcal{S}}
\begin{thm}\label{repfdsquare} Le foncteur $F_d^\square$ est repr\'esentable par un espace analytique $p$-adique 
$\cS_d^\square$. Le morphisme $\delta: \cS_d^\square \rightarrow \cT_d^{\rm reg}$ est lisse de dimension 
relative $\frac{d(d-1)}{2}$. L'espace $\cS_d^\square$ est irr\'eductible,
r\'egulier, et \'equidimensionnel de dimension $\frac{d(d+3)}{2}$.
\end{thm}

Nous aurons besoin du lemme suivant.

\begin{lemma} (Rigidit\'e) Soient $(D,\Fil_\bullet(D),\nu)$ et $(D',\Fil_\bullet(D'),\nu')$ des $\fg$-modules 
triangulins r\'eguliers rigidifi\'es sur $\Ro_A$. S'ils sont \'equivalents, alors il existe une unique \'equivalence 
entre eux. 
\end{lemma}

\begin{pf} On peut supposer $(D,\Fil_\bullet(D),\nu)=(D',\Fil_\bullet(D'),\nu')$ et il s'agit de voir que ce triplet 
a pour seule auto-\'equivalence l'identit\'e. Si $f : D \rightarrow D$ en est une, alors par hypoth\`ese $f(\Fil_i(D))=\Fil_i(D)$ et $f$ induit l'identit\'e sur chaque $\Fil_{i+1}(D)/\Fil_i(D)$. Ainsi, $u:=f-{\rm id} \in \End_{\fg/A}(D)$ a la propri\'et\'e que $u(\Fil_{i+1}(D)) \subset \Fil_i(D)$ pour tout $i<{\rm rang}_{\Ro_A}(D)$. Pour voir que $u=0$ il suffit donc de voir que $$\Hom_{\fg/A}(\Ro_A(\delta_j),\Ro_A(\delta_i))=0$$
d\`es que $j>i$, soit encore que $H^0(\Ro_A(\delta_i\delta_j^{-1})=0$ sous cette hypoth\`ese. Mais ceci vient de ce que $\delta_i\delta_j^{-1}$ est r\'egulier et du th\'eor\`eme~\ref{thmH1}.
\end{pf}

%

D\'emontrons maintenant le th\'eor\`eme. Quand $d=1$, $F_d^\square=\cT$ et le r\'esultat est \'evident. 
Pour $d\geq 2$ on proc\`ede par r\'ecurrence sur $d$. On dispose d'un morphisme de foncteurs \'evident $F_d^\square \rightarrow F_{d-1}^\square \times F_1^\square  = \cS_{d-1}^\square \times \cT$, 
associant \`a la classe de $(D,\Fil_\bullet(D),\nu) \in F_d^\square(A)$ la paire form\'ee de la classe de 
$(\Fil_{d-1}(D),(\Fil_i(D))_{i\leq d-1},(\nu_i)_{1\leq i \leq d-1})$ et de $\delta_d$. Il se factorise par 
l'ouvert Zariski $$U_d \subset \cS_{d-1}^\square \times \cT$$ 
qui est l'image inverse de l'ouvert $\cT_d^{\rm reg} \subset \cT_{d-1}^{\rm reg} \times \cT$ par le morphisme 
param\`etre $\pi_d: \cS_{d-1}^\square \times \cT \rightarrow \cT_{d-1}^{\rm reg} \times \cT$. Nous allons d\'emontrer que 
le morphisme ci-dessus $F_d^\square\rightarrow U_{d}$ est relativement repr\'esentable par un fibr\'e vectoriel de rang
$d-1$, ce qui provera le th\'eor\`eme. Le fibr\'e en question sera trivial
au dessus de tout ouvert affino\"ide de $U_d$. \ps 
\newcommand{\cM}{{\mathcal M}}

Nous aurons besoin d'un sorite pr\'eliminaire. Soient $A$ une $\Q_p$-alg\`ebre
affino\"ide et $u=(c_A,\delta_d) \in U_d(A)$ et $x=(D_A,\Fil_\bullet,\nu)$ un repr\'esentant de la classe
$c_A$. Consid\'erons le
$A$-module $M(u)=H^1(D_A(\delta_d^{-1})$. Si $y=(D'_A,\Fil_\bullet,\nu')$ est \'equivalent \`a $x$, l'unique
\'equivalence $y \rightarrow x$ identifie
donc canoniquement $M(y)$ et $M(x)$. Ainsi, il y a un sens \`a d\'efinir le
$A$-module $M(u)$ associ\'e \`a
un \'el\'ement $u \in U_d(A)$, comme \'etant par exemple la limite
inductive des $M(x)$ pour $x$
parcourant l'ensemble\footnote{Voir la note pr\'ec\'edente!} des
repr\'esentants de $c_A$ : le choix d'un repr\'esentant $x$ de $c_A$ fournit alors un isomorphisme
canonique $M(c_A) \isomo M(x)$. On v\'erifie de suite que si $A \rightarrow B$ est un
morphisme entre alg\`ebres affino\"ides, alors pour tout
$u_A \in U_d(A)$, d'image $u_B \in U_d(B)$, on dispose d'un
morphisme canonique
$$M(u_A) \otimes_A B \longrightarrow M(u_B)$$
d\'efini de mani\`ere \'evidente sur les repr\'esentants. Les
th\'eor\`emes~\ref{changementdebase}
et~\ref{thmH1} assurent que c'est un isomorphisme entre modules libres de
rang $d-1$. En particulier $\Omega \mapsto M(\Omega)$, pour $\Omega
\subset U_d$ ouvert affino\"ide, d\'efinit un
faisceau coh\'erent sur $U_d$ tel qu'en fait $M(\Omega)$ est libre de rang
$d-1$ sur $\OO(\Omega)$ pour tout $\Omega$. On note encore $\cM$ ce faisceau
coh\'erent sur $U_d$. D'apr\`es ce que nous venons de voir, si $u \in
U_d(A)$, alors on a une identification canonique
\begin{equation}\label{idcan} u^*(\cM)(A)=M(u).\end{equation}\par
Consid\'erons alors
$$\eta : \cS_d^\square : = {\rm Spec}_{U_d}^{\rm an}({\rm
Symm}\,\,\cM^\vee) \rightarrow U_d$$ le spec relatif analytique de la
$\OO_{U_d}$-alg\`ebre quasi-coh\'erente ${\rm
Symm}\,\,\cM^\vee$ (voir \cite[\S 2.2]{conradample}) : c'est le fibr\'e vectoriel sur $U_d$ associ\'e \`a $\cM^\vee$. 
Par la propri\'et\'e universelle de cette construction, pour
toute alg\`ebre affino\"ide $A$ et tout $u \in U_d(A)$, disons 
$u=([(D_A,\Fil_\bullet,\nu)],\delta_d)$, on a des identifications canoniques 
$$\{ v \in \cS_d^\square(A),
\eta(v)=u\}=\Hom_A(u^*(\cM^\vee)(A),A)=u^*(\cM)(A)=H^1(u) \isomo
 {\rm Ext}(\Ro_A(\delta_d),D_A).$$ 
Cela d\'efinit un morphisme de foncteurs $\cS_d^\square \rightarrow F_d^\square$
au dessus de $U_d$ : on associe \`a une paire form\'ee d'un \'el\'ement 
$([(D_A,\Fil_\bullet,\nu)],\delta_d) \in U_d(A)$ et d'une classe $E \in {\rm
Ext}(\Ro_A(\delta_d),D_A)$, que l'on voit comme la donn\'ee d'une suite exacte $$
0 \longrightarrow D_A \overset{\iota}{\longrightarrow}  E \overset{\pi}{\longrightarrow}
\Ro_A(\delta) \longrightarrow 0,$$
la classe du $\fg$-module triangulin $E$ avec $\Fil_i(E):=\iota(\Fil_i(D_A))$
et $\nu_i {\rm \cdot} \iota^{-1}=:\nu_i$ pour
$i\leq d-1$, et ${\rm Fil}_{d}(E)=E$ et $\nu_{d}=\pi$. Il est imm\'ediat que
$\cS_d^\square \rightarrow F_d^\square$ est un isomorphisme de foncteurs. $\square$

Notons $\mathbb{D}^r$ la boule unit\'e affino\"ide ferm\'ee de rayon $1$ sur $\Q_p$,
d'alg\`ebre $\Q_p\langle t_1,\dots,t_r\rangle$.

\begin{cor}\label{corfamcoleman} Si $x \in S_d^\square$, il existe un voisinage ouvert
affino\"ide $U$ de $x$ dans $S_d^\square$, un voisinage ouvert
affino\"ide $\Omega$ de $\delta(x)$ dans $\cT_d^{\rm reg}$, et un 
isomorphisme $$\iota :
U \isomo \Omega \times {\mathbb{D}}^{\frac{d(d-1)}{2}}$$
tels que ${\rm pr}_2 \cdot \iota = \delta$.
\end{cor}

\begin{pf}Par d\'efinition, si $E$ est un fibr\'e vectoriel sur un espace
rigide $Y$ alors pour
tout $y \in Y$ on peut trouver un voisinage ouvert $W$ de $y$ dans $Y$ tel
que 
$W \times_Y E \isomo W \times \AAA^m$ comme fibr\'e. Le corollaire d\'ecoule
alors de la construction
inductive de $\cS_d^\square$ \'etablie ci-dessus. 
\end{pf}

\subsection{$\fg$-modules triangulins r\'eguliers  non
rigidifi\'es} \label{diverstrucs} Bien que ce ne soit pas n\'ecessaire pour le th\'eor\`eme
principal de cet article, il est naturel de consid\'erer le foncteur $$F_d :
\Aff \rightarrow {\rm Ens}$$ o\`u $F_d(A)$ est l'ensemble des classes
d'\'equivalence de $\fg$-modules triangulins $(D,\Fil_\bullet(D))$ sur
$\Ro_A$ dont le param\`etre est dans $\cT_d^{\rm reg}(A)$ (sans
rigidification).  La notion d'\'equivalence utilis\'ee ici est celle dans
$\fg/A$ avec pr\'eservation de la filtration.

Pour \'eliminer les
auto-\'equivalences des objets param\'etr\'es par $F_d$ il est n\'ecessaire
de se restreindre \`a un sous-foncteur ad\'equat. Si $L$ est une extension finie de $\Q_p$ et 
si $(D,\Fil_\bullet(D))$ est un $\fg$-module triangulin sur $\Ro_L$, on dira que $D$ est {\it non scind\'e} 
si pour tout $0 \leq i < {\rm rg}_{\Ro_L}(D)$, l'extension $$0 \longrightarrow \Fil_i(D) \longrightarrow \Fil_{i+1}(D) 
\rightarrow \Fil_{i+1}(D)/\Fil_i(D) \longrightarrow 0$$
n'est pas scind\'ee dans $\fg/L$. Si $(D,\Fil_\bullet(D))$ est un $\fg$-module triangulin sur $\Ro_A$, on dira que 
$D$ est {\it partout non scind\'e} si pour tout $x \in {\rm Sp}(A)$, $(D_x,\Fil_\bullet(D)_x)$ est non scind\'ee. 
Si un $\fg$-module $D$ triangulin sur $\Ro_A$ est partout non scind\'e et si $B$ est une $A$-alg\`ebre affino\"ide, 
alors $D \widehat{\otimes}_A B$ est aussi partout non scind\'e vu comme $\fg$-module triangulin sur $\Ro_B$. On
dispose donc d'un sous-foncteur $$F_d^{\rm ns} \subset F_d$$ 
param\'etrant les $\fg$-modules triangulins r\'eguliers partout non scind\'es, et idem pour $F_d^{\square,{\rm ns}} \subset F_d^\square$. L'oubli de la rigidification d\'efinit un morphisme de foncteurs $$\eta : F_d^\square \rightarrow F_d$$ qui est surjectif
sur les points.  Soit $G_d=\mathbb{G}_m^d/\mathbb{G}_m$ (plong\'e
diagonalement) vu comme tore rigide analytique sur $\Q_p$.  On dispose enfin 
d'une action de $G_d$ sur $F_d^\square$ agissant sur les rigidifications :
$(a_i)\cdot [(D,\Fil_\bullet(D),(\nu_i))]:=[(D,\Fil_\bullet(D),(a_i \cdot
\nu_i))]$.  Cette action se factorise bien par $G$ car si $a \in A^\ast$ et
$x$ est un repr\'esentant d'une classe de $F_d^\square(A)$ alors
$(a,a,\dots,a)\cdot x$ est \'equivalent \`a $x$ via la multiplication par
$a$. Elle pr\'eserve $F_d^{\square,{\rm ns}}$.

\begin{lemma}\label{rigidite2} \begin{itemize}
\item[(i)] Si $(D_A,\Fil_\bullet)$ est un $\fg$-module triangulin r\'egulier
sur $\Ro_A$ qui est partout non-scind\'e, alors ses auto-\'equivalences sont les homoth\'eties
$A^\times$.\ps
\item[(ii)] $F_d^{\square,{\rm ns}}$ est repr\'esent\'e par un ouvert Zariski de $F_d^\square$.\ps
\item[(iii)] Pour toute alg\`ebre affino\"ide $A$, $G_d(A)$ agit librement sur
$F_d^{\square,{\rm ns}}(A)$ et l'application $\eta(A): F_d^{\square,{\rm
ns}}(A)/G_d(A) \longrightarrow F_d^{\rm ns}(A)$ est bijective.
\end{itemize}
\end{lemma}

\begin{pf} V\'erifions le (i). Nous allons montrer plus g\'en\'eralement que $$\End_{\fg_A}((D_A,\Fil_\bullet))=A.$$ Quitte \`a remplacer $A$ par $A/I$ pour un id\'eal $I$ de
codimension finie, on peut supposer que $A$ est artinien d'apr\`es le th\'eor\`eme
d'intersection de Krull. Dans ce cas, une r\'ecurrence sur la longueur de $A$ permet  de supposer que $A=L$ est un corps. Dans ce cas, on proc\`ede par r\'ecurrence sur $d$. Quand $d=1$ cela vient de ce que $H^0(\Ro_L)=L$. Pour $d\geq 1$, remarquons que si $(D,\Fil_\bullet)$ est non scind\'e, il en va de m\^eme de $(\Fil_{d-1}(D),\Fil_\bullet)$. Ainsi, un endomorphisme de $D$ pr\'eservant sa filtration agit par une homoth\'etie sur $\Fil_{d-1}$ et sur le quotient $D/\Fil_{d-1}$. Notons qu'il agit par $0$ sur ces deux $\fg$-modules si et seulement si il provient d'un morphisme $D/\Fil_{d-1}(D) \rightarrow \Fil_{d-1}(D)$ dans $\fg/L$, auquel cas il est en fait nul car si $\delta_i$ est le param\`etre de $D$ alors $H^0(\Fil_{d-1}(D)(\delta_d^{-1}))=0$ par l'hypoth\`ese de r\'egularit\'e. Nous avons donc montr\'e que le morphisme de $L$-alg\`ebres 
$$\alpha : \End_{\fg/L}((D,\Fil_\bullet)) \longrightarrow \End_{\fg/L}((\Fil_{d-1}(D),\Fil_\bullet))
\times \End(D/\Fil_{d-1}(D)) \isomo L \times L $$
est injectif. En particulier, si $\End_{\fg/L}((D,\Fil_\bullet))$ n'est pas r\'eduit aux homoth\'eties alors $\alpha$ est bijectif : il existe donc un endomorphisme idempotent de $D$ valant l'identit\'e sur $\Fil_{d-1}$ et $0$ sur $D/\Fil_{d-1}$, ce qui contredit le fait que $0 \rightarrow \Fil_{d-1}(D) \rightarrow D \rightarrow D/\Fil_{d-1}(D) \rightarrow 0$ est non scind\'ee.\ps
V\'erifions le (ii) par r\'ecurrence sur $d$. Si $d=1$, $F_d^{\square,{\rm ns}}=F_d^\square$ et il n'y a rien 
\`a d\'emontrer. En g\'en\'eral nous avons vu que $S_d^\square$ est un certain fibr\'e vectoriel de rang $d-1$ sur un ouvert $U_d \subset S_{d-1}^\square \times \cT$. Ce point de vue fait appara\^itre $F_d^{\square,{\rm ns}}$ comme l'ouvert du compl\'ementaire de la section nulle de ce fibr\'e pris au dessus de $F_{d-1}^{\square,{\rm ns}} \times \cT$, d'o\`u le r\'esultat.\ps

V\'erifions le (iii). Si $(a_i)(D,\Fil_\bullet,\nu)$ est \'equivalent \`a $(D,\Fil_\bullet,\nu)$ alors $(D,\Fil_\bullet)$ admet un automorphisme agissant sur chaque $\Fil_i(D)/\Fil_{i-1}(D)$ par le scalaire 
$a_i$. Comme les seuls automorphismes sont des homoth\'eties par le (i) il vient que tous les $a_i$ sont \'egaux : l'action de l'\'enonc\'e est libre. Le second point du (iii) vient de ce que les automorphismes dans $\fg/A$ de $\Ro_A(\delta)$ sont les homoth\'eties $A^\times$.
\end{pf}

Nous n'avons pas trouv\'e de r\'ef\'erences pour l'existence d'un quotient pour une action libre d'un tore sur un espace analytique. Ceci, combin\'e au lemme ci-dessus, entra\^inerait que le faisceau pour la topologie de Tate $\cF_d^{\rm ns}$ associ\'e \`a $F_d^{\rm ns}$ est repr\'esentable. Concr\`etement, $\cF_d^{\rm ns}(A)$ est la limite inductive (qui est filtrante injective) sur tous les  recouvrements finis de ${\rm Sp}(A)$ par des ouverts affino\"ides $U_i$ des noyaux des $\prod_i F_d^{\rm ns}(\OO(U_i)) \rightrightarrows \prod_{i,j} F_d^{\rm ns}(\OO(U_i\cap U_j))$.

Nous nous contenterons ici de traiter le cas $d=2$, pour lequel le probl\`eme se r\'esoud ais\'ement.

\begin{prop} $F_2^{\rm ns}$ est repr\'esent\'e par $\cT_2^{\rm reg}$.
\end{prop}
\begin{pf} En effet, le morphisme param\`etre $F_2^{\rm ns} \rightarrow \cT_2^{\rm reg}$ est un isomorphisme : si $\delta=(\delta_1,\delta_2) \in \cT_2^{\rm reg}(A)$ est donn\'e, alors on a vu que $H^1(\Ro_A(\delta_1\delta_2^{-1}))$ est libre de rang $1$ sur $A$, donc il existe un et un seul \'element de $F_2^{\rm ns}(A)$ de param\`etre $\delta$.
\end{pf}
\newcommand{\cTt}{{\widetilde{\cT}}}
Il se trouve que dans ce cas nous pouvons d\'ecrire aussi ce qui se passe au
voisinage des points non r\'eguliers. Pour tout $i \geq 0$
entier, notons $F_i,F'_i \subset \cT^2$ les ferm\'es d\'efinis respectivement
par les \'equations $\delta_1\delta_2^{-1}=\chi x^i$ et $\delta_1
\delta_2^{-1} = x^{-i}$. Tous ces ferm\'es sont deux \`a deux disjoints et
chaque ouvert affino\"ide de $\cT^2$ ne rencontre qu'un nombre fini d'entre
eux. On d\'esigne par $F$ la r\'eunion des $F_i$ et $F'$ celle des $F'_i$ :
ce sont encore des ferm\'es de $\cT^2$, et on a $\cT_2^{\rm reg} \coprod F
\coprod F' = \cT^2$. Soit $\pi : \cTt_2
\rightarrow \cT^2\backslash F'$ l'\'eclat\'e de $\cT^2\backslash F'$ le long de $F$. Le r\'esultat suivant confirme
l'intuition de Colmez dans~\cite{colmeztri} selon laquelle $\cTt_2$ est l'espace de module grossier des 
$\fg$-modules triangulins partout non scind\'es
de rang $2$ de param\`etre dans $\cT^2\backslash F'$ (sur l'ouvert $\cT_2^{\rm reg}$ c'est m\^eme un espace de module
fin par le r\'esultat pr\'ec\'edent).

\begin{prop} \label{thmdegal2}Pour toute extension finie $L/\Q_p$, il existe une bijection
canonique entre $\cTt_2(L)$ et l'ensemble des classes d'isomorphie de
$\fg$-modules triangulins sur $\Ro_L$ qui sont de rang $2$, non scind\'e, et
de param\`etre dans $\cT^2(L)\backslash F'(L)$, l'application $\pi$ donnant le param\`etre associ\'e. 

 De plus, il existe un recouvrement affino\"ide admissible $(U_i)$ de $\cTt_2$, et pour
chaque $i$ un $\fg$-module triangulin $D_i$ de rang $2$ sur $\Ro_{U_i}$ et de param\`etre
$\pi_{|U_i}$, tels que pour toute extension $L/\Q_p$ finie et tout $x \in
U_i(L)$, $(D_i)_x$ est isomorphe au $\fg$-module triangulin sur $\Ro_L$
associ\'e \`a $x$ par la bijection pr\'ec\'edente.
\end{prop}

La premi\`ere partie du th\'eor\`eme est d\^ue \`a Colmez :
c'est son calcul de $H^1(\Ro_L(\delta))$ pour l'ouvert $\cT_2^{\rm reg}$,
combin\'e \`a sa formule pour l'invariant $L$ (\cite{colmezinvL}) au voisinage des points de la forme $x \mapsto \chi x^i$ avec $i\geq
1$. Colmez a aussi d\'emontr\'e une version faible de la seconde partie au
voisinage de tout $x \in \cT_2^{\rm reg}$ qui
est de plus {\it $p$-r\'egulier} au sens que $\delta_1\delta_2^{-1}(p)
\notin p^\Z$.

\begin{pf} Au dessus de $\cT_2^{\rm reg}$, le th\'eor\`eme est un cas
particulier de la proposition pr\'ec\'edente. Quitte \`a tordre par une
famille de $\fg$-modules de rang $1$, on peut donc supposer que l'on
se place dans un voisinage ouvert affino\"ide $U \subset \cT$ contenant
un unique point $u \in U(\Q_p)$ tel que $\delta(u)=\chi x^i$ est non r\'egulier, et que
l'on s'int\'eresse \`a l'\'eclat\'e $\pi : \widetilde{U} \rightarrow U$ en
$u$. Dans ce cas, on dispose d'une identification naturelle donn\'ee par le th\'eor\`eme~\ref{h1universel} $$\mu : \pi^{-1}(u) \isomo {\mathbb
P}(H^1(\Ro_{\Q_p}(\chi x^i))).$$
Ce th\'eor\`eme assure aussi que $H^1(\Ro_U(\delta_U))$ est canoniquement isomorphe \`a $\OO(U)$, et que si $m \subset \OO(U)$ d\'esigne l'id\'eal des fonctions s'annulant en $u$, alors l'inclusion induit une \'egalit\'e $$H^1(m\Ro_U(\delta_U))=H^1(\Ro_U(\delta_U))=\OO(U).$$ Soit $D_U=m\Ro_U \oplus \Ro_U$ un $\fg$-module sur $\Ro_U$ dont la classe dans $H^1(m\Ro_U(\delta_U))$ en est un $\OO(U)$-g\'en\'erateur. Soit $U_i$ un recouvrement fini de $\widetilde{U}$ par
des ouverts affino\"ides sur chacun desquels $m\OO(U_i) \subset \OO(U_i)$ est libre de rang
$1$, disons engendr\'e par l'\'el\'ement $f_i$. Le transform\'e strict de $D_U$ sur $\Ro_{U_i}$, 
qui est aussi le quotient de $D_U \otimes_{\Ro_U} \Ro_{U_i}$ par sa $f_i$-torsion, est un $\fg$-
module sur $\Ro_{U_i}$ de $\Ro_{U_i}$-module sous-jacent $m\Ro_{U_i}\oplus \Ro_{U_i}=f_i\Ro_{U_i}\oplus \Ro_{U_i}$ : 
il est bien libre de rang $2$. Par construction de l'\'eclat\'e, le choix d'un $z \in (\pi^{-1}(u)\cap U_i)(L)$ d\'efinit 
un $\OO(U)$-morphisme surjectif $$(m/m^2)\otimes_{\Q_p} L \rightarrow f_i\OO(U_i)\otimes_{\OO(U_i)} L,$$ soit encore \`a 
un vecteur tangent $v_z \in T_u \otimes_{\Q_p} L$, tout vecteur tangent s'obtenant ainsi pour un certain $i$. Un tel 
choix d\'efinit donc un $\OO(U)$-morphisme surjectif $$D_U\otimes_{\OO(U)}L \rightarrow (D_i)_z.$$ La preuve du (ii) 
du th\'eor\`eme~\ref{h1universel} dit exactement que l'image de ce morphisme a pour classe dans $H^1(\Ro_{L}(\chi x^{-i}))$ 
l'\'element $\mu(v_z)$ associ\'e \`a l'isomorphisme $\mu(L) : {\mathbb P}(T_u)(L) \isomo {\mathbb P}(H^1(\Ro_L(\chi x^i)))$. 
\end{pf}

\begin{remark} {\it La situation est diff\'erente au voisinage des points dans $F'$. En effet, n\'egligeons les torsions en 
nous pla\c{c}ant dans un voisinage ouvert affino\"ide $U \subset \cT$ du point
$\delta=x^{-i}$. D'apr\`es le
th\'eor\`eme~\ref{h1universel}, pour tout \'el\'ement $E \in
H^1(\Ro(\delta))$ on peut trouver un $\fg$-module triangulin $D$ sur $\Ro_U$
de param\`etre $(\delta_U,1)$ tel que pour tout $z \in \cT^{\rm reg}\cap U$
son \'evaluation $D_z$ est non scind\'ee, et dont l'\'evaluation en $\delta$
est exactement $E$. On pourrait penser aller plus loin en
introduisant ici aussi l'\'eclat\'e de $U$ en $\delta$. Il n'y cependant pas
de mani\`ere naturelle de construire de famille de $\fg$-modules sur cet
\'eclat\'e. Disons simplement qu'un indice de ceci est que le th\'eor\`eme~\ref{h1universel} identifie 
canoniquement $T_\delta$ avec le {\rm dual} de l'espace vectoriel 
$H^1(\Ro_{\Q_p}(x^{-i}))$, plut\^ot qu'avec ce dernier. }
\end{remark}

Ainsi qu'il l'est expliqu\'e dans~\cite{berch}, notons que la proposition~\ref{thmdegal2}
entra\"ine la :

\begin{prop} La conjecture 5.1 de \cite{berch} est vraie : les
repr\'esentations potentiellement triangulines de dimension $2$ forment une
partie fine de la vari\'et\'e des caract\`eres $p$-adiques de dimension $2$ de ${\rm
Gal}(\Qpb/\Qp)$.
\end{prop}


Terminons ce paragraphe par une comparaison entre les foncteurs d\'efinis ici et les "foncteurs de d\'eformations triangulines" consid\'er\'es dans \cite[\S 2.5]{bch} et \cite[\S 3]{chU3}. Soit $F : \Aff \rightarrow {\rm Ens}$ un foncteur quelconque, $L/\Q_p$ une extension finie et $x \in F(L)$. Soit $\mathcal{C}$ la cat\'egorie des $L$-alg\`ebres locales artiniennes de corps r\'esiduel $L$ : c'est une sous-cat\'egorie pleine de celle des affino\"ides sur $L$. Le {\it compl\'et\'e formel} de $F$ en $x$ est le foncteur $\widehat{F}_x : \mathcal{C} \rightarrow {\rm Ens}$ d\'efini comme suit : pour tout objet $A$ de $\mathcal{C}$, $\widehat{F}_x(A) \subset F(A)$ est le sous-ensemble des \'el\'ements dont l'image dans $F(L)$ par l'unique morphisme $A \rightarrow L$ est l'\'el\'ement $x$ (c'est donc un sous-foncteur de $F_{|\mathcal{C}}$). Quand $F$ est repr\'esent\'e par un espace analytique
$Z$ sur $\Q_p$, $\widehat{F}_x$ est pro-repr\'esent\'e par $\widehat{\OO_{Z,x}}\otimes_{k(x)} L$.

\begin{lemma}\label{sousfonct}  Soit $x=[(D,\Fil_\bullet)] \in F_d(L)$. Si $D$ est non-scind\'e, et plus g\'en\'eralement si $\End_{\fg/L}((D,\Fil_\bullet))=L$, alors $\widehat{(F_d)}_x$ est canoniquement isomorphe au foncteur des d\'eformations triangulines de $(D,\Fil_\bullet)$ au sens de~\cite[\S 2.5]{bch}. De plus, $\widehat{(F_d)}_x$ est pro-repr\'esentable.
\end{lemma}
 
\begin{pf} En effet, le foncteur $\cX_{D,\Fil_\bullet}$ des d\'eformations triangulines de $(D,\Fil_\bullet)$ d\'efini {\it loc. cit.} param\`etre les classes d'isomorphismes de triplets $(D_A,\Fil_\bullet,\pi)$ o\`u $\pi : D_A \otimes_A L \rightarrow D$ est un isomorphisme dans $\fg/A$ envoyant $\Fil_i(D_A)$ sur $\Fil_i(D)$. L'association $(D_A,\Fil_\bullet,\pi) \mapsto [(D_A,\Fil_\bullet)]$ d\'efinit un morphisme de foncteurs $\cX_{D,\Fil_\bullet} \rightarrow \widehat{(F_d)}_x$ qui est surjectif sur les points. Pour l'injectivit\'e, il faut remarquer que si ${\rm Aut}_{\fg/L}((D,\Fil_\bullet))=L^\times$ (les homoth\'eties), alors pour tout $(D_A,\Fil_\bullet)$ dont la classe est dans $\widehat{(F_d)}_x$, alors l'application naturelle $${\rm Aut}_{\fg/A}(D_A,\Fil_\bullet) \rightarrow {\rm Aut}_{\fg/L}(D,\Fil_\bullet)$$
est surjective, car le terme de gauche contient $A^\times$. L'assertion de pro-repr\'esentabilit\'e est \cite[Prop. 3.4]{chU3}.
\end{pf}

\subsection{Densit\'e des $\fg$-modules cristallins dans $S_d^\square$.} Si $L$ est une extension finie de $\Q_p$ et $D$ un $\fg$-module sur $\Ro_L$, 
on dit que $D$ est {\it cristallin} si le $L$-espace vectoriel 
$$\cDc(D):=(D[1/t])^\Gamma$$ 
est de dimension ${\rm rg}_{\Ro_L}(D)$. Nous renvoyons 
\`a~\cite[\S 2.2.7]{bch} pour une discussion de cette d\'efinition, 
principalement motiv\'ee par des travaux de Berger : si $D$ est \'etale
alors il est cristallin si et seulement si sa repr\'esentation galoisienne
$V$ associ\'ee l'est, auquel cas $\cDc(D)$ est canoniquement isomorphe \`a
$\Dc(V)$ comme $L[\varphi]$-module filtr\'e (nous ne donnerons pas ici la
recette de la filtration naturelle sur $\cDc(D)$). Bien entendu, 
un point $x \in \cS_d^\square(L)$ est dit cristallin si le $\fg$-module 
triangulin $D_x$ sur $\Ro_L$ qui lui est associ\'e l'est. Dans ce cas, le 
param\`etre $(\delta_i)$ de $D_x$ est {\it alg\'ebrique} : pour tout $i$, il existe $k_i \in \Z$ tel que
$\delta_i(\gamma)={\gamma}^{k_i}$ pour tout $\gamma \in \Gamma$ (voir par
exemple \cite[prop. 2.4.1]{bch}). \ps\ps

On rappelle qu'une partie $A$ d'un espace rigide $Y$ est dite Zariski-dense si le 
seul ferm\'e analytique global r\'eduit de $Y$ contenant $A$ est la nilr\'eduction de $Y$. 
Si $Y$ est affino\"ide il est \'equivalent de demander que $A$ est Zariski-dense dans 
${\rm Spec}(\OO(Y))$. On dit de plus que la partie $A$ s'accumule en une partie $B \subset Y$ 
si tout \'el\'ement de $B$ admet une base de voisinages ouverts affino\"ides $U_i$ tels que 
$A\cap U_i$ est Zariski-dense dans $U_i$. On dit que $A$ est d'accumulation si $A$ s'accumule 
en $A$. On \'etend ces d\'efinitions \`a des parties de $Y(\Qpb):=\bigcup_L Y(L)$ (la r\'eunion 
portant sur les sous-extensions finies) en consid\'erant l'ensemble des poins ferm\'es sous-jacents. 
Comme exemple typique, remarquons que $\N^d$ est Zariski-dense et d'accumulation dans $\AAA^d$.\ps

Nous renvoyons \`a \cite{conradirr} pour les g\'en\'eralit\'es sur les composantes irr\'eductibles 
des espaces rigides analytiques $p$-adiques.

\begin{thm}\label{densecrisdanstri} Pour chaque extension finie $L/\Q_p$, l'ensemble des points cristallins de 
$\cS_d^\square(L)$ est Zariski-dense et s'acculume en chaque point de
$\cS_d^\square(L)$ de param\`etre alg\'ebrique.
\end{thm}

\begin{pf} Remarquons que $\cT_d^{\rm reg}$ est irr\'eductible, comme ouvert Zariski de l'espace 
irr\'eductible $\cT^d$. Comme un fibr\'e vectoriel sur une base irr\'eductible lisse est irr\'eductible 
lisse, $\cS_d^{\square}$ est irr\'eductible lisse par construction. Pour d\'emontrer la densit\'e Zariski 
de $\cS_d^\square(\Q_p)$,  il suffit donc de d\'emontrer qu'il est non vide ainsi que la propri\'et\'e 
d'accumulation. \ps
Soit $L$ une extension finie de $\Q_p$. Consid\'erons $A_d(L) \subset \cT_d^{\rm reg}(L)$ l'ensemble des 
points param\'etrant les $(\delta_i) \in \cT(L)^d$ tels que : 

(a) $\delta_i(p)/\delta_j(p) \neq p^{\pm 1}$ pour tout $i<j$,

(b) il existe une suite d'entiers $k_i \in \Z$ v\'erifiant : $\delta_i(\gamma)=\gamma^{-k_i}$ pour tout $\gamma \in \Gamma$ et tout
$i=1,\dots,d$,

(c) la suite $k_i$ est strictement croissante : $k_1 < k_2 < \dots < k_d$.

On note aussi $B_d(L) \subset \cT_d^{\rm reg}(L)$ l'ensemble des points satisfaisant uniquement la condition
(b), c'est \`a dire les param\`etres alg\'ebriques. 
Il est \'evident que $A_d(L)$ est non-vide, que $A_d(L) \subset B_d(L)$, et que
$A_d(L)$ s'accumule en $B_d(L)$ dans $\cT_d^{\rm reg}$. 

\begin{lemma}\label{critcris} Un $\fg$-module triangulin sur $\Ro_L$ qui a son param\`etre dans
$A_d(L)$ est cristallin.
\end{lemma}

\begin{pf} C'est un cas particulier de~\cite[Prop. 2.3.4]{bch}, largement pr\'ecis\'ee dans \cite{benois}, 
qui repose de mani\`ere essentielle sur des r\'esultats de Berger \cite{berger1} \cite{berger2}. En effet, 
si $D$ est comme dans l'\'enonc\'e alors il est de De Rham par cette proposition
(\`a cause de la condition (c)), et donc potentiellement 
semi-stable par le th\'eor\`eme de Berger. Comme $D$ est extension successive de $\fg$-modules cristallins par 
hypoth\`ese, et comme la formation du ${\mathcal D}_{\rm pst}$ est exacte sur les $\fg$-modules potentiellement semi-stables
d'apr\`es~\cite[Prop. 1.2.9]{benois}, il vient que $D$ est semi-stable. 
La condition (a) force alors $D$ \`a \^etre cristallin. 
\end{pf}

Retournons \`a la preuve du th\'eor\`eme~\ref{densecrisdanstri}. Supposons
donc $x \in S_d^\square(L)$ cristallin, ou plus g\'en\'eralement
tel que $\delta(x) \in B_d(L)$. Soit $U$ un voisinage ouvert affino\"ide de
$x$ dans $S_d^\square$, ainsi que $\Omega$ et $\iota$, comme dans le
corollaire~\ref{corfamcoleman}. Soit $(\Omega_i)_{i\in I}$ une base de voisinages ouverts affino\"ides de $\delta(x)$ dans $\cT_d^{\rm reg}$ 
dans lesquels $A_d(L)$ est Zariski-dense. Quand $i$ parcourt $I$ et $V$ les ouverts affino\"ides de 
${\mathbb{D}}^{\frac{d(d-1)}{2}}$, les $U_{i,V}:=\iota^{-1}(\Omega_i \times V)$ forment une base de voisinages 
ouverts affino\"ides de $x$ dans $\cS_d^\square$. Comme $A_d(L)$ est Zariski-dense dans chaque $\Omega_i$, 
il en va de m\^eme de\footnote{Il est imm\'ediat de v\'erifier que si $U$ et $V$ sont des affino\"ides, 
et si $A \subset U$ et $B \subset V$ des parties Zariski-denses, alors $A \times B$ est Zariski-dense 
dans $U \times V$.} $\iota^{-1}((A_d(L)\cap \Omega_i) \times V)$ dans $U_{i,V}$, qui est constitu\'e de 
points cristallins d'apr\`es le lemme~\ref{critcris}. 
\end{pf}

\subsection{La famille de repr\'esentations galoisiennes sur le lieu \'etale}
Pour terminer ce chapitre, consid\'erons le {\it sous-ensemble} $$S_d^{\square,0} \subset
S_d^\square$$
constitu\'e des points $x \in S_d^{\square}$ tels que le $\fg$-module $D_x$ sur
$\Ro_{k(x)}$ associ\'e est \'etale. Si $x \in S_d^{\square, 0}$ on
d\'esigne par $V_x$ la repr\'esentation de ${\rm Gal}(\Qpb/\Q_p)$ telle que
$\Drig(V_x) \simeq D_x$.

\begin{prop}\label{conskliu} Pour chaque $x \in S_d^\square$, il existe un voisinage ouvert affino\"ide
$\Omega$ de $x$ dans $S_d^\square$, un mod\`ele $\cA \subset \OO(\Omega)$, et un
$\cA$-module libre $M$ de rang
$d$ muni d'une application $\cA$-lin\'eaire continue de ${\rm
Gal}(\Qpb/\Q_p)$ tels que : \begin{itemize}
\item[(i)] $\Omega \subset S_d^{\square,0}$,\par
\item[(ii)] pour tout morphisme $Z \rightarrow \Omega$ avec $Z$ affino\"ide, 
le $\fg$-module $\Drig(M \otimes_{\cA} \OO(Z))$ sur $\Ro_{Z}$ d\'efini par 
Berger et Colmez est isomorphe au $\fg$-module triangulin
rigidifi\'e universel d\'eduit du morphisme donn\'e $Z \rightarrow S_d^\square$.
\end{itemize}
Le $\OO(\Omega)[{\rm Gal}(\Qpb/\Qp)]$-module $M \otimes_{\cA} \OO(\Omega)$ est unique \`a isomorphisme
pr\`es pour la propri\'et\'e (ii).
\end{prop}

\begin{pf} L'existence de $\cA$ et $M$ satisfaisant (i) et (ii) pour
$Z=\Omega$, ainsi que l'assertion d'unicit\'e, est un cas particulier du th\'eor\`eme 
\cite[Thm. 0.2]{kedliu} de Kedlaya-Liu appliqu\'e \`a la famille universelle de $\fg$-modules sur $S_d^\square$. 
V\'erifions le (ii) pour $Z \rightarrow \Omega$ quelconque. Il s'agit de voir
qu'avec la d\'efinition du $\Drig$ d'une famille de Berger-Colmez
l'application naturelle 
\begin{equation}\label{isombcol}\Drig(M \otimes_{\cA} \OO(\Omega))
\otimes_{\Ro_\Omega} \Ro_Z \rightarrow \Drig(M \otimes_{\cA} \OO(Z))
\end{equation}
est un isomorphisme. Le terme de droite a bien un sens car la
repr\'esentation $M \otimes_{\cA} \OO(Z)$ admet un $\cB$-r\'eseau libre
stable pour tout mod\`ele $\cB$ de $\OO(Z)$ contenant l'image de $\cA$ (de
tels mod\`eles existent et on en fixe un). L'isomorphisme ci-dessus d\'ecoule alors de l'assertion d'unicit\'e de la proposition \cite[prop.
4.2.8]{bergercolmez} comme dans la preuve de leur th\'eor\`eme 4.2.9. : les
d\'etails sont sans difficult\'e et laiss\'es au lecteur. 
\end{pf}

En particulier, $S_d^{\square,0} \subset S_d^\square$ est un ouvert pour la topologie
na\"ive. Notons qu'\`a priori, 
ce r\'esultat ne conf\`ere pas \`a $S_d^{\square,0}$ de structure naturelle d'espace rigide analytique. (On pourrait 
cependant choisir une \'ecriture de $S_d^{\square,0}$ comme r\'eunion disjointe d'ouverts affino\"ides $\mathcal{U}=
\{\Omega_i, i \in I\}$ de $S_d^{\square}$, ce qui est loisible (et l'on peut m\^eme demander que chaque $\Omega_i$ 
satisfasse les conclusions de la proposition~\ref{conskliu}), et d\'ecr\^eter que $\mathcal{U}$ est un recouvrement 
admissible de $S_d^{\square,0}$. Pour chaque telle structure l'inclusion $S_d^{\square,0} \rightarrow S_d^{\square}$ 
est une immersion ouverte.) On a en revanche une notion canonique de fonction {\it localement 
analytique $f: S_d^{\square,0} \rightarrow Y$} vers un espace analytique $Y$ quelconque : nous 
entendrons par l\`a une application telle que pour chaque $x \in S_d^{\square,0}$ la restriction 
de $f$ a un voisinage affino\"ide suffisament petit de $x$ est une fonction rigide analytique 
(noter que $S_d^\square$ est r\'eduit). Il r\'esulte de la proposition ci-dessus que 
les applications $x \mapsto {\rm trace}(g | V_x)$, pour $g \in  {\rm
Gal}(\Qpb/\Qp)$, sont localement analytiques sur 
$S_d^{\square,0}$, et simultan\'ement analytiques sur un
voisinage de chaque $x$, ce qui implique le :

\begin{cor} Soit $\mathcal{X}_d$ la vari\'et\'e des caract\`eres $p$-adiques de dimension $d$ de ${\rm Gal}(\Qpb/\Qp)$. Il existe une unique application localement analytique $$f : S_d^{\square, 0} \rightarrow \mathcal{X}_d$$
associant \`a tout $x \in S_d^{\square,0}$ la semi-simplification de la repr\'esentation $V_x$.
\end{cor}

L'image de cette application est la foug\`ere infinie r\'eguli\`ere.  Nous renvoyons \`a \cite{chdet} pour la d\'efinition de 
$\mathcal{X}_d$. Pour le lecteur peu familier avec cette th\'eorie, donnons une version plus concr\`ete concernant la $\rhob$-composante 
connexe $\mathcal{X}_d(\rhob)\simeq X \subset \mathcal{X}_d$ o\`u $\rhob$ est fix\'ee comme dans l'introduction
(sans supposer n\'ecessairement $\rhob \not\simeq \rhob(1)$). Pour cela nous allons 
d\'efinir ind\'ependamment l'ensemble $$S_d^{\square}(\rhob) \subset S_d^{\square,0}$$
pull-back par l'application $f$ du corrollaire ci-dessus de l'ouvert $\cX_d(\rhob) \subset \cX_d$. Cela nous oblige \`a quelques sorites et rappels pr\'eliminaires sur la notion de {\it repr\'esentation r\'esiduelle associ\'ee \`a une famille de repr\'esentations} et  sur la propri\'et\'e universelle de $X$ (qui rappelons-le est d\'efini comme l'espace analytique associ\'e par Berthelot \`a la fibre g\'en\'erique de l'anneau de d\'eformation universelle de $\rhob$). Nous renvoyons \`a
\cite[\S 3]{chdet} pour plus de d\'etails. \ps

Si $Y$ est un affino\"ide, une {\it famille de repr\'esentations de $G:={\rm Gal}(\Qpb/\Q_p)$
param\'etr\'ee par $Y$} est la donn\'ee d'un $\OO(Y)$-module libre de rang fini muni d'une action $\OO(Y)$-lin\'eaire
continue de $G$. Pour $y \in Y$, on note
$M_y:=M\otimes_{\OO(Y)} k(y)$ l'\'evaluation de $M$ en $y$ et $\overline{M}_y$ la semi-simplifi\'e de la r\'eduction modulo $\pi_{k(y)}$ d'un
$\OO_{k(y)}$-r\'eseau stable par $G$ dans $M_y$ (elle est
donc \`a coefficients dans le corps fini $k_y$ r\'esiduel de $k(y)$). Si $Y$
est connexe, on peut montrer (\cite[Def. 3.11]{chdet}) qu'il existe une repr\'esentation semi-simple
continue $r : G \rightarrow \GL_m(\overline{\mathbb{F}}_p)$, 
et un morphisme d'alg\`ebres $W(\mathbb{F}) \rightarrow \OO(Y)$ o\`u $\mathbb{F} \subset
\overline{\mathbb{F}}_p$ d\'esigne le corps (fini) engendr\'e par les
coefficients des polyn\^omes caract\'eristiques des $r(g)$ pour $g \in G$, tels que pour tout $y \in Y$ alors $$\forall g \in G, \, \, \det(1-T g | \overline{M}_y) = \det(1-T r(g)) \in \mathbb{F}[T].$$ 
Cette identit\'e a un sens car $\mathbb{F} \subset k_y$ pour tout $y$. Notons que modifier le plongement $W(\mathbb{F}) \rightarrow \OO(Y)$ 
par le Frobenius de $W(\mathbb{F})$ revient \`a appliquer \`a $r$ le Frobenius sur les coefficients (ou son inverse). \`A cette ind\'etermination 
pr\`es, un r\'esultat standard d\^u \`a Brauer-Nesbitt implique que la classe d'isomorphisme de $r$ est uniquement d\'etermin\'ee par $M$. On dit 
alors que $r$ est la repr\'esentation r\'esiduelle de $M$. C'est un fait relativement formel mais important que {\it l'espace analytique $X$ 
sur $\Q_p$ (oubliant la $F$-structure) repr\'esente le foncteur de la cat\'egorie des affino\"ides vers les ensembles associant \`a $Y$ 
l'ensemble des classes de $\OO(Y)[G]$-isomorphie de familles de repr\'esentations de $G$ param\'etr\'ees par $Y$ dont la repr\'esentation 
r\'esiduelle est $\rhob$} (\cite[\S 3]{chdet}).\ps

Si $r : G \rightarrow \GL_d(\overline{\mathbb{F}}_p)$ est une repr\'esentation semi-simple continue quelconque, on d\'efinit enfin 
$$S_d^{\square}(r) \subset S_d^{\square,0}$$
comme l'ensemble des $x \in S_d^{\square,0}$ tels que $V_x$ a pour r\'eduction $r$. 

\begin{cor}\label{scholieuniv} ({\rm et scholie}) Pour tout $r$, $S_d^{\square}(r)$ est un ouvert
na\"if\footnote{C'est \`a dire une r\'eunion quelconque d'ouverts
affino\"ides.} de $S_d^{\square}$. Si $r=\rhob$ et $X$ sont comme dans
l'introduction\footnote{Avec \'eventuellement $\rhob \simeq \rhob(1)$.}, il existe une unique application localement analytique $$f : S_d^{\square}(\rhob) \rightarrow X$$
associant \`a tout $x \in S_d^{\square}(\rhob)$ la repr\'esentation $V_x$.
Plus pr\'ecis\'ement, elle a la propri\'et\'e que pour tout $x \in S_d^{\square}(\rhob)$, 
il existe $\Omega \subset S_d^{\square}(\rhob)$ un voisinage ouvert affino\"ide assez petit de $x$ dans
$S_d^{\square}$ tel que : \begin{itemize}
\item[(i)] $f$ est analytique sur $\Omega$, \ps
\item[(ii)] Pour tout morphisme $Z \rightarrow \Omega$ avec $Z$ affino\"ide,
le $\fg$-module associ\'e par Berger-Colmez \`a la famille d\'eduite du
morphisme $f : \Omega \rightarrow X$ soit isomorphe au $\fg$-module associ\'e 
au morphisme $Z \rightarrow S_d^\square$ par la propri\'et\'e universelle de ce dernier.
\end{itemize}
\end{cor}

\begin{pf} En effet, le premier point r\'esulte de la discussion ci-dessus et de la proposition~\ref{conskliu} en consid\'erant la composante connexe $\Omega_x$ contenant $x$ du $\Omega$ donn\'e par la proposition. La propri\'et\'e universelle de $X$ d\'efinit alors un unique morphisme analytique $\Omega_x \rightarrow X$ qui n'est autre que l'application de l'\'enonc\'e au niveau des points et le corollaire tout entier r\'esulte de la proposition~\ref{conskliu} (ii). 
\end{pf}

\begin{remark} Par d\'efinition, $S_d^{\square,0}$ est la r\'eunion disjointe des $S_d^{\square}(r)$ sur l'ensemble des 
$r : G \rightarrow \GL_d(\overline{\mathbb{F}}_p)$ semi-simples continus consid\'er\'es modulo isomorphisme et action du 
Frobenius sur les coefficients. En fait, $S_d^{\square}(r)$ est non vide pour tout $r$
: quand $r$ est irr\'eductible cela d\'ecoule de la proposition~\ref{crislift} (i) ci-dessous, le cas g\'en\'eral s'en d\'eduit
en consid\'erant des sommes directes. On verra m\^eme plus tard que $S_d^{\square}(r)$ contient des points
cristallins absolument irr\'eductibles. 
\end{remark}

\section{D\'emonstration du th\'eor\`eme $A$}

Reprenons les notations de l'introduction. 

\subsection{Premi\`eres r\'eductions}\label{reductions}

Commen\c{c}ons la d\'emonstration du th\'eor\`eme $A$.  Si $L$ est une extension finie de $\Q_p$ et si $V$ est une
$L$-repr\'esentation cristalline de ${\rm Gal}(\Qpb/\Qp)$, on dira que $V$
est {\it g\'en\'erique} si les conditions suivantes sont satisfaites
(\cite[\S 3.18]{chU3}) :\begin{itemize}\ps
\item[(G1)] $V$ a $d:=\dim_L(V)$ poids de Hodge-Tate\footnote{Il s'agit ici
de poids de Hodge-Tate relativement \`a $L$, c'est \`a dire les racines du
polyn\^ome de Sen vu comme \'el\'ement de $L[T]$.} distincts. \ps
\item[(G2)] Les valeurs propres $\varphi_1,\dots,\varphi_d$ du Frobenius
cristallin de $\Dc(V)$ dans $\overline{L}$ satisfont $\varphi_i\varphi_j^{-1} \notin p^\Z$ pour tout $i\neq j$.\ps
\item[(G3)] Si $S \subset \Dc(V)$ est un sous-$L$-espace vectoriel
$\varphi$-stable, alors $S$ est en somme directe avec le sous-espace de la
filtration de Hodge de $\Dc(V)$ dont la dimension est $d-\dim_L(S)$. 
\end{itemize}
\ps
On dira aussi que {\it $V$ est d\'eploy\'ee sur $L$} si les valeurs propres du Frobenius cristallin de
$\Dc(V)$ dans $\overline{L}$ sont toutes dans $L$.\ps

L'ingr\'edient principal de la d\'emonstration est alors le suivant. Notons
$T_y(Y)$ l'espace tangent de Zariski au point $y$ de l'espace rigide $Y$ (c'est un espace 
vectoriel de dimension finie sur le corps r\'esiduel $k(y)$). Pour une 
extension finie $L$ de $\Q_p$, on note aussi $T_y(Y)$ l'espace
$T_{|y|}(Y) \otimes_{k(y)} L$, $|y| \in Y$ d\'esignant le point ferm\'e
sous-jacent \`a $y$ et $k(y) \rightarrow L$ l'application structurale de $y$.

\begin{prop}\label{tggen}  Soient $x \in X(L)$ cristallin g\'en\'erique d\'eploy\'e sur
$L$, $U \subset X$ un voisinage ouvert affino\"ide de $x$ dans $X$, et $W \subset U$ l'adh\'erence
Zariski des points cristallins de $U(L)$. Alors $\dim_L T_x(W)=d^2+1$. 
\end{prop}

Nous reportons la preuve de cette proposition au \S~\ref{secfin}. Le second ingr\'edient est le suivant.

\begin{prop}\label{existcrisgen}\label{crislift}
\begin{itemize} 
\item[(i)] Il existe une extension $L/F$ de degr\'e $\leq d^2$ telle que 
$X(L)$ contienne des points cristallins d\'eploy\'es sur $L$ et dont tous les poids de Hodge-Tate sont disctincts.\ps
\item[(ii)] Si $x \in X(L)$ est cristallin d\'eploy\'e sur $L$, et si $\Omega \subset X$ est
un ouvert affino\"ide contenant $x$, alors $\Omega(L)$ contient un point cristallin
g\'en\'erique et d\'eploy\'e sur $L$.
\end{itemize}
\end{prop}

\begin{pf} Montrons le (i). Il est bien connu que si $F_d=W(\FF_{p^d})[1/p]$ est l'extension non
ramifi\'ee de degr\'ee $d$ de $\Q_p$, il existe un caract\`ere continu
$\eta : G_{F_d} \rightarrow \Fpb^\ast$ tel que $$\rhob \simeq {\rm
Ind}_{G_{F_d}}^{G_{\Q_p}} \eta,$$ 
par irr\'eductibilit\'e absolue de $\rhob$. Soit $\Sigma:={\rm Hom}_{\rm corps}(\FF_{p^d},\Fpb)$, cet ensemble \`a $d$
\'el\'ements s'identifie canoniquement par r\'eduction modulo $p$ \`a ${\rm Hom}_{\rm
corps}(F_d,W(\Fpb)[1/p])={\rm Hom}_{\rm ann}(\OO_{F_d},W(\Fpb))$. Si ${\rm rec} :
\widehat{F_d^\ast} \rightarrow G_{F_d}^{\rm ab}$ est
l'isomorphisme du corps de classe local, alors $\eta \,\,{\rm o}\,\, {\rm
rec}$ est n\'ecessairement (1) trivial sur $1+p\OO_{F_d}$, (2) de la forme $x \mapsto
\prod_{\sigma \in \Sigma} \sigma(x)^{a_\sigma}$ pour
certains entiers uniques $a_\sigma \in \{0,1,\dots,p-1\}$ sur
$\OO_{F_d}^\ast$, et (3) envoie $p$ sur un certain
element $\overline{\lambda} \in \Fpb^\ast$. Notons aussi que 
$$\det(X-\rhob_{|G_{F_d}^{\rm ab}}(p))=(X-\overline{\lambda})^d \in
\mathbb{F}[X]$$ assure que $\overline{\lambda} \in
\mathbb{F}$.\ps

Il vient que pour toute collection
d'entiers $\{b_\sigma, \sigma  \in \Sigma\}$ telle que $b_\sigma \equiv a_\sigma \bmod
p^d$ pour tout $\sigma$, et pour
tout $\lambda \in \OO_F^\ast$ tel que $\overline{\lambda}\equiv \lambda
\bmod p$, le caract\`ere $\widetilde{\eta} : G_{F_d} \rightarrow
(F_d \cdot F)^\ast$
d\'efini par la formule $$\forall u \in \OO_{F_d}^\ast, \widetilde{\eta} \,\,{\rm o}\,\, {\rm
rec}(p u)=\lambda \prod_{\sigma \in \Sigma} \sigma(u)^{b_\sigma}$$
rel\`eve $\eta$, et donc que $V:={\rm
Ind}_{G_{F_d}}^{G_{\Q_p}} \widetilde{\eta}$ rel\`eve
$\rhob$. Enfin, $V$ est cristalline car $F_d$ est non
ramifi\'ee et $\eta$ est cristallin d'apr\`es un r\'esultat de Fontaine : les poids de Hodge-Tate de $V$ sont les $b_\sigma$
et le polyn\^ome caract\'eristique du Frobenius cristallin sur $\Dc(V)$ est
$X^d-\lambda p^{\sum_\sigma a_\sigma}$. Le (i) s'en d\'eduit. \ps
	Le (ii) est une cons\'equence des r\'esultats de Kisin dans~\cite[\S
3]{kisinJAMS}, expliquons comment. Soit $x \in X(L)$ un point
cristallin \`a poids de Hodge-Tate distincts $\uk=(k_1 < \dots < k_d)$. Si
$V_R$ d\'esigne la d\'eformation universelle de $\rhob$ (l'anneau $R$
\'etant comme dans l'introduction), Kisin d\'emontre l'existence d'un
quotient $$R[1/p] \rightarrow R_{\uk}$$ dont les points dans toute
$F$-alg\`ebre de dimension finie $B$ param\`etrent exactement les morphismes $R[1/p]
\rightarrow B$ tels que $V_R \otimes_R B$ est cristalline de poids de
Hodge-Tate $k_1 < \dots < k_d$ (\cite[Thm. 2.5.5]{kisinJAMS}). Notons $X_\uk
\subset X$ le ferm\'e analytique de $X$ d\'efini par l'id\'eal noyau de
$R[1/p] \rightarrow R_{\uk}$. D'apr\`es Kisin~\cite[Thm. 3.3.4]{kisinJAMS}, $X_\uk$ est lisse
de dimension $\frac{d(d-1)}{2}+1$. Rappelons que c'est un fait g\'en\'eral que l'application
naturelle $R_\uk \rightarrow \OO(X_\uk)$ induit une bijection 
$\psi : X_\uk \longrightarrow {\rm Specmax}(R_{\uk})$ et des isomorphismes
$\widehat{(R_\uk)}_{\psi(x)} \rightarrow \OO_{X_\uk,x}$ sur les annaux
locaux compl\'et\'es pour tout $x \in
X_\uk$. Le (ii) d\'ecoule du lemme suivant.\end{pf}

\begin{lemma} L'ensemble des points $x \in X_{\uk}(L)$ tels que $\rho_x$ est g\'en\'erique est 
dense dans $X_\uk(L)$. De plus, si $x \in
X_\uk(L)$ est d\'eploy\'e sur $L$, il existe un ouvert affino\"ide $O$ de $x$ dans
$X_\uk$ tel que $O(L)$ est constitu\'e de points d\'eploy\'es sur $L$.
\end{lemma}
\newcommand{\PGL}{{\mathrm{PGL}}}
\begin{pf} En effet, Kisin d\'emontre {\it loc. cit.} 
l'existence d'un $\varphi$-module filtr\'e $D$ localement libre de rang $d$
sur $R_{\uk}$ tel que pour tout quotient artinien $R_{\uk} \rightarrow B$,
$D_{\rm cris}(V_R \otimes_R B)$ est isomorphe \`a $D \otimes_{R_{\uk}} B$, 
ce qui d\'etermine donc $V_B$ par l'\'equivalence de Fontaine. \ps
Soient $P_\varphi \in R_\uk[T]$ le poyn\^ome
caract\'eristique de $\varphi$ sur $D$ et $\cF$ la vari\'et\'e des drapeaux complets de
$F^d$, disons vue comme vari\'et\'e $F$-analytique. Remarquons d\'ej\`a que l'existence
m\^eme de $P_\varphi$, ainsi que le lemme de Krasner, impliquent la seconde partie du
lemme, concentrons nous donc sur la premi\`ere.\ps 

Notons que $\GL_d \times \cF$ est muni d'une action de $\PGL_d$
d\'efinie sur les points affino\"ides par $g\cdot (\varphi,{\bf F})=(g \varphi
g^{-1},g({\bf F}))$
(strictement il faudrait ins\'erer des recouvrements pffp ou m\^eme Zariski). 
Soient $x \in X_\uk(L)$ et $\Omega$ un voisinage ouvert
affino\"ide de $x$ dans $X_\uk$ assez petit de sorte que $D
\otimes_{R_\uk} \OO(\Omega)$ soit libre de rang $d$ sur $\OO(\Omega)$. Le
choix d'une base de ce dernier d\'efinit alors un morphisme analytique $\Omega
\rightarrow \GL_d \times \cF$ associ\'e \`a la matrice de $\varphi$ dans cette base et \`a la donn\'ee de
la filtration de $D$. Il sera commode de le modifier un peu en $$\mu :
\PGL_d \times \Omega \rightarrow \GL_d \times \cF$$
d\'efini sur les points par $(g,x) \mapsto g\cdot (\varphi(x),{\bf F}(x))$. La propri\'et\'e de $D$ vis-\`a-vis des
$F$-alg\`ebres artiniennes locales $B$ \'enonc\'ee plus haut a plusieurs cons\'equences. Appliqu\'ee aux $B$ qui sont
des corps elle entra\^ine que $\mu$ est injective. Si de plus un point $z \in \PGL_d \times \Omega$ est fix\'e, et appliqu\'ee aux
\'epaississement infinit\'esimaux de $z$, elle entra\^ine que $$\mu_z^{\ast} : 
\widehat{\OO}_{\GL_d \times \cF,\mu(z)} \rightarrow
\widehat{\OO}_{\PGL_d\times\Omega,z}$$
est un isomorphisme. En effet, c'est un fait g\'en\'eral que si $\theta : A_1 \rightarrow A_2$ est
un $F$-morphisme local entre des $F$-alg\`ebres locales noeth\'eriennes
compl\`etes telles que pour toute $F$-alg\`ebre locale de dimension finie
$B$ l'application naturelle $A_2(B) \rightarrow A_1(B)$ est bijective, alors
$\theta$ est un isomorphisme. Ce crit\`ere se v\'erifie bien ici car d'une
part pour tout tel $B$-point de $X_\uk$, les automorphismes de $V_R \otimes B$, ou ce qui
revient au m\^eme de $D \otimes_{R_{\uk}}B$, sont r\'eduits aux
scalaires $B^\times$, et d'autre part $\PGL_d(B)=\GL_d(B)/B^\ast$ car ${\rm
Pic}(B)=0$. Ainsi, $\mu$ est injectif et un isomorphisme
local formel en tout point de sa source. En particulier, le morphisme $\mu$
est plat. Comme sa source et son but sont quasi-compacts, un r\'esultat de
Bosch-L\"utckebohmert (existence de mod\`eles formels plats) assure que
$\mu(\Omega) \subset \GL_d \times \cF$ est un ouvert admissible (quasi-compact).  
Comme $\mu$ induit des isomorphismes en chaque point sur les corps
r\'esiduels, suffit de voir pour conclure que si $$(\varphi,{\bf F}) \in \GL_d(L) \times
\cF(L)$$ est donn\'e, et si $U$ en est un voisinage ouvert pour la
topologie $p$-adique, alors $U$ contient une paire $(\varphi',{\bf F}')$
o\`u $\varphi'$ est semi-simple \`a valeurs propres $\lambda_i \in
\overline{L}$ distinctes, telles que $\lambda_i/\lambda_j \notin p^\Z$ pour
$i\neq j$, et dont les espaces propres dans $\overline{L}^d$ sont en position
g\'en\'erale par rapport \`a ${\bf F'}$.

Choisir un $\varphi'$ proche de
$\varphi$ satisfaisant les deux premi\`eres conditions est imm\'ediat,
fixons le. Soit $V \subset \cF(L)$ un voisinage ouvert assez petit de ${\bf F}$ de sorte
que que $\{\varphi'\} \times V \subset U$. Notons que $V$ est dense dans
$\cF(\overline{L})$ pour la topologie de Zariski de ce dernier (qui est
irr\'eductible). Il suffit pour conclure de remarquer que l'ensemble des
drapeaux complets de $\overline{L}^d$ qui sont en position g\'en\'erale avec
un nombre fini fix\'e de drapeaux complets forment un ouvert Zariski :
c'est standard pour un drapeau ("grosse cellule de Bruhat"), et le cas g\'en\'eral s'en d\'eduit par
intersection finie. 
\end{pf}

Remarquons que les relevements cristallins exhib\'es au (i) ci-dessus ne sont pas
Zariski-denses dans $X$, car ils restent dans le sous-espace ferm\'e des
repr\'esentations induites d'un caract\`ere de $F_d$, dont on v\'erifie
facilement qu'il est de dimension $d < d^2+1$. \ps

D\'emontrons enfin le th\'eor\`eme $A$ de
l'introduction. Fixons une fois pour toutes un point cristallin $x \in
X(L)$ d\'eploy\'e sur $L$ donn\'e par la
proposition~\ref{existcrisgen} (i), ainsi qu'un ouvert affino\"ide $U
\subset X$ connexe et contenant $x$. Notons que $X$ \'etant normal, $U$ est
irr\'eductible. Soit $W \subset U$ le ferm\'e analytique r\'eduit obtenu en
prenant l'adh\'erence Zariski des points cristallins de $U(L)$. Les
affino\"ides \'etant excellents et de Jacobson (voir \cite[\S 1]{conradirr}), le
lieu r\'egulier de $W$ est un ouvert Zariski qui est dense dans $W$. En
particulier, il existe des points cristallins $y \in U(L)$ qui sont
r\'eguliers dans $W$. Quitte \`a remplacer $x$ par un tel point, nous
pouvons donc supposer que $x$ est r\'egulier dans $W$. Appliquant la
proposition~\ref{existcrisgen} (ii) \`a un voisinage ouvert affino\"ide 
$\Omega$ de $x$ dans $U$ qui soit assez petit de sorte que
$\Omega \cap W$ soit r\'egulier, on peut finalement supposer que $x \in W(L)$ est
cristallin g\'en\'erique, d\'eploy\'e sur $L$, et r\'egulier dans $W$. La proposition~\ref{tggen} assure alors que
$\dim W = \dim_F T_x(W) = d^2+1 = \dim X$, et donc que $W$ contient un
ouvert de $X$ contenant $x$. Comme $U$ est irr\'eductible, cet ouvert est
Zariski-dense dans $U$, ce qui conclut la preuve du th\'eor\`eme.\ps

\subsection{Preuve de la proposition~\ref{tggen}.}\label{secfin}
\newcommand{\cC}{\mathcal{C}}
 Fixons $L$ une extension finie de $\Q_p$ et $x \in X(L)$ tels que
$$V_0:=\rho_x$$
est cristalline g\'en\'erique et d\'eploy\'ee sur $L$. On notera
$k_1<\dots<k_d$ les poids de Hodge-Tate de $V_0$ rang\'es par ordre
croissant. (Nous prenons pour convention que $\Q_p(1)$ a pour poids de
Hodge-Tate $-1$.) On rappelle que le
$\fg$-module $$D_0:=\Drig(V_0)$$ est triangulin sur $\Ro_L$ d'exactement $d!$
mani\`eres diff\'erentes. Plus pr\'ecis\'ement, fixons $\Ref$ un {\it
raffinement} de $V_0$, c'est \`a dire un drapeau complet $L[\varphi]$-stable dans
$\Dc(V_0)$
$$\Ref=\left(\Ref_0=\{0\} \subsetneq \Ref_1 \subsetneq \cdots \subsetneq
\Ref_d=\Dc(V_0)\right).$$
Il est \'equivalent ici de se donner un ordre
$\varphi_1,\dots,\varphi_d$ sur les valeurs propres de $\varphi$ sur
$\Dc(V_0)$, $\varphi$ agissant sur $\Ref_i/\Ref_{i-1}$ par la
multiplication par $\varphi_i$. \`A chaque tel $\Ref$ est associ\'ee une
triangulation $$\mathcal{T}=(\Fil_i(D_0))$$ de $D_0$ d\'efinie comme suit. Nous avons d\'ej\`a dit que d'apr\`es
Berger le $L[\varphi]$-module $(D_0[1/t])^{\Gamma}$ est canoniquement
isomorphe \`a $\Dc(V_0)$, on pose alors $$\Fil_i(D_0)=(\Ro_L[1/t].\Ref_i)\cap
D_0.$$ L'application $\Ref \mapsto {\mathcal{T}}$ induit alors la bijection
cherch\'ee entre raffinements de $V_0$ et triangulations de $D_0$, d'apr\`es
\cite[\S 2.4]{bch}. Dans cette bijection, le param\`etre $(\delta_i)$ de
$\mathcal{T}$ est reli\'e \`a $\Ref$ par les formules $$\delta_i(p)=\varphi_i
p^{-k_i},  \, \, \, \delta_i(\gamma)=\gamma^{-k_i} \, \,\, \,\forall \gamma \in \Z_p^\ast,$$ (cela
d\'ecoule de \cite[prop. 2.4.1]{bch} et de l'hypoth\`ese (G3) de
g\'en\'ericit\'e). \ps

Fixons $\Ref$ un raffinement de $V_0$ et consid\'erons le $\fg$-module triangulin sur 
$\Ro_L$ associ\'e $(D_0,\cT)$. L'hypoth\`ese (G2) de g\'en\'ericit\'e entra\^ine qu'il 
est r\'egulier.  Fixons enfin une
rigidification $\nu_0$ de $(D_0,\cT)$, ce qui nous fournit donc un point $$x_0 =(D_0,\cT,\nu_0) \in S_d^\square(L).$$
Soient $\Omega \subset S_d^{\square}$ un voisinage ouvert affino\"ide de $x_0$ et $$f : \Omega \rightarrow X$$ comme 
dans le corollaire~\ref{scholieuniv}. Quitte \`a remplacer $\Omega$ par un voisinage plus petit, on peut supposer 
d'une part que $f(\Omega) \subset U$ (l'ouvert $U$ \'etant celui de l'\'enonc\'e de la proposition~\ref{tggen}), 
et d'autre part que les points cristallins sont Zariski-dense dans $\Omega(L)$ d'apr\`es le 
th\'eor\`eme~\ref{densecrisdanstri}. Si $W$ est l'adh\'erence Zariski des points cristallins de $U(L)$, 
l'analyticit\'e de $f$ entra\^ine que \begin{equation}\label{inclfond} f(\Omega) \subset
W.\end{equation}
D'autre part, si $L[\varepsilon]$ d\'esigne les nombres duaux sur $L$
(de sorte que $\varepsilon^2=0$) alors la fonction analytique $f$ induit une application tangente (voir la fin du~\S~\ref{diverstrucs}) 
$$df_{x_0} : \widehat{(F_d^\square)}_{x_0}(L[\varepsilon]) \longrightarrow
T_x(X)$$
qui n'est autre, d'apr\`es le corollaire~\ref{scholieuniv} appliqu\'e aux
morphismes 
${\rm Sp}(L[\varepsilon]) \rightarrow \Omega$ d'image $x_0$, o\`u ce qui
revient au m\^eme \`a tous les morphismes ${\rm Sp}(L[\varepsilon])
\rightarrow S_d^\square$ d'image $x_0$ car $\Omega \subset S_d^\square$ est
ouvert, que l'application associant \`a un 
$c=[(D,\Fil_\bullet,\nu)]$ dans $F_d^\square(L[\varepsilon])$ et d'image $[x_0]$ dans $F_d^\square(L)$ l'unique 
d\'eformation de $V_0$ sur $L[\varepsilon]$ dont le $\Drig$ est isomorphe \`a $D$.
Il vient que $$T_{V_0,\Ref}:={\rm Im}(df_{x_0}) \subset T_x(X)$$
est exactement le sous-espace des d\'eformations $\Ref$-triangulines de
$V_0$ au sens de~\cite[\S 2.5]{bch} et~\cite[\S 3]{chU3}. En effet,
$df_{x_0}$ se factorise \'evidemment par le morphisme d'oubli 
$$\widehat{(F_d^\square)}_{x_0}(L[\varepsilon]) \rightarrow
\widehat{(F_d)}_{[(D_0,\cT)]}(L[\varepsilon]),$$ et $V_0$ \'etant g\'en\'erique
l'alg\`ebre $\End_{\fg/L}(D_0)$ est r\'eduite aux homoth\'eties d'apr\`es~\cite[lemma 3.21]{chU3}, de sorte que 
$\widehat{(F_d)}_{[(D_0,\cT_0)]}$ est le foncteur des d\'eformations triangulines de $(D_0,\cT)$ d'apr\`es le 
lemme~\ref{sousfonct}. La relation~\eqref{inclfond}, combin\'ee au fait que $W$ est ferm\'e et $\Omega$
r\'eduit, assure donc que $$T_{V_0,\Ref} \subset T_x(W)$$ pour tout
raffinement $\Ref$ de $V_0$. On conclut alors par le r\'esultat crucial
suivant, d\'emontr\'e dans~\cite{chU3} : "Toute d\'eformation \`a l'ordre
$1$ d'une repr\'esentation cristalline g\'en\'erique est combinaison
lin\'eaire de d\'eformations triangulines".

\begin{thm} {\rm (\cite[Thm. C]{chU3}) } $T_x(X) = \sum_{\Ref}
T_{V_0,\Ref}$.
\end{thm}

\subsection{Une g\'en\'eralisation} Terminons par un \'enonc\'e plus g\'en\'eral que nous d\'emontrons 
par la m\^eme m\'ethode. Soit $d\geq 1$ un entier. Nous renvoyons \`a \cite{chdet} pour la notion de vari\'et\'e des caract\`eres 
$p$-adiques de dimension $d$ d'un groupe profini $G$. Disons simplement ici que
si l'on suppose que tout sous-groupe ouvert $U \subset G$ a la propri\'et\'e
que $\Hom(U,\F_p)| < \infty$, alors le foncteur ${\rm Aff} \rightarrow {\rm
Ens}$ associant \`a $A$ l'ensemble des pseudo-caract\`eres continus $G
\rightarrow A$ de dimension $d$ est repr\'esentable par un espace analytique
$p$-adique $\cX$ : c'est la vari\'et\'e des caract\`eres $p$-adiques de $G$ de
dimension $d$. En particulier, les points de cette vari\'et\'e sont en
bijection avec les classes de conjugaison de repr\'esentations continues semi-simples $G \rightarrow
GL_d(\Qpb)$. Dans le cas particulier $G={\rm Gal}(\Qpb/\Qp)$ on
peut d\'emontrer que $\cX$ est \'equidimensionnel de dimension $d^2+1$, et que
son lieu singulier co\"incide exactement avec le lieu des repr\'esentations
r\'eductibles (avec une exception quand $d=2$) : voir~\cite{chnp}. Des arguments similaires \`a ceux de la preuve 
du th\'eor\`eme $A$ d\'emontrent :

\begin{thm} Soit $\cX_d$ la vari\'et\'e des caract\`eres $p$-adiques de
dimension $d$ de ${\rm Gal}(\Qpb/\Qp)$. Soit $Y$ une composante
irr\'eductible de $\cX_d$ et $L$ une extension finie de $\Q_p$ telle que
$Y(L)$ contiennent un point cristallin d\'eploy\'e absolument irr\'eductible. Alors les
points cristallins d\'eploy\'es de $Y(L)$ sont Zariski-dense et d'accumulation dans $Y$.
\end{thm}

Il semble raisonnable de formuler la conjecture suivante, que l'on peut voir comme un analogue local de 
la conjecture de modularit\'e de Serre.

\begin{conj}\label{serrelocal}  Toute composante irr\'eductible de $\cX_d$ contient un point cristallin
absolument irr\'eductible. 
\end{conj}

\begin{cor} Si la conjecture \ref{serrelocal} est vraie, alors les points cristallins sont Zariski-dense dans
$\cX_d$.
\end{cor}

L'application "repr\'esentation r\'esiduelle" r\'ealise l'espace $\cX$ comme r\'eunion disjointe admissible d'espaces $\cX(\rhob)$
index\'es par l'ensemble des repr\'esentations continues
semi-simples $\rhob : G_{\Q_p} \longrightarrow \GL_n(\Fpb)$ prises modulo
isomorphisme et action du Frobenius sur l'image. Lorsque $\rhob$ est
absolument irr\'eductible, la composante $\cX(\rhob)$ est simplement l'espace
d\'efini dans l'introduction. Quand $\rhob \not \simeq \rhob(1)$ il est
irr\'eductible (c'est une boule!) et la conjecture ci-dessus est facile :
c'est la proposition~\ref{existcrisgen} (i). \`A d\'efaut de pouvoir
proposer une argument plausible pour la conjecture ci-dessus, d\'emontrons
le r\'esultat suivant. Il assure que pour tout $\rhob$, au moins une (de l'ensemble fini) des composantes
irr\'eductibles de $\cX(\rhob)$ contient un point cristallin absolument
irr\'eductible.

\begin{prop}\label{liftgeneral} Soit $\rhob : G_{\Q_p} \rightarrow \GL_d(\F_q)$ une
repr\'esentation semi-simple continue. Il existe une extension finie
$L/\Q_p$ et une repr\'esentation $G_{\Q_p} \rightarrow \GL_d(L)$ cristalline
absolument irr\'eductible dont la repr\'esentation r\'esiduelle est $\rhob$.
\end{prop}

Avant de proc\'eder \`a la d\'emonstration, il est raisonnable de
commencer par la proposition plus simple suivante, qui est 
une application d\'ej\`a frappante des r\'esultats de cet article 
`a l'existence de repr\'esentations cristallines ayant certaines
propri\'et\'es. La seconde partie est \`a
comparer avec des observations ant\'erieures de l'auteur et Bella\"iche
(voir par exemple \cite[\S 4]{bch} ou encore \cite[lemme 3.3]{bchsign}).
Dans cet \'enonc\'e, $r$ est une repr\'esentation semi-simple quelconque et
$L/\Q_p$ une extension finie de $\Q_p$.

\begin{prop}\label{famcolemandeux} Soient $x \in S_d^{\square}(r)(L)$ tel que $\delta(x) \in B_d(L)$ et $U \subset
S_d^\square(r)$ un voisinage de $x$. Pour tout r\'eel $C>0$, il existe un
ouvert affino\"ide $V \subset U$ et un ensemble Zariski-dense de $y \in V(L)$
tels que la repr\'esentation $V_y$ soit cristalline et tels que
$\delta(y)=(\delta_i) \in \cT(L)^d$ ait les propri\'et\'es suivantes :
\begin{itemize}
\item[(a)] il existe une suite strictement croissante d'entiers $k_1,\dots,k_d$ telle que
$\delta_i(\gamma)=\gamma^{-k_i}$ pour tout $i$ et telle que pour toute paire
de parties distinctes $I,J \subset \{1,\dots,d\}$ avec $1\leq |I|=|J| < d$, on ait $|\sum_{i
\in I} k_i -\sum_{j \in J} k_j| > C$.
\item[(b)] pour tout $i$, $v(\delta_i(p))=v(\delta(x)_i(p))$,
\end{itemize}
\noindent On peut de plus supposer que les $V_y$ sont absolument
irr\'eductibles dans chacun des cas suivants: 
\begin{itemize}
\item[(i)] $V_x$ est absolument irr\'eductible, \ps
\item[(ii)] $x$ a la propri\'et\'e que pour toute partie $I
\subset \{1,\dots,n\}$ telle que $1\leq |I| < n$ on ait $\sum_{i \in I}
v(\delta(x)_i(p)) \neq 0$.
\end{itemize}
\end{prop}

\begin{pf} Quitte \`a r\'etr\'ecir $U$ on peut supposer qu'il est dans
l'overt affino\"ide donn\'e par la proposition~\ref{conskliu}, en
particulier $U \subset S_d^\square(r)$. 
Consid\'erons un voisinage $U$ de $x$, ainsi que $\iota$ et
$\Omega$ comme dans le
corollaire~\ref{corfamcoleman}. Quitte \`a r\'etr\'ecir $\Omega$ on peut
supposer que les $\delta_i(p)$, vues comme fonctions analytiques sur $U$, sont
chacune de valuation constante. L'existence de $y$ satisfaisant (a) et (b) r\'esulte alors
de la proposition~\ref{critcris} et de ce que l'ensemble des suites croissantes $(k_1,\dots,k_d)$ vues comme
\'el\'ements de $\Hom(\Z_p^\times,L)^d$ qui satisfont la condition la
condition (b) (pour un $C$ fix\'e) est Zariski-dense dans
$\Hom(\Z_p^\times,\mathbb{G}_m)^d$ et s'accumule en tous les points de
$\Z^d$. \ps
Pour le second point, il d\'ecoule dans le cas (i) de la
proposition~\ref{conskliu} et du fait classique que si $\rho : G \rightarrow \GL_d(A)$
est une repr\'esentation d'un groupe $G$ \`a valeurs dans un anneau commutatif $A$
telle que pour un $x \in {\rm Spec}(A)$ la repr\'esentation $\rho_x
: G \rightarrow \GL_d(k(x))$ \'evalu\'ee en $x$ est absolument
irr\'eductible, alors il en va de m\^eme pour tout $\rho_y$ pour $y$ dans un
voisinage ouvert de $x$ dans ${\rm Spec}(A)$ : en effet, il existe $g_1,\dots g_{d^2} \in G$ tels que $\det({\rm
trace}(\rho(g_ig_j))) \in A_x^\times$ (Wedderburn) et donc dans $A_f^\times$ pour
$x \in D(f)$ et $f$ bien choisie.\ps
Dans le cas (ii), il suffit de voir que le $D_{\rm cris}(V_y)$ n'a pas de
sous-$\varphi$-module filtr\'e admissible non trivial. Mais si on a un tel
sous-module, disons de rang $1 \leq r < d$, l'\'egalit\'e des extr\'emit\'es de ses polygones de Hodge et Newton 
implique qu'il existe deux parties
$I,J \subset \{1,\dots,d\}$ avec $|I|=|J|=r$ telles que $$\sum_{i \in I}k_i
= \sum_{j \in J} v(\varphi_j)$$
o\`u $\varphi_1,\dots,\varphi_d \in L^\times$ d\'esignent les valeurs propres du Frobenius
cristallin. Mais par construction (et quitte \`a renum\'eroter les
$\varphi_i$), on a $\delta_i(p)=\varphi_i p^{-k_i}$ pour tout $i$, de sorte
que l'\'egalit\'e ci-dessus s'\'ecrive aussi $$\sum_{i \in I} v(\delta_i(p))
= \sum_{i \in I} k_i - \sum_{j \in J} k_j.$$
Mais par le (b), le terme de gauche est aussi $\sum_{i \in I}
v(\delta(x)_i(p))$, qui est un nombre fix\'e disons $C'$. Ainsi, si on choisit $y$ de sorte que le (a) est
v\'erifi\'e pour $C > C'$, il vient que $I=J$, et on obtient la contradiction voulue.
\end{pf}

\begin{pf} (de la proposition~\ref{liftgeneral}) Quitte \`a grossir $\F_q$
on peut supposer que $\rhob=\oplus_{i=1}^s \rhob_i$ o\`u chaque $\rhob_i$
est absolument irr\'eductible.  Nous allons raisonner par r\'ecurrence sur
le nombre $s$ de facteurs, le cas $s=1$ r\'esultant de la
proposition~\ref{existcrisgen} (i).  Pour $s>1$ on peut donc trouver des
repr\'esentation $L$-cristallines $V_1$ et $V_2$ de repr\'esentations
r\'esiduelles respectives $\rhob_1$ et $\oplus_{i=2}^r \rhob_i$.  On pose
$$a=\dim_L V_1, \,\,\,\,\,b=\dim_L V_2.$$ On note $k_1 \leq k_2 \leq \dots k_a$ (resp. 
$k_{a+1} \leq k_{a+2} \leq \dots \leq k_d$) les poids de Hodge-Tate de $V_1$
(resp.  $V_2$) et $\varphi_1,\dots,\varphi_a$ (resp. 
$\varphi_{a+1},\dots,\varphi_{d}$) les valeurs propres du Frobenius
cristallin de $V_1$ (resp. $V_2$).  Quitte \`a appliquer la
proposition ci-dessus (plus exactement, sa d\'emonstration) \`a $V_1$ et
$V_2$, munis de leur raffinements respectifs associ\'es \`a
$(\varphi_1,\dots,\varphi_a)$ et $(\varphi_{a+1},\dots,\varphi_d)$ on peut supposer
que $(k_i)$ est strictement croissante et que si $C={\rm Sup}_i |v(\varphi_i)-k_i|$ alors pour
toute paire de parties $I \neq J \subset \{1,\dots,d\}$ avec $|I|=|J| < d$ on ait
\begin{equation}\label{eqsomme} |\sum_{i \in I} k_i-\sum_{j \in J} k_j | > dC.\end{equation}
En particulier, les $\varphi_i$ sont distincts. Consid\'erons une extension non triviale cristalline $$ 0 \longrightarrow V_1
\longrightarrow V \longrightarrow V_2 \longrightarrow 0.$$
Il y a plusieurs fa\c{c}ons de voir qu'une telle extension existe, cela
vient par exemple de ce que $H^1_f(G_{\Qp},V_1 \otimes_L V_2^\vee)$ est de
dimension au moins \'egale au nombre de poids de Hodge-Tate strictement n\'egatifs de
$V_1 \otimes_L V_2^\vee$ par un r\'esultat classique de Bloch-Kato (ici ces poids sont les $k_i - k_j$ 
avec $i\leq a$ et $j > a$, qui sont en fait tous
strictement n\'egatifs). Consid\'erons la permutation $\sigma \in \got{S}_d$
d\'efinie par $\sigma(i)=i+a$ si $i=1,\dots,b$ et $\sigma(i+b)=i$ si
$i=1,\dots, a$, et consid\'erons le raffinement
$$(\varphi_{\sigma(1)},\varphi_{\sigma(2)},\dots,\varphi_{\sigma(d)})$$
de la repr\'esentation $V$. Soit $\tau \in \got{S}_d$ la permutation
telle que $(k_{\tau(1)},k_{\tau(2)},\dots,k_{\tau(d)})$ est la suite des sauts de la 
filtration de Hodge de $D_{\rm cris}(V)$ d\'efinie par $\Ref$
(voir~\cite[\S 2.4]{bch}). Il nous suffira ici de dire que si $(D,\cT)$ est le $\fg$-module triangulin associ\'e \`a
$(V,\Ref)$, alors ses param\`etres $\delta_i$ v\'erifient 
$$\delta_i(p)=\varphi_{\sigma(i)}p^{-k_{\tau(i)}}, \, \,
\delta_i(\gamma)=\gamma_i^{-k_{\tau(i)}} \,\, \,\,\,\forall \gamma \in
\Z_p^\times.$$
Soit $D' \subset D$ l'unique sous-$\fg$-module
cristallin satur\'e dont les valeurs propres du Frobenius
cristallin sont les $\varphi_i$ avec $a+1 \leq i \leq b$. L'application
naturelle $\eta : D' \rightarrow D_{\rm rig}(V_2)$ est injective (c'est un
isomorphisme sur les $\Dc$) et on constate que l'on a l'in\'egalit\'e pour l'ordre lexico-graphique $$(\tau(1),\tau(2),\dots,\tau(b))
\leq (a+1,a+2,\dots,d).$$
Cette in\'egalit\'e est une \'egalit\'e si et seulement si $\eta$
est surjective, ce qui ne se produit pas car $V$ est
non scind\'ee. C'est donc une in\'egalit\'e stricte, de sorte qu'il existe un $1 \leq i_0 \leq b$ tel que
$\tau(i_0) \leq a$, et en particulier $$\tau(i_0) \leq a < \sigma(i_0).$$ \ps
Nous allons appliquer une variante de la proposition pr\'ec\'edente cas (i)
et (ii) au point
$x=(D,\cT,\nu)$ pour un choix
quelconque de $\nu$. Cette proposition permettrait de conclure si nous
savions que pour toute partie $I \subset \{1,\dots,d\}$ telle que $1\leq |I| \leq d-1$
alors $\sum_{i \in I} v(\delta_i(p))  \neq 0$. Cela n'est pas toujours
satisfait dans notre cas, mais $U$ n'ayant que deux facteurs irr\'eductibles
c'est un exercice de v\'erifier que si
aucune des repr\'esentations $V_y$ donn\'ee par cette proposition n'est
absolument ir\'eductible (pour
$y$ assez proche de $x$), alors elles ont toutes un constituent de
dimension $\dim(V_1)=a$. La th\'eorie de Sen en famille, appliqu\'ee \`a la
famille donn\'ee par la proposition~\ref{conskliu}, permet d'identifier ses
poids de Hodge-Tate :  si ceux de $V_y$ rang\'es en
ordre croissant sont les $k'_i$, le (ou un si $a=b$) constituent de
dimension $a$ a pour poids de Hodge-Tate les $k'_{\tau^{-1}(i)}$ pour
$i=1,\dots,a$. Par cons\'equent, l'argument de la preuve de la
proposition ci-dessus montre qu'il suffit dans notre cas
de v\'erifier
que $\sum_{i \in I} v(\delta_i(p))  \neq 0$ quand $I=\tau^{-1}(\{1,\dots,a\})$. Mais $$|\sum_{i \in I}
v(\delta_i(p))- (\sum_{i \in \sigma(I)} k_i-\sum_{i \in \tau(I)} k_i)| < dC.$$
Par la relation~\eqref{eqsomme}, il suffit pour conclure de voir que $\sigma(I)
\neq \tau(I)$, c'est \`a dire que $\sigma \tau^{-1} (\{1,\dots,a\}) \neq
\{1,\dots,a\}$. On conclut car $\tau(i_0) \leq a$ et $\sigma(i_0) > a$. 
\end{pf}

\begin{remark} Notons qu'il ressort de la d\'emonstration que le corps de 
d\'efinition $L$ peut \^etre pr\'ecis\'e dans
l'\'enonc\'e de la proposition~\ref{liftgeneral}. Par exemple, si $\rhob$
est la repr\'esentation triviale, alors on peut choisir $L=\Q_p$.\end{remark}


\end{document}